\ifdefined\pdfoutput\pdfoutput=1\fi
\documentclass[aos,preprint]{imsart}
\ifdefined\pdfminorversion
\fi

\usepackage{amsmath,amssymb,mathtools,amsthm}
\usepackage{graphicx}
\usepackage{booktabs}
\usepackage{subcaption}
\usepackage{tabularx}
\usepackage{placeins}
\usepackage{enumitem}
\usepackage[numbers,sort&compress]{natbib}
\usepackage[colorlinks=true,linkcolor=blue,citecolor=blue,urlcolor=blue]{hyperref}
\usepackage[capitalise]{cleveref}

\startlocaldefs
\theoremstyle{plain}
\newtheorem{theorem}{Theorem}[section]
\newtheorem{lemma}[theorem]{Lemma}
\newtheorem{proposition}[theorem]{Proposition}
\newtheorem{corollary}[theorem]{Corollary}

\theoremstyle{definition}
\newtheorem{definition}[theorem]{Definition}
\newtheorem{assumption}[theorem]{Assumption}

\newtheorem{remark}[theorem]{Remark}

\newcommand{\R}{\mathbb{R}}

\newcommand{\E}{\mathbb{E}}
\newcommand{\Prob}{\mathbb{P}}
\newcommand{\dd}{\mathrm{d}}
\newcommand{\FWER}{\mathrm{FWER}}
\newcommand{\FDR}{\mathrm{FDR}}
\newcommand{\Power}{\mathrm{Power}}
\DeclareMathOperator{\supp}{supp}
\DeclareMathOperator{\spn}{span}
\DeclareMathOperator{\SO}{SO}
\DeclareMathOperator{\Cov}{Cov}
\DeclareMathOperator{\Var}{Var}
\DeclareMathOperator*{\argmax}{arg\,max}

\endlocaldefs

\begin{document}

\begin{frontmatter}
\pdfsubject{Preprint}
\title{Multiple testing for detecting rotated templates in Gaussian noise}
\runtitle{Detecting rotated templates in noise}

\begin{aug}
\author[A]{\fnms{Keren}~\snm{Mor Waknin}\ead[label=e1]{kerenmor@tauex.tau.ac.il}}
\and
\author[A]{\fnms{Yoel}~\snm{Shkolnisky}\ead[label=e2]{yoelsh@tauex.tau.ac.il}}
\address[A]{Department of Applied Mathematics, Tel-Aviv University, Israel\printead[presep={,\ }]{e1,e2}}
\runauthor{K. Mor Waknin and Y. Shkolnisky}
\end{aug}

\begin{abstract}
	We study the detection of many copies of known templates in a large noisy image. The number of copies is unknown, as are their positions, rotations and amplitudes. A copy is detected when a location within a prescribed distance of its position is reported; any other reported location is a false detection. Deciding which locations to report is a nonstandard multiple testing problem. The tested locations are selected from the image, and the test statistics are spatially dependent. Moreover, each test must cope with the unknown rotation of a copy. We therefore use a test statistic that is invariant to rotations, so no search over rotations is needed. It is obtained by projecting the patch around each location onto a finite-dimensional space of functions closed under rotations. We show how to choose this space for detecting the templates, rather than for approximating them. A location is reported when its test statistic exceeds a threshold. The threshold is obtained by adapting either the Bonferroni or the Benjamini--Hochberg procedure. We prove that the first adaptation asymptotically controls the family-wise error rate, and the second the false discovery rate (FDR). Both detect all copies with probability tending to one, provided the copies are sufficiently far apart. These guarantees hold for stationary Gaussian noise of known distribution. They require the signal-to-noise ratio of the projected copies to grow only logarithmically with the image volume. We show that this logarithmic growth cannot be improved in general by any procedure that controls the FDR, when the number of copies grows more slowly than the image volume.
\end{abstract}


\end{frontmatter}

	\section{Introduction}
	\label{sec:introduction}
	
	We study the problem of detecting many object instances embedded in a
	large noisy image on $\mathbb{R}^d$, for $d\ge 2$.
	To fix intuition, consider $d=2$ and imagine a single template
	(say, a ``duck'') appearing several times in the image, with each instance having
	an unknown location, in-plane rotation, and positive amplitude.
	Figure~\ref{fig:intro_illustration} shows the corresponding clean
	and noisy images; only the noisy image is observed.
	\begin{center}
		\begin{minipage}{0.34\linewidth}
			\centering
			\includegraphics[width=\linewidth]{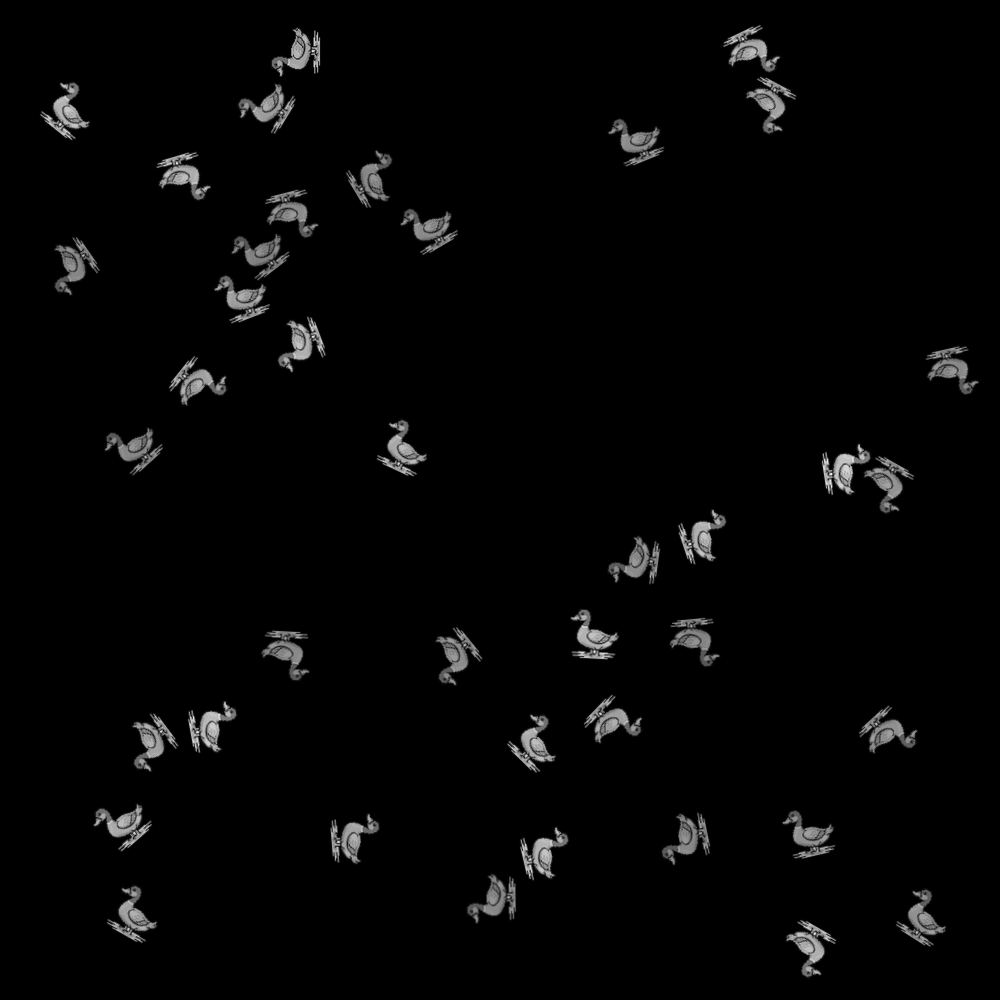}
		\end{minipage}
		\hspace{0.04\linewidth}
		\begin{minipage}{0.34\linewidth}
			\centering
			\includegraphics[width=\linewidth]{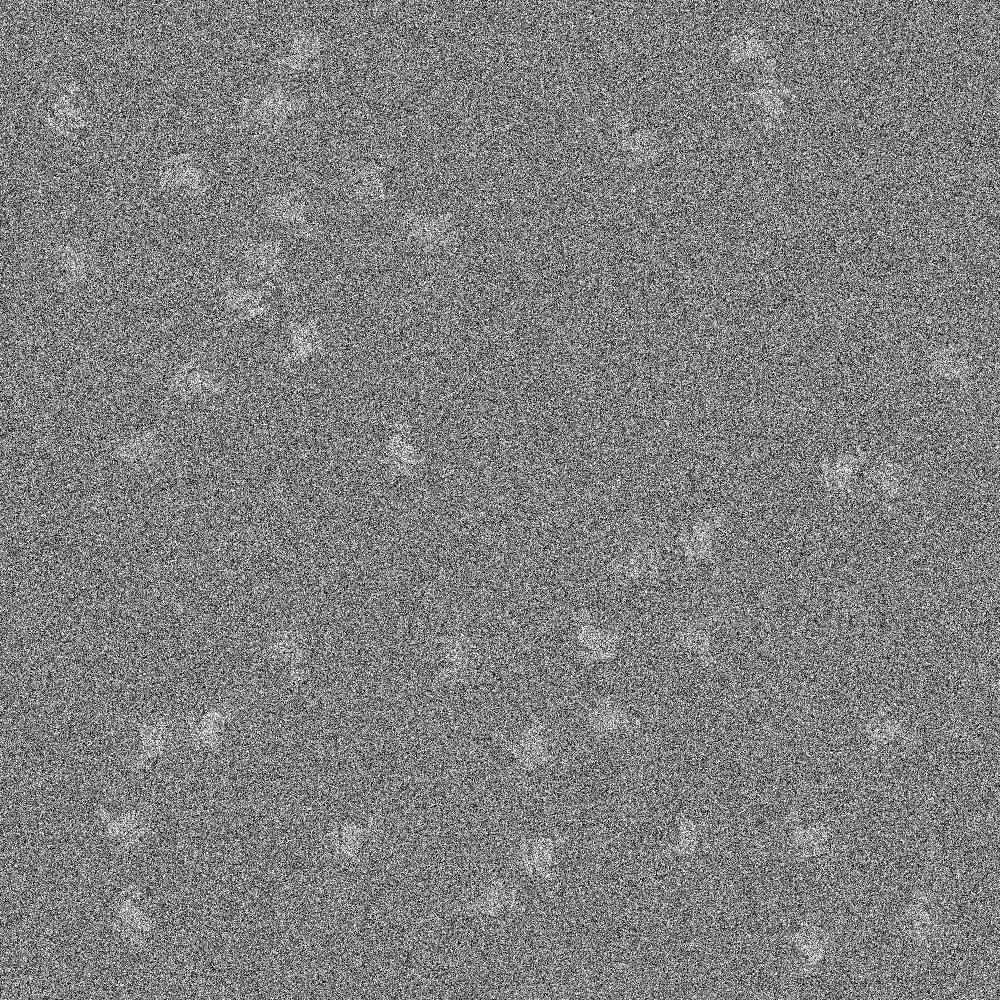}
		\end{minipage}
		\captionof{figure}{The clean image (left) contains several translated and
			rotated copies of a single template, shown here as a ``duck,'' each
			multiplied by an unknown positive amplitude. Only the noisy image
			(right) is observed.}
		\label{fig:intro_illustration}
	\end{center}
	More generally, we consider a fixed finite collection of known
	templates and a specified stationary Gaussian noise model.
	We model each object instance as a translated and rotated copy of
	one of the templates, multiplied by a positive amplitude.
	In signal processing terms, each object instance is a signal to be
	detected in noise.
	The number of instances is unknown; for each instance, its location,
	rotation, amplitude, and the template from the collection used to
	generate it are also unknown.
	The output of an object detection algorithm is a set of locations, and the statistical goal is to
	locate the object instances while controlling false detections.
	We fix a distance tolerance in advance and count a
		reported location as a true detection if it lies within the prescribed tolerance of
		the location of an object instance, and as a false detection otherwise.
	Power is the expected fraction of object instances whose
	locations are within the prescribed tolerance of at least one
	reported location.
	We aim to control either the probability of at least one false
	detection, namely the family-wise error rate
	FWER~\cite{holm1979}, or the expected proportion of false detections
	among all reported locations, namely the false discovery
	rate FDR~\cite{benjaminiHochberg1995}.
	
	A classical starting point is to ask whether a single
	template (a known signal) is present at a given location
	in the noisy observed signal.
	Matched filtering computes the correlation between the observed
	image and the known template at the given
	location.
	The correlation is compared with a threshold chosen for a
	prescribed false-alarm probability.
	In white Gaussian noise, among all tests with the prescribed
	false-alarm probability, the matched filter test maximizes the
	probability of detecting the template when it is present~\cite{kay1998detection,vanTreesBellTian2013}.
When the template is unknown but lies in a given linear subspace,
	Scharf and Friedlander~\cite{scharfFriedlander1994} analyze a
	detector called the matched subspace detector, which uses the
	energy (squared norm) of the data projected onto the subspace.
	In white Gaussian noise with known variance, they derive explicit
	expressions for its false-alarm and detection probabilities.
	
	Assuming a single known  template is present, we can ask how
	accurately its location can be estimated.
	For one-dimensional location estimation in white Gaussian noise,
	the Ziv-Zakai bound~\cite{zivZakai1969} gives a lower bound on
	the mean squared error averaged over possible locations.
	Zakai and Ziv~\cite{zakaiZiv1969threshold} and
	Merhav~\cite{merhav2011threshold} study how the probability of
	large location errors depends on the signal-to-noise ratio.
For orientation estimation in two-dimensional data,
	P\"usp\"oki et al.~\cite{puspokiEtAl2019} study limits on the
	accuracy of estimating the rotation angle of a known
	two-dimensional template at a given location.
	
	When the number of object instances in an image is unknown, the task is to locate all of them.
	Template matching computes the correlation of the image with the
	template at every location within the image, with high values indicating likely object
	instances.
	When rotations of the object instances are unknown, the score at each location is the
	largest correlation over rotated copies of the template.
	Almira et al.~\cite{almiraEtAl2024} present a computationally
	efficient method for template matching over rotations, avoiding
	a dense search over angles.
	The notion that enables fast rotation search is steerability.
	A template is steerable if all its rotations lie in a fixed
	finite-dimensional space. This allows correlations with rotated
	copies of the template to be computed without repeating the full
	calculation for each angle~\cite{freemanAdelson1991}.
	For template matching, Zhao and Blu~\cite{zhaoBlu2020} construct
	steerable approximations of a given template.
	The background noise can also be taken into account when designing the detector: Fageot et al.~\cite{fageot2021steerability} use steerability to design a correlation-based detector for a single known template by maximizing a signal-to-noise criterion under a model that allows spatially correlated noise. Another approach is to train models on annotated images to locate object instances.
	The methods in~\cite{ren2015faster,redmon2016yolo,carion2020detr}
	predict the locations and classes of object instances and assign each prediction
	a score, but do not provide guarantees of FWER or FDR control
	for the reported detections.
	
	In template matching, the correlation with the template
		gives a score at each location. A~statistical question is how to
		choose a threshold on this score to control false detections.
	When the question is whether a known template appears
	anywhere in the image, we find the location in the image
		with the largest correlation and compare that correlation with a
	threshold calibrated from the distribution of the maximum
	of the correlation under noise alone.
	Siegmund and Worsley~\cite{siegmundWorsley1995} approximate the distribution of the maximum of the noise-only correlation when the size of the template is also
	unknown, and later work extends the approximation to searches over rotations~\cite{shafieSigalSiegmundWorsley2003,adlerSubagTaylor2012};
	see~\cite{adlerTaylor2007,piterbarg1996} for the theory of maxima of random fields.
	
	When there are multiple occurrences of the template in the image, we need to first select a set of candidate locations, and then test for the presence of the template at each candidate location. This leads to a multiple testing problem
		in which the goal is to decide which locations to report while
		controlling FWER or FDR.
	Schwartzman, Gavrilov, and Adler~\cite{schwartzmanGavrilovAdler2011}
	study the detection of one-dimensional signals in stationary
	ergodic Gaussian noise.
	After smoothing the observed noisy signal by convolution, they compute the p-values of the detected local maxima 
	using the distribution of local-maximum heights under noise alone,
	and apply Bonferroni or Benjamini--Hochberg (BH) procedures to
	decide which maxima to report.
	Under suitable conditions, Bonferroni controls FWER and BH controls FDR asymptotically, with power tending to one for both procedures.
	Cheng and Schwartzman~\cite{chengSchwartzman2017} extend the
	local-maxima testing approach to higher-dimensional Gaussian random fields, with asymptotic FDR control and power tending to one. In both works, a detection is counted as true if it falls anywhere within the
	support of the object instance, rather than within a prescribed distance of its
	center. For object instances modeled as arbitrary elements of a given linear subspace,
	Mor Waknin et al.~\cite{wakninEtAl2026object} develop a detection
	procedure that, under their model assumptions, controls FWER
	or FDR asymptotically and detects every instance with probability
	tending to one.
	
In this work, we propose an object detection scheme
with statistical guarantees for detecting translated
and rotated copies of known templates in a noisy image. Our detection scheme uses a fixed finite-dimensional space of functions, derived from the templates,
		which we call the detection subspace. This subspace is steerable, meaning that it contains every
rotation of each of its functions.
	At each location in the image, we compute a score by summing the
	squared correlations of the image with the functions of an orthonormal
	basis of the detection subspace.
		By rotation invariance, the score at any location stays the same
		if we rotate the image patch about that location, so no search
over rotation angles is needed.
	We calibrate detection thresholds from the distribution of the
	maximum score over a region of fixed size under noise alone.
	Using the maximum score distribution for threshold calibration
	was proposed in~\cite{wakninEtAl2026object} but without rigorous
	performance proofs. Here we also give the missing proofs.
	Furthermore, for a fixed detection subspace, we derive
	nonasymptotic bounds on detection probabilities and error rates
	and give conditions for locating all or almost all object instances
	within the prescribed distance tolerance, with asymptotic
	FWER control using Bonferroni or FDR control using BH.
These conditions are stated in terms of the ratio of the
	smallest score of an object instance to the expected
	maximum noise score over the fixed region used to calibrate the
	detection thresholds. We show how large this ratio must be for our
	bounds to guarantee that all object instances are located with
	probability tending to one. For BH, we give a weaker condition on
	the same ratio under which the expected fraction of object instances
	located (power) tends to one. When the number of object instances grows more slowly
than the image volume, the logarithmic factor in our
bounds for BH power cannot be improved in general
under asymptotic FDR control
(Proposition~\ref{prop:retained_energy_lower_bounds}).
	
	    Our detection and error-control guarantees hold for any fixed 
		detection subspace. We use these guarantees to derive a criterion for
		choosing the detection subspace from the known templates and the noise
		model, so that the detection subspace is chosen to improve
        detection rather than to approximate the templates
        accurately.
        All simulations were implemented in Python. The code is available at
        \url{https://github.com/ShkolniskyLab/steerable_detection}. 

	\paragraph{Organization.}
	Section~\ref{sec:model_detection} presents the observation model and detection
	procedure. Section~\ref{sec:theoretical_foundations} develops the guarantees
	and lower bounds for a fixed detection subspace.
	Section~\ref{sec:repr_framework} constructs steerable subspaces from the known templates and develops a criterion for
	selecting the detection subspace.
	Section~\ref{sec:simulations} reports the numerical experiments, and
	Section~\ref{sec:discussion} presents the main conclusions and directions for
	future work.
	
	\section{Model and detection procedure}
	\label{sec:model_detection}
	
	\subsection{Observation model}\label{subsec:background_setting}
\begin{figure}[t]
		\centering
		\begin{minipage}[t]{0.48\linewidth}
			\vspace{0pt}
			\centering
			\includegraphics[height=6.2cm]{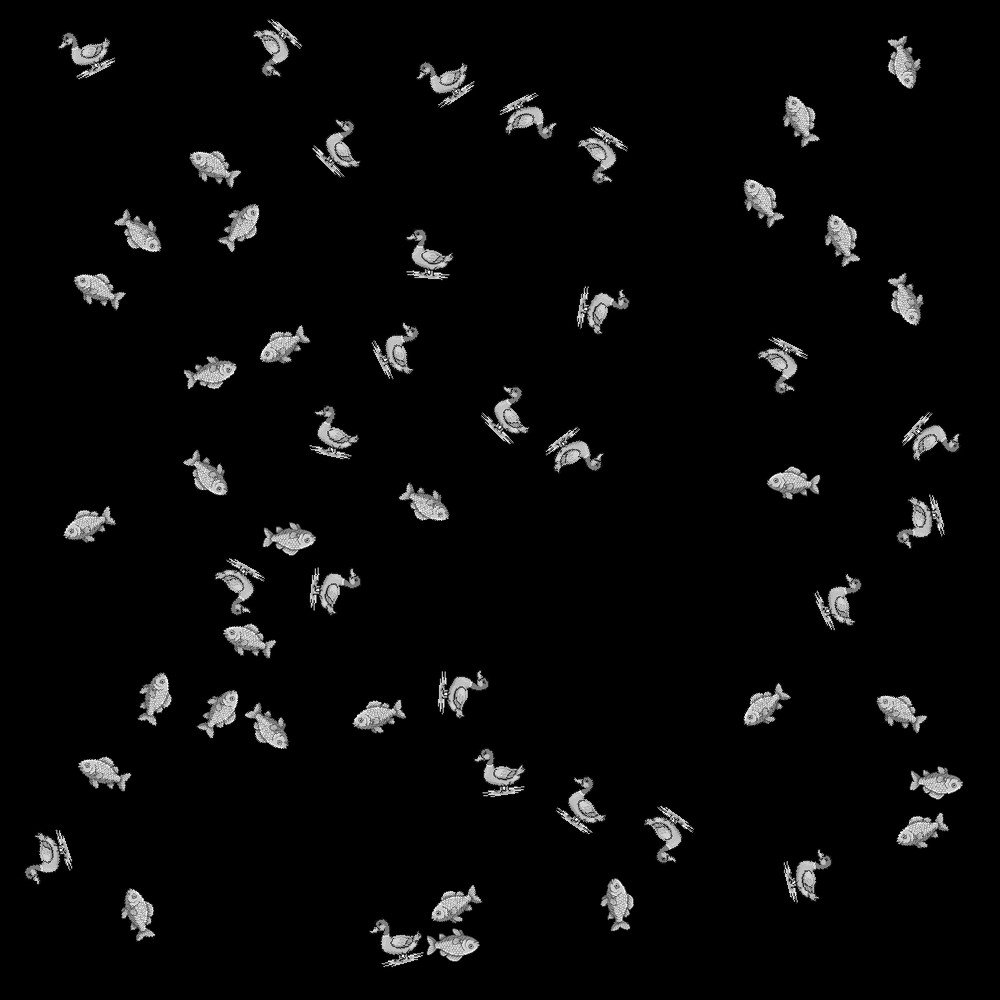}\\[2mm]
			\(x\)
		\end{minipage}\hfill
		\begin{minipage}[t]{0.48\linewidth}
			\vspace{0pt}
			\centering
			\includegraphics[height=6.2cm]{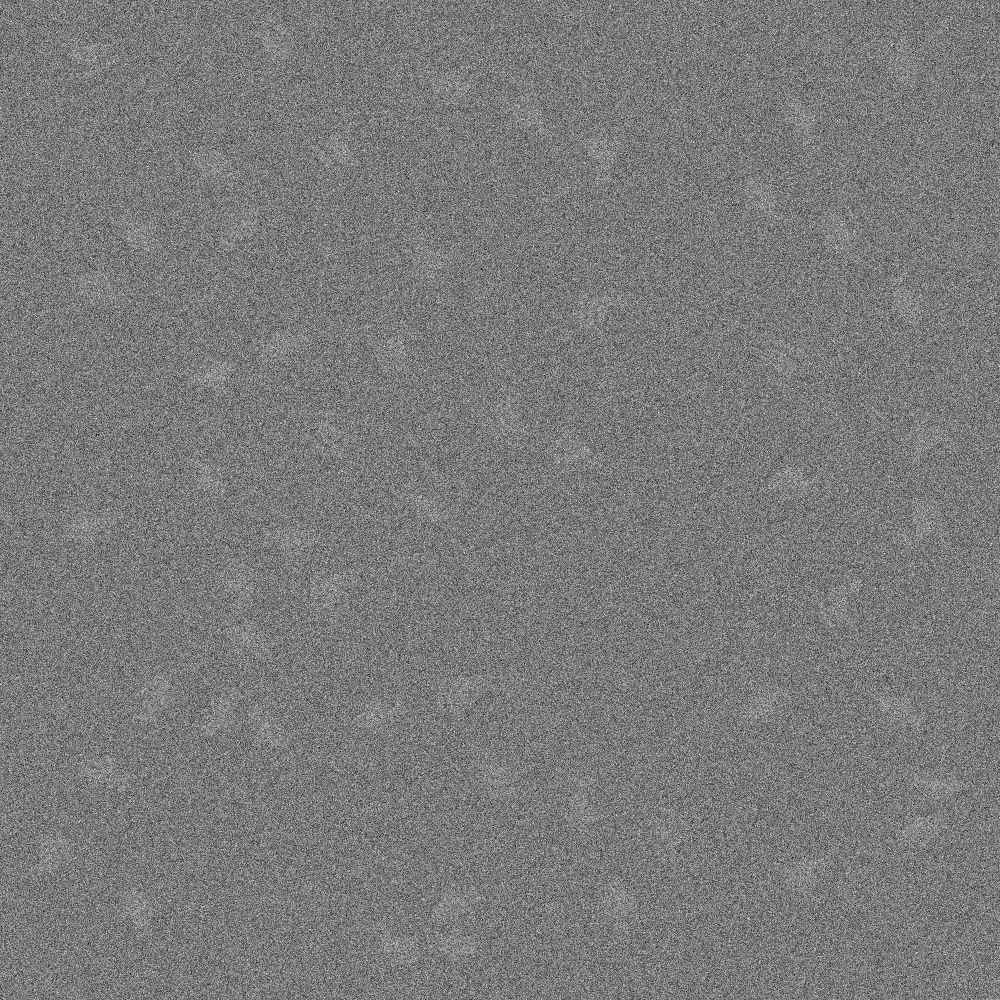}\\[2mm]
			\(y\)
		\end{minipage}
		\caption{Two-dimensional illustration of the model in
			\eqref{eq:template_instance_model}. The clean image \(x\), shown on the
			left, contains multiple copies of two known templates, a duck and a
			fish, with unknown locations, rotations, and amplitudes. The
			observed noisy image \(y\) is shown on the right.}
		\label{fig:toy_setting}
	\end{figure}

	Fix \(d\ge2\).
	Let \(y=\{y(t):t\in\mathbb R^d\}\) be a random field given by
	\begin{equation}\label{eq:the_model}
		y(t)=x(t)+z(t),
	\end{equation}
	where \(x\) is the clean image formed by the object instances to be detected,
	and \(z\) is a jointly measurable centered stationary Gaussian noise
	field on \(\mathbb R^d\). For \(v=(v_1,\dots,v_d)\in\mathbb R^d\), define
	\(\|v\|_\infty:=\max_{1\le a\le d}|v_a|\). For
	$t_0\in\mathbb{R}^d$ and $L>0$, define the open
	\(\ell_\infty\)-cube of side length $L$ centered at $t_0$, 
	\[
	\mathcal{B}(t_0,L):=\Big\{t\in\mathbb{R}^d:\ \|t-t_0\|_\infty< \tfrac{L}{2}\Big\},
	\]
	and let
	\begin{equation}\label{eq:closed_cube}
		\overline{\mathcal B}(t_0,L)
		:=
		\Big\{
		t\in\mathbb{R}^d:
		\|t-t_0\|_\infty\le \tfrac{L}{2}
		\Big\}
	\end{equation}
	denote the corresponding closed cube.
	We are given a collection of \(n\) known real-valued templates
	\(\{u_\ell\}_{\ell=1}^n\subset L^2(\mathbb R^d)\), which we normalize so that
	\[
	\|u_\ell\|_{L^2(\mathbb R^d)}=1,
	\qquad \ell=1,\dots,n.
	\]
	Let \(B>0\) be such that
	every \(u_\ell\) is supported in the closed ball
	\[
	D_{B/2}:=\{t\in\mathbb R^d:\ \|t\|\le B/2\},
	\]
	where \(\|\cdot\|\) denotes the Euclidean norm.
	For each \(L>0\), we use \(\overline{\mathcal B}(0,L)\) as the search
	domain. To apply the detection procedure throughout this domain,
	including near its boundary, we assume that \(y\) is observed on the larger
	cube $\overline{\mathcal B}(0,L+B)$.
	To model rotations, let
	\[
	\SO(d):=\{R\in\mathbb{R}^{d\times d}:\ R^\top R=I,\ \det(R)=1\},
	\]
	and let \(\SO(d)\) act on functions by \((Rf)(t):=f(R^{-1}t)\).
	We model the clean image~$x$ in~\eqref{eq:the_model} as a sum of object instances \(x_i\), namely,
	\begin{equation}\label{eq:template_instance_model}
		x(t)=\sum_{i=1}^{N} x_i(t),\qquad
		x_i(t)= b_i\,(R_i u_{\ell_i})(t-\tau_i),
	\end{equation}
	where \(1\le N=N(L)<\infty\) is the number of object instances,
	\(\tau_i\in\overline{\mathcal B}(0,L)\) is the location of the object instance,
	\(R_i\in\SO(d)\) is its rotation,
	\(\ell_i\in\{1,\dots,n\}\) is its template index, and
	\(b_i>0\) is its amplitude. Since rotations and translations preserve the
	\(L^2(\mathbb R^d)\) norm and the templates are normalized,
	\begin{equation}\label{eq:bi}
		\|x_i\|_{L^2(\mathbb R^d)}=b_i.
	\end{equation}
	Since each template is supported in \(D_{B/2}\), and every rotation maps this
	ball onto itself, \(x_i\) is supported in \(\tau_i+D_{B/2}\).
	For each \(L\), the number \(N\) and the object instance parameters
	\(\{\tau_i,R_i,\ell_i,b_i\}_{i=1}^N\) are deterministic but unknown and may vary with \(L\).
	The spatial dimension \(d\), the template collection
	\(\{u_\ell\}_{\ell=1}^n\), the support bound \(B\), and the law of \(z\)
	remain fixed as \(L\to\infty\).
	See Figure~\ref{fig:toy_setting} for an illustration.
	
	Our goal is to recover the unknown locations
	\(\{\tau_i\}_{i=1}^N\subset \overline{\mathcal B}(0,L)\) from the observation of \(y\) on
	\(\overline{\mathcal B}(0,L+B)\), given the template collection
	\(\{u_\ell\}_{\ell=1}^n\), the common support bound \(B\), and the noise law (the distribution of \(z\)).
	Table~\ref{tab:parameters_summary} summarizes the roles of the main quantities
	in the model.
	
	\begin{table}[t!]
		\centering
		\small
		\setlength{\tabcolsep}{5pt}
		\renewcommand{\arraystretch}{1.15}
		\begin{tabularx}{\linewidth}{|c|X|c|}
			\hline
			\textbf{Quantity}
			&
			\textbf{Role}
			&
			\textbf{Dependence on \(L\)}
			\\ \hline
			\(L\)
			&
			Known side length of the search domain
			&
			Varies
			\\ \hline
			\(d,\ B,\ n,\ \{u_\ell\}_{\ell=1}^n,\ \text{law of }z\)
			&
			Known model inputs
			&
			Fixed
			\\ \hline
			\(N,\ \{\tau_i,R_i,\ell_i,b_i\}_{i=1}^N\)
			&
			Unknown object instance quantities
			&
			May depend on \(L\)
			\\ \hline
		\end{tabularx}
		\caption{Roles of the main quantities in the model.}
		\label{tab:parameters_summary}
	\end{table}

	\subsection{Detection procedure}
	\label{subsec:algorithm_outline}
	
	We now present the detection procedure for the object instances \(x_i\) in~\eqref{eq:template_instance_model}, given an observation~\(y\) satisfying~\eqref{eq:the_model}. The procedure computes a score map, uses it to select
	candidate locations separated by at least a fixed distance, and
	determines which candidate locations to report using
	either Bonferroni or BH thresholding.
	We start by fixing \(0<\delta\le 2B\), which specifies
	a distance tolerance of \(\delta/2\) in the \(\ell_\infty\) norm:
	a reported location is counted as a true detection if it lies in the cube
	\(\mathcal B(\tau_i,\delta)\) around some object instance location \(\tau_i\);
	see Subsection~\ref{subsec:performance_criteria}. We call
	\(\mathcal B(\tau_i,\delta)\) the detection cube of instance \(i\),
	call \(\delta\) the localization parameter, and set
	\[
	r:=2B+\delta.
	\]
	Both \(\delta\) and \(r\) remain fixed as \(L\to\infty\) and are chosen
	independently of the observed image and the unknown parameters of the object instances described in Table~\ref{tab:parameters_summary}.
	Let \(M\ge1\), and let
	\(\{\psi_j\}_{j=1}^M\subset L^2(\mathbb R^d)\) be a fixed orthonormal
	family of real-valued functions with
	\(\supp(\psi_j)\subseteq D_{B/2}\). Define the detection subspace by
	\begin{equation}\label{eq:detection_subspace_def}
		\mathcal V:=\spn\{\psi_1,\dots,\psi_M\}.
	\end{equation}
In this subsection, the orthonormal family and its span \(\mathcal V\) are
treated as fixed and given, and the choice of \(\mathcal V\)
is addressed in Section~\ref{sec:repr_framework}.
The detection procedure consists of the following steps: 
	\begin{enumerate}
		\item \textbf{Compute the score map.}
		Define the observed score map
		\begin{equation}\label{eq:smoothed_process_intro}
			S^y(t) = \sum_{j=1}^{M}\big(y \ast\tilde{\psi}_j \big)^2(t),
		\end{equation}
		and the noise score map
		\begin{equation}\label{eq:noise_score_map}
			S^z(t) = \sum_{j=1}^{M}\big(z \ast\tilde{\psi}_j \big)^2(t),
		\end{equation}
		where \(\tilde{\psi}_j(t)=\psi_j(-t)\) and \(\ast\) denotes convolution
		on \(\R^d\). 
		\item \textbf{Select candidate locations.}\label{item:candidate_peaks}
		Using the value \(r=2B+\delta\) fixed above, proceed as follows.
		First choose a maximizer
		\[
		t_1\in\argmax_{t\in \overline{\mathcal B}(0,L)} S^y(t).
		\]
		After selecting \(t_1,\ldots,t_{k-1}\), remove the cubes
			\(\mathcal B(t_j,r)\), \(1\le j<k\), from
			\(\overline{\mathcal B}(0,L)\). If the remaining set
        \(\overline{\mathcal B}(0,L)\setminus
        \bigcup_{j=1}^{k-1}\mathcal B(t_j,r)\)
        is empty, stop.
			Otherwise, choose a maximizer
			
		\[
		t_k\in\argmax_{t\in\overline{\mathcal B}(0,L)
			\setminus\bigcup_{j=1}^{k-1}\mathcal B(t_j,r)} S^y(t).
		\]
		This produces the candidate set $T=\{t_1,\dots,t_{m_L}\}$, where
		$m_L$ is the (random) number of iterations.
		\item \textbf{Compute test values.}\label{step:compute_test_values}
		Let
		\(\overline{\mathcal B}(0,r/2)\) be the closed \(\ell_\infty\)-cube of
		side length \(r/2\) centered at the origin, which we call the
		comparison cube, and define
		\begin{equation}\label{eq:test_value_intro}
			p(u)=\Prob\!\left[
			\max_{t\in\overline{\mathcal B}(0,r/2)} S^z(t)>u
			\right].
		\end{equation}
		For each \(t\in T\), we use \(p(S^y(t))\) as a test value. 
		
		\item \textbf{Apply multiple testing.}\label{step:multiple_testing}
		Fix a target error level \(0<\alpha<1\).
		Apply a multiple-testing procedure to
		\(\{p(S^y(t)):\ t\in T\}\) as follows.
		Divide \(\mathbb R^d\) into   cubes \(Q_m\) of side length \(r/2\), where
		\begin{equation}\label{eq:Qm}
			Q_m:=\frac r2m+[0,r/2)^d,\qquad m\in\mathbb Z^d.
		\end{equation}
		Let
		\begin{equation}\label{eq:ML_def}
			\mathcal I_L
			:=\{m\in\mathbb Z^d:
			Q_m\cap\overline{\mathcal B}(0,L)\neq\emptyset\},
			\qquad
			M_L:=|\mathcal I_L|.
		\end{equation}
		Thus \(M_L\) counts the cubes \(Q_m\) that intersect
			the search domain \(\overline{\mathcal B}(0,L)\).
		The Bonferroni threshold is
		\begin{equation}\label{eq:u_bon_def}
			u_{\rm Bon}
			:=
			\inf\left\{u\in\mathbb R:\ p(u)\le \frac{\alpha}{M_L}\right\}.
		\end{equation}
		The BH threshold is obtained by applying the Benjamini--Hochberg
		procedure~\cite{benjaminiHochberg1995}, with \(M_L\) replacing
		the number of tests \(|T|\) in the denominator. Specifically, let
		\(p_{(1)}\le\cdots\le p_{(|T|)}\) denote the ordered values of
		\(\{p(S^y(t)):\ t\in T\}\), and define
		\begin{equation}\label{eq:k_BH_def}
			k
			=
			\max\left\{
			i:\ p_{(i)}\le \frac{i\alpha}{M_L}
			\right\},
		\end{equation}
		with the convention that \(k=0\) if the set is empty. If \(k=0\), set
		\(u_{\rm BH}=+\infty\) and report no detections. If \(k\ge1\), define
		\begin{equation}\label{eq:u_bh_def}
			u_{\rm BH}
			:=
			\inf\left\{u\in\mathbb R:\ p(u)\le \frac{k\alpha}{M_L}\right\}.
		\end{equation}
		
		For any threshold \(u\in\mathbb R\cup\{+\infty\}\), with the convention
		\(T(+\infty)=\emptyset\), the set of reported locations is
		\[
		T(u)=\{t\in T:\ S^y(t)\ge u\}.
		\]
	\end{enumerate}
	The intuition behind our detection procedure is as follows.
	The score \(S^y\) is designed to take large values near object
	instances. A single object instance may produce large scores at several
	nearby locations. To avoid repeatedly selecting locations corresponding to the same
	object instance, at each iteration of Step~\ref{item:candidate_peaks}
	we select a location \(t_j\) with the largest score in the remaining
	search domain and then remove the cube \(\mathcal B(t_j,r)\)
	from that domain. Consequently, distinct candidate
	locations satisfy
	\[
	\|t_j-t_k\|_\infty\ge\frac r2,
	\qquad j\ne k.
	\]
	Each candidate location lies in some \(Q_m\) with
		\(m\in\mathcal I_L\) by \eqref{eq:ML_def}. Since any two points
		in the same \(Q_m\) have \(\ell_\infty\)-distance less than \(r/2\),
		distinct candidates lie in distinct cubes.
	Thus Step~\ref{item:candidate_peaks} stops after at most \(M_L\) selections,
	and
	\[
	m_L\le M_L.
	\]
	The closure of each \(Q_m\) is a translate of the comparison
		cube $\overline{\mathcal B}(0,r/2)$ used in Step~\ref{step:compute_test_values}. Since \(S^z\) is
		stationary (since \(z\) is assumed to be), the probability that the maximum noise score
		over \(\overline Q_m\) exceeds \(u\) is therefore \(p(u)\).
		The Bonferroni and BH procedures use the resulting test values
		\(p(S^y(t))\) at the candidate locations.
	
	\section{Theoretical guarantees for a fixed detection subspace}
	\label{sec:theoretical_foundations}
	
	The detection procedure in Subsection~\ref{subsec:algorithm_outline}
	uses the orthonormal family \(\{\psi_j\}_{j=1}^M\), whose span is the
	detection subspace \(\mathcal V\) defined in~\eqref{eq:detection_subspace_def}.
	We choose this family before observing the image and keep it fixed
	as \(L\to\infty\). Section~\ref{sec:repr_framework} addresses the
	choice of the detection subspace using the known templates and noise law.
	
	We first define the statistical performance criteria for our algorithm. We then study
	the distribution of the maximum noise score over the comparison cube.
	The tail probability \(p(u)\) of this maximum, defined
	in~\eqref{eq:test_value_intro}, is used to calibrate the Bonferroni
	and BH thresholds. We also give conditions under which, in the absence
	of noise, the first \(N\) candidate locations selected in
	Step~\ref{item:candidate_peaks} include exactly one point in each
	detection cube \(\mathcal B(\tau_i,\delta)\).
	Then, we combine these results to obtain bounds for each \(L\) and
	give sufficient conditions on the smallest score of an
		object instance for power, error control, and the
	probability of reporting a location in every detection cube
	as \(L\to\infty\). Finally, Proposition~\ref{prop:retained_energy_lower_bounds}
		shows that the growth rates that these conditions impose on the
		smallest score cannot be improved in general.
	\subsection{Performance criteria}
	\label{subsec:performance_criteria}
	
	We now define the performance criteria for a
	detection procedure (with the orthonormal family
		\(\{\psi_j\}_{j=1}^M\) fixed). Recall the value
	\(0<\delta\le2B\) fixed at the beginning of
	Subsection~\ref{subsec:algorithm_outline}. A reported location is counted as
	a true detection if it lies in \(\mathcal B(\tau_i,\delta)\) for some object instance
	location~\(\tau_i\) in~\eqref{eq:template_instance_model}, and as a false
	detection otherwise. Define the true detection region and the false detection
	region by
	\begin{equation}\label{eq:true_noise_regions}
		\mathbb B_1^\delta
		:=
		\bigcup_{i=1}^N \mathcal B(\tau_i,\delta),
		\qquad
		\mathbb B_0^\delta
		:=
		\overline{\mathcal B}(0,L)\setminus \mathbb B_1^\delta .
	\end{equation}
	For \(u\in\mathbb R\cup\{+\infty\}\), recall that \(T(u)\) is the set of
	reported locations defined in Step~\ref{step:multiple_testing}.
	The numbers of reported true and false detections are
	\[
	W(u)
	:=
	\#\bigl(T(u)\cap \mathbb B_1^\delta\bigr),
	\qquad
	V(u)
	:=
	\#\bigl(T(u)\cap \mathbb B_0^\delta\bigr).
	\]
	All probabilities and expectations in the following performance criteria are
	taken with the clean image \(x\) and the object instance parameters
	(see Table~\ref{tab:parameters_summary}) held fixed. The family-wise error rate is
	\[
	\FWER(u)
	:=
	\Prob\bigl(V(u)\ge 1\bigr),
	\]
	which is the probability of making at least one false detection. The false
	discovery rate is
	\begin{equation}\label{def:fdr}
		\FDR(u)
		:=
		\E\!\left[
		\frac{V(u)}{\max\{W(u)+V(u),1\}}
		\right].
	\end{equation}
The	power of the detection procedure is
	\begin{equation}\label{def:power}
		\Power(u)
		:=
		\frac{1}{N}
		\sum_{i=1}^N
		\Prob\!\left(
		T(u)\cap \mathcal B(\tau_i,\delta)\neq\emptyset
		\right),
	\end{equation}
	which is the expected proportion of object instances for which the procedure
	reports at least one location in \(\mathcal B(\tau_i,\delta)\). The expected proportion of missed object instances is thus
	\(1-\Power(u)\). As a measure of false detections, we use \(\FDR(u)\),
		which is the expected proportion of false detections among the
	reported detections.
	Finally, we say that simultaneous recovery occurs at threshold \(u\) if \[
	T(u)\cap\mathcal B(\tau_i,\delta)\neq\emptyset
	\qquad
	\text{for every }i=1,\dots,N.
	\]
	We call the probability of this event the simultaneous recovery probability. 
	\subsection{Noise maxima}
	\label{subsec:atomlessness}
	Noise maxima are relevant both to the selection of candidate
	locations in Step~\ref{item:candidate_peaks} and to the
	thresholds used to report detections.
	The Bonferroni and BH thresholds in
	Step~\ref{step:multiple_testing} use the distribution of the
	maximum noise score over the comparison cube. Denote this maximum by
	\[
	Z_{\max}
	:=
	\max_{t\in\overline{\mathcal B}(0,r/2)}S^z(t),
	\]
	so that \(p(u)=\Prob(Z_{\max}>u)\) by \eqref{eq:test_value_intro}.
	For a prescribed probability level \(c\in(0,1)\), define
	the corresponding threshold
	\begin{equation}\label{eq:uc_def_inline}
		u_c:=\inf\{u\in\R:\ p(u)\le c\}.
	\end{equation}
	
	For a candidate location \(t\in T\), the reporting rule in
		Step~\ref{step:multiple_testing} is
	\[
	t\in T(u_c)\quad\Longleftrightarrow\quad
	S^y(t)\ge u_c.
	\]
To calibrate the reporting threshold \(u_c\), we ask how often noise
		alone produces a score of at least \(u_c\) somewhere in the comparison
		cube.
		When the image contains only noise (i.e., \(y=z\)),
		\(S^y=S^z\), and this probability is
	\[
	\Prob(Z_{\max}\ge u_c)
		= \Prob(Z_{\max}>u_c)+\Prob(Z_{\max}=u_c)
		= p(u_c)+\Prob(Z_{\max}=u_c).
	\]
	A jump of \(p\) at \(u_c\) would mean that
		\(\Prob(Z_{\max}=u_c)>0\). Hence the probability that \(Z_{\max}\)
		reaches or exceeds the threshold \(u_c\) would be larger than \(p(u_c)\)
		and could exceed \(c\) even when \(p(u_c)\le c\). Thus,
	\(	p(u_c)\le c\) does not by itself imply
		\(\Prob(Z_{\max}\ge u_c)\le c.\)
	Under the following nondegeneracy condition,
Theorem~\ref{thm:continuity_test} shows that
\(Z_{\max}\) has an absolutely continuous distribution. Theorem~\ref{thm:continuity_test} shows that \(Z_{\max}\) has
		an absolutely continuous distribution. In particular,
		\(\Prob(Z_{\max}=u_c)=0\).
	
	\begin{assumption}[Nondegeneracy of the projected noise]
		\label{ass:noise_reg_main}
		Let \(z\) be a jointly measurable centered stationary Gaussian
		field on \(\R^d\), and let \(\psi_1,\dots,\psi_M\) from Subsection~\ref{subsec:algorithm_outline} be the fixed
		orthonormal family in \(L^2(\R^d)\), supported in \(D_{B/2}\).
	For vectors in \(\mathbb R^M\), let \(\|\cdot\|\)
		denote the Euclidean norm. Define
		\[
		X_j(t):=(z*\tilde\psi_j)(t),
		\qquad
		\mathbf X(t):=(X_1(t),\dots,X_M(t))^\top,
		\qquad
		S^z(t):=\|\mathbf X(t)\|^2.
		\]
We call \(\mathbf X\) the projected noise field. Assume that
		\[
		\Sigma_{\mathbf X}
		:=
		\Cov\bigl(\mathbf X(0)\bigr)
		\]
		is positive definite, that is,
			\(v^\top\Sigma_{\mathbf X}v>0\) for every nonzero
			\(v\in\mathbb R^M\).
	\end{assumption}
	For a nonzero \(v\in\R^M\), let \(f=\sum_{j=1}^M v_j\psi_j\).
		The variance of the correlation of the noise with~\(f\) at the
		origin is
	\[
		\Var\!\left(\int_{\R^d} z(x)f(x)\,\dd x\right)
		= \Var\!\left(\sum_{j=1}^M v_jX_j(0)\right)
		= \Var\!\left(v^\top\mathbf X(0)\right)
		= v^\top\Sigma_{\mathbf X}v.
	\]
	By stationarity, this variance is the same at every
		location. Thus, Assumption~\ref{ass:noise_reg_main} rules out
a nonzero function in the detection subspace whose
correlation with the noise is zero with probability one. This assumption holds for the noise model used in
the simulations of Section~\ref{sec:simulations}:
for \(d=2\),
\(\Cov(z(t),z(s))=\sigma^2\,10^{-3\|t-s\|}\)
with \(\sigma>0\)
(Proposition~\ref{S-prop:noise_assumptions_isotropic_exponential}).
	
	\begin{theorem}[Absolute continuity of \(Z_{\max}\)]
		\label{thm:continuity_test}
		Under Assumption~\ref{ass:noise_reg_main}, the distribution of
		\(Z_{\max}\) is absolutely continuous with respect to Lebesgue
		measure. Consequently, \(\Prob(Z_{\max}=u)=0\) for every
		\(u\in\R\), and the tail function \(p\) of \(Z_{\max}\)
		is continuous on \(\R\).
	\end{theorem}
	
The proof is given in
Section~\ref{S-app:continuity_noise_test} of the supplementary material.
Corollary~\ref{S-cor:generalized_inverse_attainment} shows that, for
every \(c\in(0,1)\), the threshold \(u_c\)
defined in~\eqref{eq:uc_def_inline} satisfies
\(0<u_c<\infty\) and
\[
\Prob(Z_{\max}\ge u_c)=p(u_c)=c.
\]
The supplementary material proves the corresponding absolute continuity
result for a broader class of scores; see
Corollary~\ref{S-cor:seminorm_score_continuity}.

	At each iteration of Step~\ref{item:candidate_peaks}, we select
	a location with the largest observed score in the remaining domain.
	To analyze how noise affects the selected locations, we
	bound the largest noise score over the entire search domain
	\(\overline{\mathcal B}(0,L)\).
	Under Assumption~\ref{ass:noise_reg_main},
	Lemma~\ref{S-lem:noise_maxima_fixed_growing} in the supplementary
	material gives, for every \(L>0\),
	\begin{equation}\label{eq:search_domain_noise_bound}
	\Prob\!\left\{
	\max_{t\in\overline{\mathcal B}(0,L)}S^z(t)
	>
	\eta_{\mathrm{peak}}\bigl(1+2\sqrt{\log M_L}\bigr)^2
	\right\}
	\le \frac{1}{M_L},
\end{equation}
	where
	\begin{equation}\label{eq:eta_peak_def}
		\eta_{\mathrm{peak}}
		:=
		\E[Z_{\max}]
		\in(0,\infty)
	\end{equation}
	is the expected maximum noise score over the comparison cube.
	The comparison cube, the detection subspace, and the noise law
	are fixed, so \(\eta_{\mathrm{peak}}\) does not depend on \(L\).
Thus, the high-probability upper bound
	in~\eqref{eq:search_domain_noise_bound} on the maximum noise
score over the search domain \(\overline{\mathcal B}(0,L)\)
scales like \(\log M_L\).

	\subsection{Simultaneous recovery, error control, and power}
	\label{subsec:error_control_power}
	In this subsection we address two questions.
		(1) Under what conditions does the probability of simultaneous
		recovery (reporting a location in every detection cube) tend to one
		as \(L\to\infty\)?
	(2) Under what conditions does the power of the detection procedure
	at the BH threshold tend to one as \(L\to\infty\)?
		See Subsection~\ref{subsec:performance_criteria} for the definitions. 
	We first consider how Step~\ref{item:candidate_peaks}
		selects candidate locations in the absence of noise.

	Recall from Subsection~\ref{subsec:algorithm_outline} that
	\(\mathcal B(\tau_i,\delta)\) is the detection cube of
	object instance \(i\), where \(0<\delta\le2B\),
and that selecting a candidate location \(t\) removes
	the cube \(\mathcal B(t,r)\) from the search domain
	\(\overline{\mathcal B}(0,L)\), where \(r=2B+\delta\).
	Define the clean score map by
	\begin{equation}\label{eq:Sx}
		S^x(t)
		:=
		\sum_{j=1}^{M}
		\bigl(x\ast\tilde{\psi}_j\bigr)^2(t).
	\end{equation}
	
Let \(P_{\mathcal V}\) denote the orthogonal projection in
\(L^2(D_{B/2})\) onto the detection subspace \(\mathcal V\)
defined in~\eqref{eq:detection_subspace_def}.
By~\eqref{eq:template_instance_model}, translating object
	instance \(x_i\) to the origin gives \(b_iR_i u_{\ell_i}\), which
	is supported in \(D_{B/2}\).
Define the coefficient vector \(a_i=(a_{i,1},\dots,a_{i,M})\in\mathbb R^M\)
of \(P_{\mathcal V}(b_iR_i u_{\ell_i})\) with respect to the
orthonormal basis \(\{\psi_j\}_{j=1}^M\) by
\begin{equation}\label{eq:x}
	a_{i,j}
	:=
	\big\langle b_iR_i u_{\ell_i},\psi_j
	\big\rangle_{L^2(D_{B/2})}
	=
	\bigl(x_i\ast\tilde{\psi}_j\bigr)(\tau_i),
	\qquad j=1,\dots,M,
\end{equation}
where the second equality follows from the definition
	of \(x_i\) in~\eqref{eq:template_instance_model}.

 Let \(S^{x_i}\) denote the score map
	in~\eqref{eq:Sx} with \(x_i\) in place of \(x\).
	By~\eqref{eq:x} and orthonormality,
\[
S^{x_i}(\tau_i)
= \|a_i\|^2
= \bigl\|P_{\mathcal V}(b_iR_i u_{\ell_i})\bigr\|_{L^2(D_{B/2})}^2.
\]
Set
\begin{equation}\label{eq:aL_def}
	a_L:=\min_{1\le i\le N}\|a_i\|.
\end{equation}

\begin{assumption}[Localization property]
	\label{def:localization_property}
	For each \(L\), label the object instances so that
	\[
	\|a_1\|\ge\|a_2\|\ge\cdots\ge\|a_N\|.
	\]
	For the prescribed localization parameter \(0<\delta\le2B\),
		there exists a constant \(\rho>0\), independent of \(L\), with the
		following property. For every \(L\) and every \(1\le i\le N\),
		there exists a representative point
	\begin{equation}\label{eq:qi}
		q_i\in\mathcal B(\tau_i,\delta)
		\cap\overline{\mathcal B}(0,L)
	\end{equation}
	such that, for every
	\begin{equation}\label{eq:comparison set}
		t\in\overline{\mathcal B}(0,L)\setminus
		\left(
		\bigcup_{j=1}^{i-1}\mathcal B(\tau_j,2B)
		\ \cup\
		\bigcup_{j=i}^{N}\mathcal B(\tau_j,\delta)
		\right),
	\end{equation}
	we have
\begin{equation}\label{eq:localization_score_gap}
	S^x(q_i)-S^x(t)\ge\rho\|a_i\|^2.
\end{equation}		
	We call \(\rho\) a localization constant and the set
	in~\eqref{eq:comparison set} the comparison set for object
	instance \(i\).
\end{assumption}

Equation~\eqref{eq:localization_score_gap} implies that if \(\|a_i\|>0\), then the clean score at \(q_i\), given by $S^x(q_i)$, is higher 
	than the score at every point in the comparison set of object instance \(i\).
Proposition~\ref{prop:template_localization_condition} gives a
sufficient condition on the known templates and the
	detection subspace for Assumption~\ref{def:localization_property}
	to hold at the prescribed \(\delta\). Under condition~\eqref{eq:template_localization_condition},
	\eqref{eq:rho_chi_def} below gives an explicit choice of \(\rho\).

\begin{assumption}[Separation of object instances]
	\label{assumptions}
	The locations \(\tau_i\in\overline{\mathcal B}(0,L)\) satisfy,
	for every pair of distinct indices \(i,j\),
	\begin{equation}\label{eq:separation}
		\|\tau_i-\tau_j\|_\infty>B+\frac{3\delta}{2}.
	\end{equation}
\end{assumption}

When \(a_L>0\), we use~\eqref{eq:localization_score_gap}
	and~\eqref{eq:separation} to show that, in the absence of noise,
	Step~\ref{item:candidate_peaks} selects its first \(N\) candidate
	locations from the detection cubes \(\mathcal B(\tau_i,\delta)\),
	one from each cube. The proof is given in
	Lemma~\ref{S-lem:clean_candidate_selection} of the supplementary
	material.

In the noisy image \(y=x+z\), detecting object instance
	\(i\) requires that Step~\ref{item:candidate_peaks} selects a candidate
	location in \(\mathcal B(\tau_i,\delta)\) whose observed score is at
	least the reporting threshold in Step~\ref{step:multiple_testing}.
We next bound the probability of detecting all object
instances, whose rotations are unknown, in terms of
\(a_L^2\), where \(a_L\) is the smallest projection norm among
the object instances, as defined in~\eqref{eq:aL_def}. To this end, we first specify a condition on the detection
subspace \(\mathcal V\) under which \(a_L\) does not depend on the
unknown rotations \(R_i\) in~\eqref{eq:template_instance_model};
see Lemma~\ref{S-lem:proj_rot_appendix} of the supplementary material.
Specifically, for \(R\in\SO(d)\), write
\(R\mathcal V:=\{Rf:f\in\mathcal V\}\).
We call the finite-dimensional detection subspace \(\mathcal V\)
steerable if it is invariant under rotations, meaning that
\begin{equation}\label{def:rotation_invariant}
	R\mathcal V=\mathcal V
	\qquad\text{for every }R\in\SO(d).
\end{equation}

	\begin{proposition}[Simultaneous recovery and Bonferroni error bounds]
		\label{prop:asymptotic_error_power}
		Fix \(0<\alpha<1\). Suppose that \(\mathcal V\)
		is steerable in the sense of~\eqref{def:rotation_invariant}
		and that Assumptions~\ref{ass:noise_reg_main},
			\ref{def:localization_property}, and~\ref{assumptions} hold.
		There exists a constant \(\Gamma_{\mathcal V}>0\), depending
		only on \(\mathcal V\), the template collection, and the
		localization constant \(\rho\), for which the following
		bounds hold. Set
		\[
		\varepsilon_L
		:=
		M_Lp\!\left(\frac{2a_L^2}{\Gamma_{\mathcal V}}\right),
		\]
		where \(M_L\) and \(p\) are defined in~\eqref{eq:ML_def}
		and~\eqref{eq:test_value_intro}, respectively.
		For every \(L\),
		\[
		\FWER(u_{\rm Bon})
		\le \min\{1,\alpha+\varepsilon_L\}.
		\]
		For every \(L\) with \(\varepsilon_L\le\alpha\), the
		simultaneous recovery probability
		satisfies
        \begin{equation}\label{eq:simultaneous_recovery_finite_bound}
			\Prob\!\left\{
			T(u)\cap\mathcal B(\tau_i,\delta)\neq\emptyset
			\text{ for every }i=1,\dots,N
			\right\}
			\ge 1-\varepsilon_L
		\end{equation}
		for \(u\in\{u_{\rm Bon},u_{\rm BH}\}\).
	\end{proposition}
	
	\begin{corollary}[Asymptotic simultaneous recovery]
		\label{cor:simultaneous_recovery_asymptotic}
		Under the assumptions of
		Proposition~\ref{prop:asymptotic_error_power},
		let \(\Gamma_{\mathcal V}\) be the constant given there. If
        \begin{equation}\label{eq:object_strength_condition}
			\liminf_{L\to\infty}
			\frac{a_L^2}{\eta_{\mathrm{peak}}\log M_L}
			>
			\Gamma_{\mathcal V},
		\end{equation}
		where \(\eta_{\mathrm{peak}}\) is defined
		in~\eqref{eq:eta_peak_def}, then, for 
		\(u\in\{u_{\rm Bon},u_{\rm BH}\}\), the
		simultaneous recovery probability at threshold \(u\)
		tends to one as \(L\to\infty\). Moreover,
		\[
		\limsup_{L\to\infty}\FWER(u_{\rm Bon})\le\alpha.
		\]
	\end{corollary}
	
	The proofs, including an explicit choice of
	\(\Gamma_{\mathcal V}\), are given in
	Section~\ref{S-app:proofs_theoretical_foundations}
	of the supplementary material.
	
Corollary~\ref{cor:simultaneous_recovery_asymptotic} gives
	a sufficient condition for detecting all object instances with
	probability tending to one. We next show that a weaker condition on the smallest
squared projection norm suffices for BH to detect
nearly all object instances on average.

	\begin{proposition}[Power at the BH threshold]
		\label{prop:bh_density_strength}
		Fix \(0<\alpha<1\). Suppose that \(\mathcal V\)
		is steerable in the sense of~\eqref{def:rotation_invariant}
		and that Assumptions~\ref{ass:noise_reg_main},
			\ref{def:localization_property}, and~\ref{assumptions} hold. There are an absolute constant
		\(C_{\rm BH}>0\) and constants \(A_{\rm BH},C,c>0\), depending only
		on the templates, \(\mathcal V\), \(d\), \(B\), \(\delta\),
		\(\rho\), and the noise law, such that, for every \(L\)
		satisfying both inequalities
		\begin{equation}\label{eq:bh_finite_power_conditions}
			a_L^2\ge A_{\rm BH},
			\qquad
			\frac{a_L^2}{\eta_{\mathrm{peak}}}
			\ge C_{\rm BH}\log\!\left(\frac{2M_L}{\alpha N}\right),
		\end{equation}
		we have
		\begin{equation}\label{eq:bh_average_power_rate}
			1-\Power(u_{\rm BH})
			\le C\exp\!\left\{-c\frac{a_L^2}{\eta_{\mathrm{peak}}}\right\}.
		\end{equation}
	\end{proposition}
	
	\begin{corollary}[Asymptotic power at the BH threshold]
		\label{cor:bh_power_asymptotic}
		Under the assumptions of
		Proposition~\ref{prop:bh_density_strength}, if
		\begin{equation}\label{eq:bh_log_density_balance}
			a_L^2\longrightarrow\infty
			\quad\text{and}\quad
			\frac{a_L^2}{\eta_{\mathrm{peak}}}
			\ge C_{\rm BH}\log\!\left(\frac{2M_L}{\alpha N}\right)
			\ \text{for all sufficiently large }L,
		\end{equation}
		then \(\Power(u_{\rm BH})\to1\). In particular, this holds if
		\begin{equation}\label{eq:positive_density_weak_strength}
			a_L^2\longrightarrow\infty,
			\qquad
			\liminf_{L\to\infty}\frac{N}{M_L}>0.
		\end{equation}
	\end{corollary}
	
The BH selection rule in~\eqref{eq:k_BH_def}
explains the role of the ratio \(N/M_L\)
in Corollary~\ref{cor:bh_power_asymptotic}.
Indeed, to detect at least \(i\) distinct object instances,
where \(1\le i\le N\), it suffices that each of~\(i\) distinct detection cubes contains a candidate location \(t\)
with \(p(S^y(t))\le\alpha i/M_L=\alpha(i/N)(N/M_L)\).
Recall that \(M_L\) counts the cubes \(Q_m\) intersecting
the search domain \(\overline{\mathcal B}(0,L)\), as defined
in~\eqref{eq:ML_def}.
By the BH rule in~\eqref{eq:k_BH_def}, all these candidate
locations are reported, so at least a fraction \(i/N\) of the
object instances is detected. If \(N/M_L\to0\) as \(L\to\infty\),
the number of object instances grows more slowly than the image
volume. For detecting a given fraction \(i/N\) of the object
instances, the bound \(\alpha i/M_L\) then tends to zero.
Satisfying \(p(S^y(t))\le\alpha i/M_L\) therefore requires smaller
test values and hence larger observed scores. If the number of
object instances grows with the image volume so that \(N/M_L\)
stays bounded away from zero as \(L\to\infty\), the bound
\(\alpha i/M_L\) for a fixed fraction \(i/N\) also stays bounded
away from zero. Consequently, observed scores above a fixed
threshold, independent of \(L\), suffice for these \(i\) candidate
locations to satisfy \(p(S^y(t))\le\alpha i/M_L\) and be reported
by BH. In this case, Corollary~\ref{cor:bh_power_asymptotic}
guarantees power tending to one even when \(a_L^2\) tends to
infinity arbitrarily slowly.

The proofs, including explicit constants, are given in
Subsection~\ref{S-app:bh_density_strength} of the supplementary
material.

	To obtain FDR control at the BH threshold when power
		tends to one, we also bound the number of false detections.
		As the image grows, there are more cubes \(Q_m\) in which noise
		can produce scores at or above a given threshold, while the
		number of object instances may grow more slowly. We use the
		following variance bound to control the number of such cubes.
		When \(N/M_L\) stays bounded away from zero, ergodicity of the
		noise field under translations provides an alternative.
	
	\begin{assumption}[Variance bound]
		\label{ass:noise_exceedance_variance}
		There exists a constant \(0<C_{\rm var}<\infty\), independent
		of \(L\) and \(u\), such that
		\begin{equation}\label{eq:noise_exceedance_variance}
			\operatorname{Var}\!\left(
			\#\left\{m\in\mathcal I_L:
			\max_{t\in\overline{Q_m}}S^z(t)\ge u\right\}
			\right)
			\le C_{\rm var}M_Lp(u),
			\qquad L>0,\quad u\in\mathbb R.
		\end{equation}
	\end{assumption}
	
To explain Assumption~\ref{ass:noise_exceedance_variance},
note that under Assumption~\ref{ass:noise_reg_main}, the maximum
noise score over each cube \(\overline Q_m\), \(m\in\mathcal I_L\),
reaches or exceeds \(u\) with probability \(p(u)\), by stationarity and
Theorem~\ref{thm:continuity_test}.
The number of cubes counted
in~\eqref{eq:noise_exceedance_variance} equals
\[
\sum_{m\in\mathcal I_L}
\mathbf 1_{\left\{
\max_{t\in\overline{Q_m}}S^z(t)\ge u
\right\}}.
\]
Each indicator has expectation \(p(u)\), so the sum
has expectation \(M_Lp(u)\).
Assumption~\ref{ass:noise_exceedance_variance} requires
the variance in~\eqref{eq:noise_exceedance_variance}
to be at most a constant multiple of the mean \(M_Lp(u)\).
For example, if the noise maxima over distinct cubes
were independent,
the bound would hold with \(C_{\rm var}=1\), because
the variance would be
\[
\sum_{m\in\mathcal I_L}p(u)(1-p(u))
=
M_Lp(u)(1-p(u))
\le M_Lp(u).
\]
The variance bound in~\eqref{eq:noise_exceedance_variance}
can also hold for correlated noise, with a suitable constant
\(C_{\rm var}\).
Under Assumption~\ref{ass:noise_reg_main},
one sufficient condition is that
\(\Cov(z(0),z(h))=0\) whenever \(\|h\|_\infty\) exceeds a fixed
constant (Lemma~\ref{S-lem:noise_exceedance_variance_finite_range}).

	\begin{proposition}[Asymptotic FDR control at the BH threshold]
		\label{prop:bh_fdr_control}
		Fix \(0<\alpha<1\). Suppose that
		Assumptions~\ref{ass:noise_reg_main} and~\ref{assumptions} hold
		and that \(\Power(u_{\rm BH})\to1\).
		If \(N\to\infty\) and
		Assumption~\ref{ass:noise_exceedance_variance} holds, then
		\[
		\limsup_{L\to\infty}\FDR(u_{\rm BH})\le\alpha.
		\]
		The same conclusion holds if
we replace \(N\to\infty\) and
Assumption~\ref{ass:noise_exceedance_variance} by
\(\liminf_{L\to\infty}N/M_L>0\) and
ergodicity of \(z\) under translations.
	\end{proposition}
	
	The proof is given in
	Subsection~\ref{S-app:bh_density_strength} of the supplementary
	material, and Lemma~\ref{S-lem:bh_fdr_general_bound} gives an FDR
	bound for every \(L\). The noise field used in the simulations, with covariance
in~\eqref{eq:simulation_correlated_noise}, is ergodic
under translations, as shown in
Proposition~\ref{S-prop:noise_assumptions_isotropic_exponential}.
	
\begin{remark}
	\label{rem:threshold_orders}
	Fix \(0<\alpha<1\), and suppose that
	Assumption~\ref{ass:noise_reg_main} holds.
The Bonferroni threshold \(u_{\rm Bon}\)
	in~\eqref{eq:u_bon_def} is deterministic and depends on the
	image side length \(L\) only through the number \(M_L\) of cubes
	in~\eqref{eq:ML_def}. As the image grows, \(u_{\rm Bon}\)
	grows logarithmically with \(L\), by
Corollary~\ref{S-cor:calibrated_thresholds} of the
supplementary material.
	
	The BH threshold \(u_{\rm BH}\) in~\eqref{eq:u_bh_def}
		is random because the index \(k\) in~\eqref{eq:k_BH_def} is determined by the observed test values at the candidate
		locations.
	Suppose in addition that Assumption~\ref{assumptions} holds
	and that \(\Power(u_{\rm BH})\to1\). If \(N/M_L\to0\) and
	Assumption~\ref{ass:noise_exceedance_variance} holds, then
	with probability tending to one, \(u_{\rm BH}\)
		grows logarithmically with \(M_L/N\), the number of cubes
		per object instance.
	See Lemma~\ref{S-lem:bh_threshold_sparse} of the supplementary
	material.
	If instead \(\liminf_{L\to\infty}N/M_L>0\), then, with
	probability tending to one, the BH threshold is bounded above
	by a constant independent of \(L\).
	See Lemma~\ref{S-lem:bh_threshold_positive_density} of the
	supplementary material.
\end{remark}

\subsection{A lower bound on the smallest squared projection norm}
\label{subsec:object_strength_lower_bound}

When \(N\to\infty\) and \(N/M_L\to0\),
	Corollary~\ref{cor:bh_power_asymptotic} shows that BH has power
	tending to one if \(a_L^2/\eta_{\mathrm{peak}}\) eventually exceeds
	a sufficiently large multiple of \(\log(M_L/N)\).
	Under the variance bound in Assumption~\ref{ass:noise_exceedance_variance},
	Proposition~\ref{prop:bh_fdr_control} also gives asymptotic FDR control.
We now ask whether another procedure that reports locations
	at least \(r/2\) apart in the \(\ell_\infty\) norm could give both
	guarantees under a weaker condition on the smallest squared projection
	norm \(a_L^2\).

Fix \(0<\alpha<1\). For \(L>B\), let \(\mathcal A_{L,N}\) be the set of
all \(N\)-tuples of object instance locations
\(\boldsymbol\tau=(\tau_1,\dots,\tau_N)\) satisfying~\eqref{eq:separation},
with \(\tau_i\in\overline{\mathcal B}(0,L-B)\) for every \(i\).
Consider any detection procedure that observes the image
	\(y\) of~\eqref{eq:the_model} on \(\overline{\mathcal B}(0,L+B)\).
	Let \(\widehat T_L\) denote the finite set of locations it reports
	in \(\overline{\mathcal B}(0,L)\). Write
\(\FDR_{\boldsymbol\tau}(\widehat T_L)\) and
\(\Power_{\boldsymbol\tau}(\widehat T_L)\) for its FDR and power
when the object instances have locations \(\boldsymbol\tau\).
The goal is a procedure with
\begin{equation}\label{eq:uniform_fdr}
	\limsup_{L\to\infty}
	\sup_{\boldsymbol\tau\in\mathcal A_{L,N}}
	\FDR_{\boldsymbol\tau}(\widehat T_L)\le\alpha
\end{equation}
and
\begin{equation}\label{eq:uniform_power}
	\lim_{L\to\infty}
	\inf_{\boldsymbol\tau\in\mathcal A_{L,N}}
	\Power_{\boldsymbol\tau}(\widehat T_L)=1,
\end{equation}
that is, asymptotic FDR control and power tending to one
	for every allowed choice of the object instance locations.

\begin{proposition}[Optimality of the logarithmic order]
	\label{prop:retained_energy_lower_bounds}
	There exist a model of the form~\eqref{eq:template_instance_model}
		and a constant \(c_{\rm low}>0\) with the following properties.
		The templates, the steerable detection subspace, and the noise law
		are fixed as \(L\) varies. Assumptions~\ref{ass:noise_reg_main}
		and~\ref{ass:noise_exceedance_variance} hold. For every
		\(\boldsymbol\tau\in\mathcal A_{L,N}\), the clean score map
		satisfies Assumption~\ref{def:localization_property} with \(\rho=1\).
	
	Fix \(0<\alpha<1\), let \(N=N(L)\ge1\) satisfy \(N/M_L\to0\),
	and let every object instance have squared projection
	norm \(a_L^2>0\). Suppose that
	\begin{equation}\label{eq:intermediate_lower_strength}
		\limsup_{L\to\infty}
		\frac{a_L^2}{\eta_{\mathrm{peak}}\log(M_L/N)}
		<c_{\rm low}.
	\end{equation}
	Consider a sequence of detection procedures, one for
each \(L\), such that any two locations reported by
each procedure are at \(\ell_\infty\)-distance at least
\(r/2\).
		If their reported sets \(\widehat T_L\) satisfy~\eqref{eq:uniform_fdr},
	then
	\[
	\limsup_{L\to\infty}
	\inf_{\boldsymbol\tau\in\mathcal A_{L,N}}
	\Power_{\boldsymbol\tau}(\widehat T_L)\le\alpha,
	\]
	so~\eqref{eq:uniform_power} fails.
\end{proposition}

The proof is given in
Section~\ref{S-app:intermediate_density_lower_bound} of the
supplementary material.
	
	\section{Steerable subspaces and pSNR selection}
	\label{sec:repr_framework}
	
	Section~\ref{sec:theoretical_foundations} established guarantees for
	simultaneous recovery, power, and error control for a fixed detection
	subspace. We now address the choice of this subspace using the known
	templates and noise law. We consider detection subspaces that are steerable in the
	sense of~\eqref{def:rotation_invariant}, so the projection norm of an
	object instance does not depend on its rotation.
Once we choose the detection subspace, we use an orthonormal
	basis of that subspace to compute the observed score map \(S^y\)
	in~\eqref{eq:smoothed_process_intro}.
	To choose this subspace, we first construct a
    finite-dimensional steerable ambient space that approximates
    every template and its rotations uniformly. We then form a
    finite family of candidate detection subspaces by
    choosing subsets of the ambient space.
    For the prescribed localization parameter \(\delta\), we
    retain as candidate detection subspaces those satisfying
    the sufficient conditions for localization in
    Proposition~\ref{prop:template_localization_condition}.
	If this family is nonempty, we select a subspace from it by maximizing the
	projected signal-to-noise ratio (pSNR) in~\eqref{eq:psnr_V}, computed from
	the templates and noise law. Finally,
	Corollary~\ref{cor:asymptotic_selected} gives detection and error control
	guarantees for the selected detection subspace.
	
	\subsection{Constructing a steerable ambient space}
	\label{subsec:construction_setting}
	
	For a steerable subspace \(\mathcal V\), as defined in
	\eqref{def:rotation_invariant}, projection and rotation satisfy, for every
	\(f\in L^2(D_{B/2})\) and \(R\in\SO(d)\),
	\[
	P_{\mathcal V}(Rf)=R P_{\mathcal V}f;
	\]
	see Lemma~\ref{S-lem:proj_rot_appendix}. 
	To construct steerable subspaces, we expand functions in polar
	coordinates and group the angular components by frequency. For
	\(d=2\), expand \(f(\varrho,\cdot)\), as a function of the angle
	\(\omega\in[0,2\pi)\), in the Fourier basis
	\(\{1\}\cup\{\cos(j\omega),\sin(j\omega):j\ge1\}\). A rotation by an
	angle \(\theta\) replaces \(\omega\) by \(\omega-\theta\), and maps
	\(\cos(j\omega)\) and \(\sin(j\omega)\) to linear combinations of
	themselves. The span of the frequency-\(j\) terms is therefore
	invariant under rotations.
	For \(d>2\), spherical harmonics play the role of the
	Fourier basis.
	The spherical harmonics of degree \(j\) are the restrictions to
	\(S^{d-1}\) of homogeneous harmonic polynomials of degree \(j\).
	In spherical coordinates \((\varrho,\omega)\), with \(\omega\in S^{d-1}\), let \(\{Y_{j,m}\}_{m=1}^{h_{d,j}}\) be a
	real orthonormal basis, with respect to surface area measure
	\(d\sigma\), of spherical harmonics of degree~\(j\), where
	\(h_{d,j}\) is the dimension of this space. Rotations map spherical
	harmonics of degree \(j\) to spherical harmonics of the same degree,
	and harmonics of different degrees are orthogonal; see
	Lemma~\ref{S-lem:real_angular_expansion}.
	For a nonnegative integer \(J\) and finite-dimensional spaces
	\(V_j\subset L^2([0,B/2],\varrho^{d-1}\,d\varrho)\) of radial
	functions, define the linear spaces
	\[
	\mathcal S_j
	:=\left\{\sum_{m=1}^{h_{d,j}}a_m(\varrho)Y_{j,m}(\omega):
	a_m\in V_j\text{ for every }m\right\},
	\qquad j=0,\ldots,J.
	\]
	We call \(\mathcal S_j\) the angular component of degree \(j\).
		Each of its elements combines degree-\(j\) angular functions with
	radial coefficients in \(V_j\). For \(d=2\) and
	\(j\ge1\), these are the functions
	\(a_1(\varrho)\cos(j\omega)+a_2(\varrho)\sin(j\omega)\) with
	\(a_1,a_2\in V_j\); a rotation changes only the coefficients, which
	remain in \(V_j\). Define
	\begin{equation}\label{eq:ambient_space}
		\mathcal S(J;V_0,\ldots,V_J):=\bigoplus_{j=0}^{J}\mathcal S_j,
	\end{equation}
	the space of functions whose angular expansion contains only degrees
	\(0,\ldots,J\), with the degree-\(j\) radial coefficients in \(V_j\).
	Each \(\mathcal S_j\) is steerable, so their direct sum
		in~\eqref{eq:ambient_space} is a finite-dimensional steerable space.
	See Lemma~\ref{S-lem:real_angular_expansion}.
	\begin{proposition}[Uniform approximation by a steerable ambient space]
		\label{prop:ambient_uniform_coverage}
		For every \(\varepsilon>0\), one can choose \(J\) and
		\(V_0,\ldots,V_J\) from the known templates $u_{1},\ldots, u_{n}$, so that the space
		\(\mathcal S=\mathcal S(J;V_0,\ldots,V_J)\) satisfies
		\begin{equation}\label{eq:uniform_template_approximation}
			\sup_{1\le\ell\le n}\sup_{R\in\SO(d)}
			\|Ru_\ell-P_{\mathcal S}(Ru_\ell)\|_{L^2(D_{B/2})}
			\le\varepsilon.
		\end{equation}
	\end{proposition}
	
	The proof is given in Section~\ref{S-app:rotation_invariant} of
	the supplementary material. We fix \(0<\varepsilon<1\) and a space
	\(\mathcal S\) of the form~\eqref{eq:ambient_space} satisfying
	\eqref{eq:uniform_template_approximation}, independent of \(L\).
	Since the templates have unit norm, \(\varepsilon<1\) ensures that every
	template has nonzero projection onto \(\mathcal S\). 
	
	\subsection{Selecting the detection subspace}
	\label{subsec:candidate_score_maps}
	
The ambient space \(\mathcal S\)
of Proposition~\ref{prop:ambient_uniform_coverage}
approximates all rotated templates. However, our goal is not to approximate the templates. Rather, for
 the detection procedure of Section~\ref{sec:model_detection},
	we want a basis that represents the templates well, so that \(S^x\)
	is large near the object instances, and at the same time represents
	the noise poorly, so that a peak in \(S^y\) due to the signal \(x\)
	can be distinguished from a peak due to the noise \(z\).
	We therefore consider subspaces formed by taking sums of angular
	components of \(\mathcal S\), and call them candidate detection
	subspaces. Specifically, for \(I\subseteq\{0,\ldots,J\}\), let
	\(\mathcal S(I):=\bigoplus_{j\in I}\mathcal S_j\).
	Let \(\mathcal I\) be a collection of index sets, independent
	of \(L\). Define the family of candidate detection subspaces by
\begin{equation}\label{eq:candidate_family}
	\mathfrak S
	:=\{\mathcal S(I): I\in\mathcal I,\ \mathcal S(I)\ne\{0\}\}.
\end{equation}
	We first use Proposition~\ref{prop:template_localization_condition}
		to identify candidate detection subspaces that satisfy the localization property
		in Assumption~\ref{def:localization_property} at the prescribed \(\delta\).
		Among these subspaces, we seek one that best separates the object instances from the noise
(see~\eqref{eq:psnr_V}).

		To check that the localization property holds, consider a candidate subspace
        \(\mathcal V\in\mathfrak S\) onto which every template has nonzero
        projection. The localization property requires a representative
        point \(q_i\) as in~\eqref{eq:qi} for each object instance \(i\),
        such that
        \(S^x(q_i)-S^x(t)\ge\rho\|a_i\|^2\)
        for every \(t\) in the comparison set in~\eqref{eq:comparison set}.
   The comparison set for object instance \(i\)
consists of points in the search domain
outside the union of
the cubes \(\mathcal B(\tau_j,2B)\) for \(j<i\)
and the detection cubes
\(\mathcal B(\tau_j,\delta)\) for \(j\ge i\),
where the instances are ordered so that
\(\|a_1\|\ge\cdots\ge\|a_N\|\).
		The clean score and the comparison sets depend on the
unknown instance parameters in~\eqref{eq:template_instance_model}.
The following proposition gives a sufficient condition based
only on the known templates, \(\mathcal V\), \(B\), and \(\delta\)
that guarantees the localization property whenever
Assumption~\ref{assumptions} holds.

\begin{proposition}[Template condition for localization]
	\label{prop:template_localization_condition}
	Fix \(0<\delta\le2B\), and let \(\mathcal V\in\mathfrak S\) satisfy
	\(\|P_{\mathcal V}u_\ell\|_{L^2(D_{B/2})}>0\) for every \(\ell\).
	Suppose Assumption~\ref{assumptions} holds for every \(L\).
	
    Define
    \begin{equation}\label{eq:radial_projection_envelope}
		e_{\mathcal V}(\lambda)
		:=\max_{1\le\ell\le n}\sup_{\|h\|=\lambda}
		\frac{\|P_{\mathcal V}[u_\ell(h+\cdot)]\|_{L^2(D_{B/2})}}
		{\|P_{\mathcal V}u_\ell\|_{L^2(D_{B/2})}},
		\qquad 0\le\lambda\le B.
	\end{equation}
	For \(1\le k\le2^d\), let
	\(\mathcal C_k(B,\delta)\) be the set of tuples
	\((v_1,\ldots,v_k)\in(\mathbb R^d)^k\) satisfying
	\begin{equation}\label{eq:admissible_displacement_tuples}
		\begin{aligned}
			\|v_m\|\le B,\qquad \|v_m\|_\infty\ge\delta/2
			&\qquad (1\le m\le k),\\
			\|v_m-v_{m'}\|_\infty\ge B+3\delta/2
			&\qquad (1\le m<m'\le k).
		\end{aligned}
	\end{equation}
	Set
	\begin{equation}\label{eq:chi_V_def}
		\chi_{\mathcal V}(\delta)
		:=\max_{1\le k\le2^d}
		\sup_{(v_1,\ldots,v_k)\in\mathcal C_k(B,\delta)}
		\sum_{m=1}^k e_{\mathcal V}(\|v_m\|),
	\end{equation}
with the supremum over an empty set taken to be zero.
	If	\begin{equation}\label{eq:template_localization_condition}
		\chi_{\mathcal V}(\delta)<1,
	\end{equation}
	then the clean score map in~\eqref{eq:Sx}, computed using an
	orthonormal basis of \(\mathcal V\), satisfies
	the localization property in
	Assumption~\ref{def:localization_property}, with
	\(q_i=\tau_i\) and
	\begin{equation}\label{eq:rho_chi_def}
		\rho=1-\chi_{\mathcal V}(\delta)^2.
	\end{equation}
\end{proposition}
	
The intuition behind
Proposition~\ref{prop:template_localization_condition}
is as follows.
At a point \(t\) in the comparison set~\eqref{eq:comparison set}
for instance \(i\), several object instances may contribute
to the clean score.
The condition \(\chi_{\mathcal V}(\delta)<1\)
ensures that the clean score at \(\tau_i\) exceeds
the clean score throughout the comparison set for
instance \(i\), with the score gap required
in~\eqref{eq:localization_score_gap}.
The proof is given in
Section~\ref{S-app:psnr_selection_proofs} of the supplementary material.

The candidate detection subspaces in \(\mathfrak S\),
defined in~\eqref{eq:candidate_family},
are steerable by construction. To apply the simultaneous recovery
and power guarantees of Section~\ref{sec:theoretical_foundations},
we also require the localization property of
Assumption~\ref{def:localization_property}.
The conditions \(P_{\mathcal V}u_\ell\ne0\) for every
\(\ell\) and \(\chi_{\mathcal V}(\delta)<1\) in
Proposition~\ref{prop:template_localization_condition}
depend only on the known templates, and on \(\mathcal V\), \(B\),
and \(\delta\).
We therefore restrict the candidate family
\(\mathfrak S\) in~\eqref{eq:candidate_family}
to the family \(\mathfrak S_{\rm loc}\) defined by
\begin{equation}\label{eq:retained_detection_subspaces}
		\mathfrak S_{\rm loc}
		:=\left\{
		\mathcal V\in\mathfrak S :
		P_{\mathcal V}u_\ell\ne 0
		\text{ for every } \ell=1,\ldots,n,
		\quad \chi_{\mathcal V}(\delta)<1
		\right\}.
\end{equation}
	If \(\mathcal S\in\mathfrak S\), then
\(\mathfrak S_{\rm loc}\ne\emptyset\) for \(\delta=2B\)
and for \(\delta<2B\) sufficiently close to \(2B\), by
Remark~\ref{S-rem:nonempty_retained_family} of the
supplementary material.

We now choose a detection subspace from
\(\mathfrak S_{\rm loc}\).
The BH power bound in~\eqref{eq:bh_average_power_rate} is
\[
	1-\Power(u_{\rm BH})
	\le C\exp\!\left\{
	-c\,\frac{a_L^2}{\eta_{\mathrm{peak}}}
	\right\}.
\]
The sufficient condition for simultaneous recovery
in~\eqref{eq:object_strength_condition} is
\[
	\liminf_{L\to\infty}
	\frac{a_L^2}{\eta_{\mathrm{peak}}\log M_L}
	>
	\Gamma_{\mathcal V}.
\]
For fixed \(C\) and \(c\), the upper bound on
\(1-\Power(u_{\rm BH})\) is decreasing in
\(a_L^2/\eta_{\mathrm{peak}}\).
Condition~\eqref{eq:object_strength_condition} requires
\(a_L^2/\eta_{\mathrm{peak}}\) to exceed a sufficiently
large multiple of \(\log M_L\) for all sufficiently large \(L\).
We therefore seek a detection subspace for which
\(a_L^2/\eta_{\mathrm{peak}}\) is large.

We now write \(a_L(\mathcal V)\) and
\(\eta_{\mathrm{peak}}(\mathcal V)\) to make the dependence
of \(a_L\) and \(\eta_{\mathrm{peak}}\) on the detection
subspace explicit.
For each \(\mathcal V\in\mathfrak S_{\rm loc}\), the smallest
projection norm among the object instances, as defined
in~\eqref{eq:aL_def}, is
\begin{equation}\label{eq:aL_V_def}
	a_L(\mathcal V)
	:=\min_{1\le i\le N}
	\|P_{\mathcal V}(b_iR_i u_{\ell_i})\|_{L^2(D_{B/2})}.
\end{equation}
The expected maximum noise score over the comparison
cube, as defined in~\eqref{eq:eta_peak_def}, is
\begin{equation}\label{eq:eta_peak_V}
	\eta_{\mathrm{peak}}(\mathcal V)
	:=\E\!\left[\max_{t\in\overline{\mathcal B}(0,r/2)}
	S_{\mathcal V}^z(t)\right],
\end{equation}
where \(S_{\mathcal V}^z\) is the noise score map
in~\eqref{eq:noise_score_map} computed using an orthonormal
basis of \(\mathcal V\).

The ratio
\(a_L(\mathcal V)^2/\eta_{\mathrm{peak}}(\mathcal V)\)
depends on the noise law through
\(\eta_{\mathrm{peak}}(\mathcal V)\)
of~\eqref{eq:eta_peak_V}.
The definition of \(a_L(\mathcal V)\)
in~\eqref{eq:aL_V_def} involves the unknown amplitudes \(b_i\),
rotations \(R_i\), and template indices \(\ell_i\)
(see~\eqref{eq:template_instance_model}).
Since \(\mathcal V\) is steerable,
Lemma~\ref{S-lem:proj_rot_appendix} of the supplementary
material shows that \(a_L(\mathcal V)\) does not depend
on the rotations \(R_i\).
The amplitudes \(b_i\) and template indices \(\ell_i\)
remain unknown, so we use the lower bound in~\eqref{S-eq:psnr_retained_bound} of the supplementary material
\[
	\frac{a_L(\mathcal V)^2}
	{\eta_{\mathrm{peak}}(\mathcal V)}
	\ge
	\left(\min_{1\le i\le N}b_i\right)^2
	\frac{\displaystyle\min_{1\le\ell\le n}
		\|P_{\mathcal V}u_\ell\|_{L^2(D_{B/2})}^2}
	{\eta_{\mathrm{peak}}(\mathcal V)}.
\]
Since the amplitudes \(b_i\)
in~\eqref{eq:template_instance_model} do not depend on
the choice of detection subspace \(\mathcal V\),
we define the projected signal-to-noise ratio (pSNR) of
\(\mathcal V\) by
\begin{equation}\label{eq:psnr_V}
	\mathrm{pSNR}(\mathcal V)
	:=\frac{\displaystyle\min_{1\le\ell\le n}
		\|P_{\mathcal V}u_\ell\|_{L^2(D_{B/2})}^2}
	{\eta_{\mathrm{peak}}(\mathcal V)}.
\end{equation}
For each \(\mathcal V\in\mathfrak S_{\rm loc}\), we can
compute \(\mathrm{pSNR}(\mathcal V)\)
in~\eqref{eq:psnr_V} from the known templates and noise law.
Maximizing \(\mathrm{pSNR}(\mathcal V)\)
in~\eqref{eq:psnr_V} over \(\mathfrak S_{\rm loc}\)
therefore maximizes the lower bound on
\(a_L(\mathcal V)^2/\eta_{\mathrm{peak}}(\mathcal V)\),
where \(a_L(\mathcal V)\) and
\(\eta_{\mathrm{peak}}(\mathcal V)\) are defined
in~\eqref{eq:aL_V_def} and~\eqref{eq:eta_peak_V},
respectively.
Whenever \(\mathfrak S_{\rm loc}\ne\emptyset\), we choose
the detection subspace by taking
\begin{equation}\label{eq:Sstar_def}
	\mathcal S^\star\in
	\arg\max_{\mathcal V\in\mathfrak S_{\rm loc}}
	\mathrm{pSNR}(\mathcal V).
\end{equation}

The following corollary gives sufficient conditions
for BH power and simultaneous recovery when using the
detection subspace selected in~\eqref{eq:Sstar_def}.
\begin{corollary}[Detection guarantees for the pSNR criterion]
	\label{cor:asymptotic_selected}
	Fix \(0<\alpha<1\), and let \(b_L=\min_{1\le i\le N}b_i\),
	where \(b_i\) is the amplitude of object instance \(i\)
	in~\eqref{eq:template_instance_model}.
	Suppose that \(\mathfrak S_{\rm loc}\ne\emptyset\), that
	Assumption~\ref{ass:noise_reg_main} holds for an orthonormal
	basis of \(\mathcal S\) of~\eqref{eq:ambient_space}, and that
	Assumption~\ref{assumptions} holds.
	
	For every \(\mathcal V\in\mathfrak S_{\rm loc}\), the following
	guarantees hold for the detection procedure using \(\mathcal V\).
	
	\textup{(i) BH power.}
	There exist constants \(A_\star,C_\star,c_\star>0\), common to all
	\(\mathcal V\in\mathfrak S_{\rm loc}\), such that for every
	\(0<\beta<1\), we have \(\Power(u_{\rm BH})\ge 1-\beta\) whenever
	\begin{equation}\label{eq:selected_bh_target_power}
		b_L^2 \ge \frac{1}{\mathrm{pSNR}(\mathcal V)}
		\max\!\left\{
		A_\star,\,
		C_{\rm BH}\log\!\left(\frac{2M_L}{\alpha N}\right),\,
		\frac{1}{c_\star}\log\!\left(\frac{C_\star}{\beta}\right)
		\right\},
	\end{equation}
	where \(C_{\rm BH}\) is as in Proposition~\ref{prop:bh_density_strength}.
	
\textup{(ii) Simultaneous recovery.}
The probability
of simultaneous recovery tends to one at both the Bonferroni and BH
thresholds if
\begin{equation}\label{eq:selected_simultaneous_recovery}
	\liminf_{L\to\infty}\frac{b_L^2}{\log M_L}
	> \frac{\Gamma_\star}{\mathrm{pSNR}(\mathcal V)},
\end{equation}
where \(\Gamma_\star := \max_{\mathcal V\in\mathfrak S_{\rm loc}}
\Gamma_{\mathcal V}\), with \(\Gamma_{\mathcal V}\) as in
Proposition~\ref{prop:asymptotic_error_power}.
	
Among all \(\mathcal V\in\mathfrak S_{\rm loc}\), the
	subspace \(\mathcal S^\star\) of~\eqref{eq:Sstar_def} givesthe weakest requirements on \(b_L^2\) in~\eqref{eq:selected_bh_target_power} for BH power of at least \(1-\beta\)
	 and for simultaneous recovery
	in~\eqref{eq:selected_simultaneous_recovery}.
\end{corollary}

Condition~\eqref{eq:selected_simultaneous_recovery} also gives
asymptotic Bonferroni FWER control by
Corollary~\ref{cor:simultaneous_recovery_asymptotic}.
Moreover, when BH power tends to one,
the additional assumptions in
Proposition~\ref{prop:bh_fdr_control}
ensure asymptotic FDR control.
The proof of Corollary~\ref{cor:asymptotic_selected} is given in
Section~\ref{S-app:psnr_selection_proofs} of the supplementary material.

\section{Simulations}
\label{sec:simulations}

We evaluate the detection procedure of
Subsection~\ref{subsec:algorithm_outline} for \(d=2\).
We report the power and error rates defined in
Subsection~\ref{subsec:performance_criteria},
using the detection subspace \(\mathcal S^\star\)
selected by the pSNR selection rule
in~\eqref{eq:Sstar_def}.
The experiments cover several signal-to-noise
ratios (SNRs), defined in~\eqref{eq:simulation_snr}.
We also test the main point of
Section~\ref{sec:repr_framework}, that the detection
subspace should be chosen for detection rather than for
approximating the templates. To this end, we compare
\(\mathcal S^\star\) with the ambient space \(\mathcal S\),
which is constructed to approximate the templates and
their rotations, as
in~\eqref{eq:uniform_template_approximation}, and with
the other candidate detection subspaces that pass the
numerical check of the localization conditions
in~\eqref{eq:retained_detection_subspaces}.

\begin{figure}[htbp]
    \centering
    \begin{tabular}{@{}c@{\hspace{0.02\textwidth}}c@{}}
        \includegraphics[width=0.15\textwidth]
        {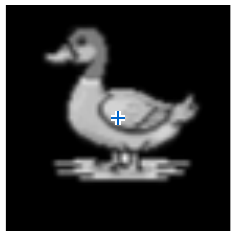}
        &
        \includegraphics[width=0.15\textwidth]
        {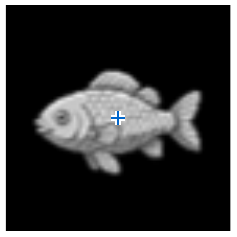}
        \\[3pt]
        \(u_1\) & \(u_2\)
    \end{tabular}
    \caption{The unit-norm templates \(u_1\) (duck) and
\(u_2\) (fish), with plus signs marking their
centers of mass.}
    \label{fig:duck_fish_templates}
\end{figure}

We use two templates, a duck \(u_1\) and a fish \(u_2\),
shown in Figure~\ref{fig:duck_fish_templates}.
Both templates are normalized to unit norm, as in
Subsection~\ref{subsec:background_setting}.
For each template, we take its center of mass,
marked by a plus sign in the figure, as its origin.
Thus, \(\tau_i\) in~\eqref{eq:template_instance_model}
is the location of the center of mass of object instance \(i\).

We set the prescribed localization parameter to \(\delta=20\).
Following~\eqref{eq:true_noise_regions}, a reported location~\(t\) is counted as a true detection of object instance
\(i\) if it lies in the detection cube
\(\mathcal B(\tau_i,\delta)\), that is, if
\(\|t-\tau_i\|_\infty<10\) pixels.
Reports outside all detection cubes are counted
as false detections.
We set the side length of the search domain
\(\overline{\mathcal B}(0,L)\) to \(L=1024\) pixels.
We take \(B=72\), with both templates supported in
the disk \(D_{B/2}\), as defined in
Subsection~\ref{subsec:background_setting}.
Each clean image \(x\) in~\eqref{eq:template_instance_model}
contains \(N=80\) object instances, with~40 copies of~\(u_1\) and~40 copies of~\(u_2\).
We generate the object locations at random, with
\(\|\tau_i-\tau_j\|_\infty>B+3\delta/2=102\) pixels
for every \(i\ne j\), as required by
Assumption~\ref{assumptions}.
For each object instance, we independently draw a rotation
angle uniformly from \([0,2\pi)\) and an amplitude \(b_i\)
uniformly from \([0.8,1.2]\).
Each rotated copy is scaled to have norm \(b_i\),
as in~\eqref{eq:bi}.

\begin{figure}[htbp]
    \centering
    \begin{tabular}{@{}c@{\hspace{0.02\textwidth}}c@{}}
        \includegraphics[width=0.46\textwidth]
        {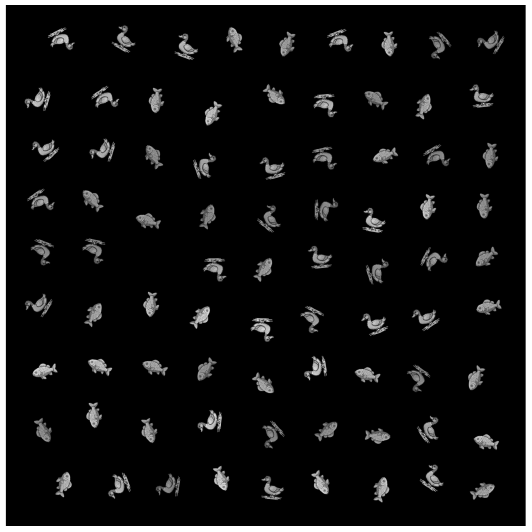}
        &
        \includegraphics[width=0.46\textwidth]
        {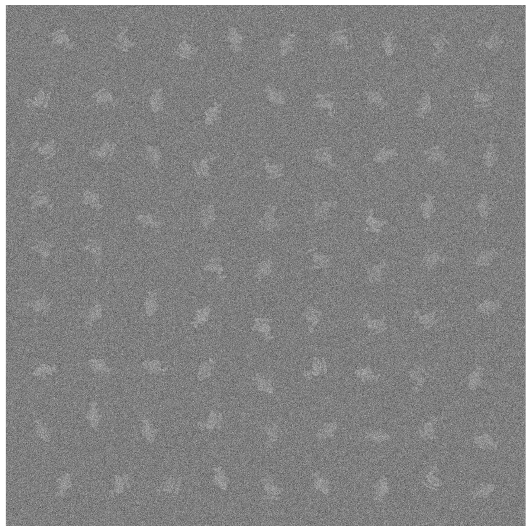}
        \\[3pt]
        \(x\) & \(y\)
    \end{tabular}
    \caption{A clean image \(x\) (left) and the corresponding
    noisy observation \(y\) at \(\mathrm{SNR}=0.1\) (right).}
    \label{fig:duck_fish_clean_noisy}
\end{figure}

We generate the observed image \(y=x+z\), as
in~\eqref{eq:the_model}, by adding a centered stationary
Gaussian noise field \(z\), independent of \(x\), with
known covariance
\begin{equation}
    \label{eq:simulation_correlated_noise}
    \operatorname{Cov}(z(t),z(s))
    =\sigma^2\,10^{-3\|t-s\|}.
\end{equation}
We define the signal-to-noise ratio by
\begin{equation}
    \label{eq:simulation_snr}
    \mathrm{SNR}
    =
    \frac{\|u_\ell\|_{L^2(\mathbb R^2)}^2}
    {\E\!\left[\|z\|_{L^2(\mathcal B(0,B))}^2\right]}
    =
    \frac{1}{\sigma^2 B^2},
    \qquad \ell=1,2,
\end{equation}
where \(\mathcal B(0,B)\) is the square of side length \(B\)
defined in Subsection~\ref{subsec:background_setting}.
The last equality follows because both templates have
unit norm and the noise has variance \(\sigma^2\)
at every location.
We consider
\[
    \mathrm{SNR}\in
    \{0.01,0.015,0.02,0.05,0.075,0.1,0.2\}.
\]
We generate 500 clean images independently.
For each clean image, we generate a Gaussian noise
field with the covariance
in~\eqref{eq:simulation_correlated_noise} at \(\sigma=1\),
independently of the clean image and of the noise
generated for the other images.
For each clean image, we use this same noise field
at all SNR levels, scaling it to obtain the prescribed
SNR according to~\eqref{eq:simulation_snr}.
Figure~\ref{fig:duck_fish_clean_noisy} shows a clean image
(left) and its noisy observation at \(\mathrm{SNR}=0.1\)
(right).

   The observed score map \(S^y\)
in~\eqref{eq:smoothed_process_intro} is computed using
an orthonormal basis of the detection subspace.
To choose the detection subspace, we first construct
a steerable ambient space~\(\mathcal S\) to approximate the duck template \(u_1\),
the fish template \(u_2\), and all their rotations.
Steerability means that \(\mathcal S\) contains every
rotation of each of its functions, as
in~\eqref{def:rotation_invariant}.
We use Fourier-Bessel functions on \(D_{B/2}\)
to construct \(\mathcal S\) in the
form~\eqref{eq:ambient_space}.
We choose \(\mathcal S\) so that the orthogonal projection
of each template onto \(\mathcal S\) retains at least
\(90\%\) of its energy (squared norm).
By steerability, the same bound holds for every rotation
of each template; see
Lemma~\ref{S-lem:proj_rot_appendix} of the supplementary
material.
The resulting ambient space \(\mathcal S\), which is
the direct sum of the angular components \(\mathcal S_j\)
in~\eqref{eq:ambient_space}, has dimension \(123\).

To reduce computational cost, we restrict the candidate
family in~\eqref{eq:candidate_family} to \(18\) nested
candidate detection subspaces of \(\mathcal S\).
Starting from \(\{0\}\), we construct each candidate subspace by successively adding 
angular components~\(\mathcal S_j\) of~\(\mathcal S\).
At each step, for each component \(\mathcal S_j\)
not yet included, we project \(u_1\) and \(u_2\)
onto the sum of the current subspace and \(\mathcal S_j\).
We add the component~\(\mathcal S_j\)  that maximizes the smaller
of the two squared projection norms.
Each addition defines a candidate detection subspace.
The final candidate detection subspace is the
ambient space \(\mathcal S\).

For each candidate detection subspace, we apply
a numerical localization check based on the conditions defining
\(\mathfrak S_{\rm loc}\)
in~\eqref{eq:retained_detection_subspaces},
with \(\delta=20\).
The three candidate detection subspaces that pass the numerical
localization check have dimensions \(7\), \(19\),
and \(31\).

We select the detection subspace using the pSNR
criterion in~\eqref{eq:psnr_V}.
Specifically, for a candidate detection subspace \(\mathcal V\),
let \(P_{\mathcal V}\) denote the orthogonal projection
onto \(\mathcal V\) in \(L^2(D_{B/2})\).
The criterion is maximizing
\[
    \mathrm{pSNR}(\mathcal V)
    =
    \frac{\displaystyle\min_{\ell=1,2}
    \|P_{\mathcal V}u_\ell\|_{L^2(D_{B/2})}^2}
    {\eta_{\mathrm{peak}}(\mathcal V)},
\]
where \(\eta_{\mathrm{peak}}(\mathcal V)\) (see~\eqref{eq:eta_peak_V}) is the
expected maximum noise score over the comparison
cube \(\overline{\mathcal B}(0,r/2)\).
Since \(\mathcal V\) is steerable, rotating either
template does not change its projection norm;
see Lemma~\ref{S-lem:proj_rot_appendix}
of the supplementary material.

We estimate \(\eta_{\mathrm{peak}}(\mathcal V)\)
using 5,000 independent images containing only
Gaussian noise, with the covariance
in~\eqref{eq:simulation_correlated_noise} at \(\sigma=1\).
The 5,000 noise images are independent of the
500 images used to evaluate power and error rates.
For each noise image, we compute the noise score map
in~\eqref{eq:noise_score_map} using an orthonormal
basis of \(\mathcal V\), and find its maximum over
the comparison cube \(\overline{\mathcal B}(0,r/2)\)
from Step~\ref{step:compute_test_values}.
Here \(r=2B+\delta=164\), so the comparison cube
has side length \(82\) pixels.
At noise variance \(\sigma^2\), we define
\(\widehat{\eta}_{\mathrm{peak}}(\mathcal V)\)
as \(\sigma^2\) times the average of the
5,000 maxima.
We therefore define the estimated criterion by
\begin{equation}
    \label{eq:simulation_psnr_estimate}
    \widehat{\mathrm{pSNR}}(\mathcal V)
    =
    \frac{\displaystyle\min_{\ell=1,2}
    \|P_{\mathcal V}u_\ell\|_{L^2(D_{B/2})}^2}
    {\widehat{\eta}_{\mathrm{peak}}(\mathcal V)}.
\end{equation}

\begin{figure}[htbp]
    \centering
    \begin{tabular}{@{}c@{\hspace{0.03\textwidth}}c@{}}
        \includegraphics[width=0.23\textwidth]
        {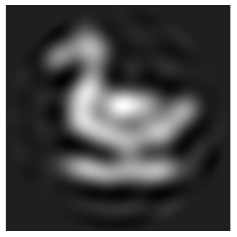}
        &
        \includegraphics[width=0.23\textwidth]
        {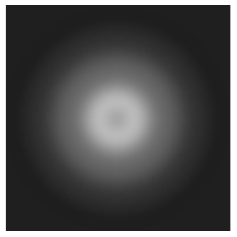}
        \\[3pt]
        \(P_{\mathcal S}u_1\), \(91.83\%\)
        &
        \(P_{\mathcal S^\star}u_1\), \(62.85\%\)
        \\[8pt]
        \includegraphics[width=0.23\textwidth]
        {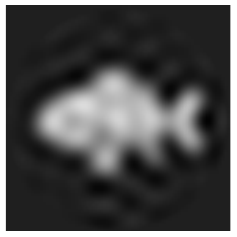}
        &
        \includegraphics[width=0.23\textwidth]
        {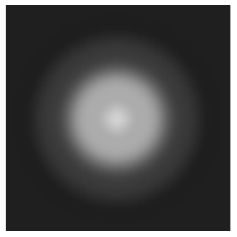}
        \\[3pt]
        \(P_{\mathcal S}u_2\), \(92.83\%\)
        &
        \(P_{\mathcal S^\star}u_2\), \(70.43\%\)
    \end{tabular}
    \caption{Orthogonal projections of the duck template
    \(u_1\) (top) and the fish template \(u_2\) (bottom).
    The left column uses the ambient space \(\mathcal S\)
    of dimension \(123\); the right column uses the selected
    detection subspace \(\mathcal S^\star\) of dimension \(7\).
    The percentages indicate the fraction of each
    template's energy (squared norm) retained by the projection.}
    \label{fig:duck_fish_projections}
\end{figure}

We implement the selection rule in~\eqref{eq:Sstar_def}
by maximizing
\(\widehat{\mathrm{pSNR}}(\mathcal V)\)
of~\eqref{eq:simulation_psnr_estimate}
over the three candidate detection subspaces that pass
the numerical localization check.
The selected detection subspace \(\mathcal S^\star\)
has dimension \(7\). For these templates, the selection rule happens to
choose the radial component \(\mathcal S_0\), of
dimension \(7\) (recall that all candidate detection
subspaces are steerable and yield rotation-invariant
scores).
The factor \(\sigma^2\) in
\(\widehat{\eta}_{\mathrm{peak}}(\mathcal V)\)
is common to all candidate detection subspaces,
so \(\mathcal S^\star\) is the same at every SNR level.
Figure~\ref{fig:duck_fish_projections} shows the
projections of \(u_1\) and \(u_2\) onto the ambient
space \(\mathcal S\) and on the selected detection
subspace \(\mathcal S^\star\).

We now apply the detection procedure of
Subsection~\ref{subsec:algorithm_outline} using
\(\mathcal S^\star\).
We use the Bonferroni and BH procedures in
Step~\ref{step:multiple_testing} at the target
error level \(\alpha=0.05\).

To compute the test values in
Step~\ref{step:compute_test_values}, we estimate
the tail probability \(p(u)\)
in~\eqref{eq:test_value_intro} by Monte Carlo.
We generate 50,000 independent images containing
only Gaussian noise, with the covariance
in~\eqref{eq:simulation_correlated_noise}
at \(\sigma=1\).
The 50,000 noise images are independent of both
the 5,000 noise images used to select
\(\mathcal S^\star\) and the 500 images used
to evaluate power and error rates.
For each noise image, we compute the noise score map
in~\eqref{eq:noise_score_map} using an orthonormal
basis of \(\mathcal S^\star\), and find its maximum
over the comparison cube
\(\overline{\mathcal B}(0,r/2)\).
At each SNR, we multiply the 50,000 maxima by
the corresponding \(\sigma^2\) and use them to obtain a Monte Carlo
estimate \(\widehat p(u)\) of the probability that
the maximum noise score over the comparison cube
exceeds \(u\).
For each candidate location \(t\) selected in
Step~\ref{item:candidate_peaks}, we use
\(\widehat p(S^y(t))\) as the test value.

\begin{figure}[htbp]
    \centering
    \begin{tabular}{@{}c@{\hspace{0.02\textwidth}}c@{}}
        \includegraphics[width=0.46\textwidth]
        {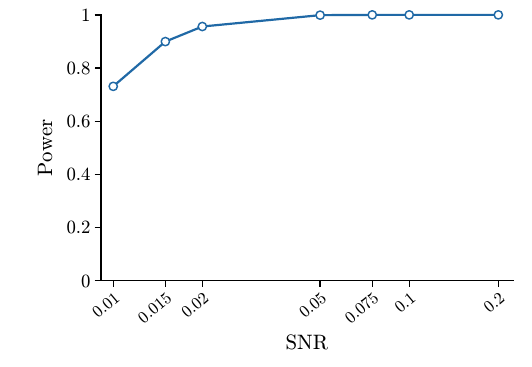}
        &
        \includegraphics[width=0.46\textwidth]
        {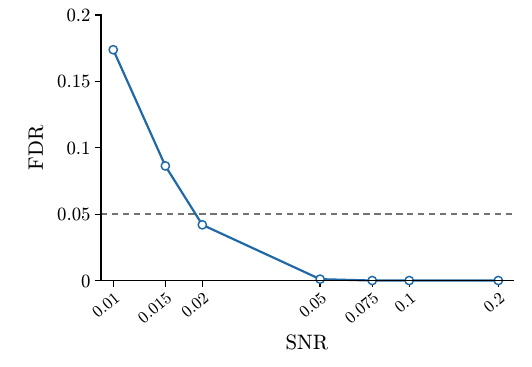}
        \\[6pt]
        \includegraphics[width=0.46\textwidth]
        {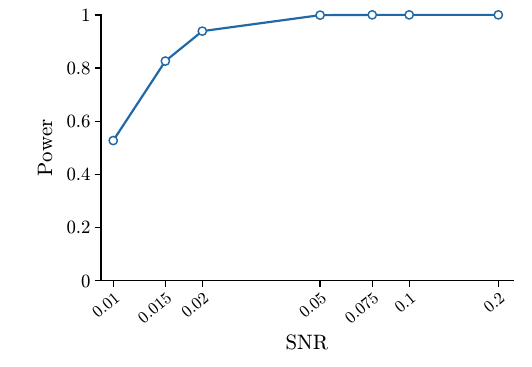}
        &
        \includegraphics[width=0.46\textwidth]
        {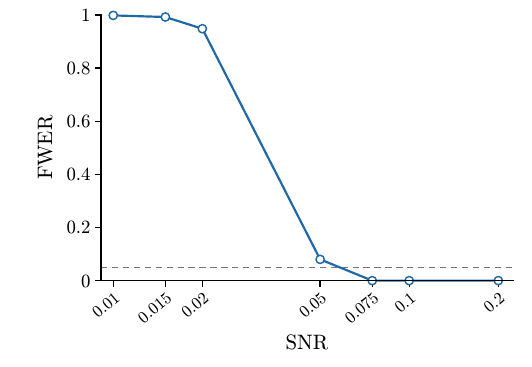}
    \end{tabular}
    \caption{Estimated power and error rates using
    \(\mathcal S^\star\), over 500 independently
    generated images.
    The top row shows the power (left) and FDR (right)
    of BH.
    The bottom row shows the power (left) and FWER
    (right) of Bonferroni.
    Dashed lines in the right panels indicate
the target error level \(\alpha=0.05\).
The horizontal axes use a logarithmic scale.}
    \label{fig:duck_fish_selected_performance}
\end{figure}

We estimate the power and error rates defined in
Subsection~\ref{subsec:performance_criteria},
averaging over both the Gaussian noise and the
randomly generated clean images.
For each SNR, we estimate power by averaging
the fraction of detected object instances over
the 500 images.
We estimate FDR by averaging the fraction of
reported locations that are false over the
same 500 images, assigning zero to images
with no reported locations.
We estimate FWER by the fraction of the 500 images
with at least one false detection.
Figure~\ref{fig:duck_fish_selected_performance} shows
the estimated power and FDR of BH, and the estimated
power and FWER of Bonferroni, at all seven SNR levels.

At the three highest SNR levels, both procedures
detect all object instances in each of the 500 images.
The estimated BH FDR exceeds the target level
\(\alpha=0.05\) at the two lowest SNR levels,
\(0.01\) and \(0.015\).
The estimated Bonferroni FWER exceeds this level
at the four lowest SNR levels,
\(0.01\), \(0.015\), \(0.02\), and \(0.05\).

Next, we compare \(\mathcal S^\star\), of
dimension \(7\), with the ambient space \(\mathcal S\),
of dimension \(123\), and with the candidate detection
subspaces of dimensions \(19\) and \(31\), which pass
the numerical localization check.
The ambient space \(\mathcal S\) does not pass
the numerical localization check, so the guarantees of
Section~\ref{sec:theoretical_foundations} are not
established for it.

\begin{figure}[htbp]
    \centering
    \includegraphics[width=0.96\textwidth]
    {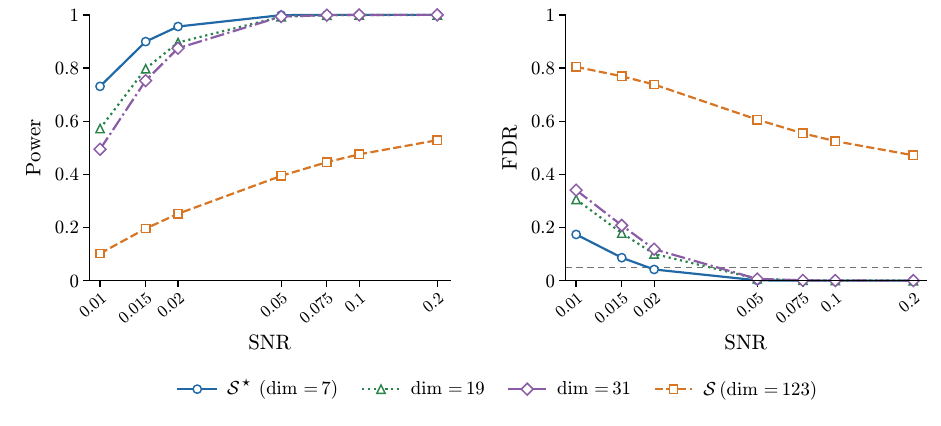}
    \caption{Estimated BH power (left) and FDR (right)
    for the selected detection subspace \(\mathcal S^\star\),
    the candidate detection subspaces of dimensions \(19\) and \(31\)
    that pass the numerical localization check,
    and the ambient space \(\mathcal S\).
    All estimates use the same 500 images at each SNR.
    The dashed line in the right panel marks the target level \(\alpha=0.05\).}
    \label{fig:duck_fish_subspace_comparison}
\end{figure}

For each of the four subspaces, we apply the
detection procedure of
Subsection~\ref{subsec:algorithm_outline} with BH
at \(\alpha=0.05\), using the same 500 observed
images at each SNR.
We estimate the tail probability \(p(u)\) in~\eqref{eq:test_value_intro} separately for each
subspace, using 50,000 noise images and the
Monte Carlo procedure described for
\(\mathcal S^\star\).
Figure~\ref{fig:duck_fish_subspace_comparison}
shows the estimated BH power and FDR for the
four detection subspaces.

\begin{figure}[htbp]
    \centering
    \includegraphics[width=0.98\textwidth]
{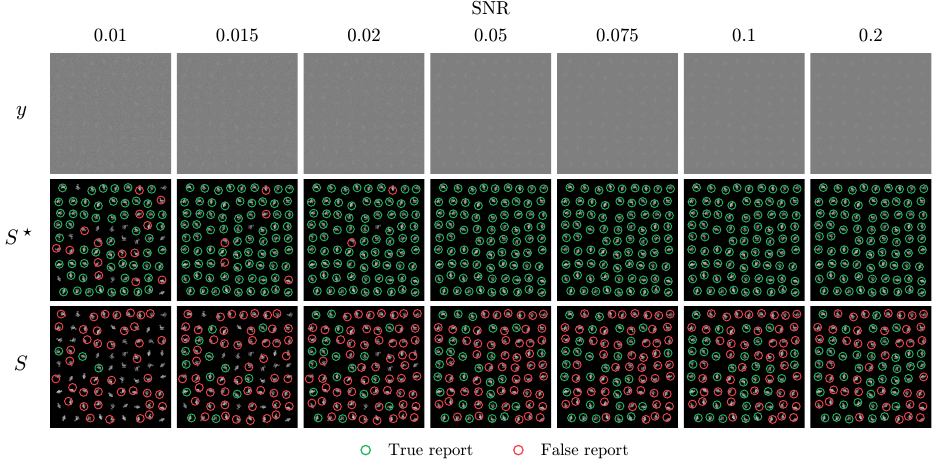}
    \caption{BH detections at \(\alpha=0.05\).
The top row shows the noisy observations \(y\)
at the seven SNR levels indicated.
The middle and bottom rows show the corresponding
detections using \(\mathcal S^\star\) and
\(\mathcal S\), respectively, displayed on the same
clean image \(x\).
Each circle is centered at a reported location.
Green and red indicate true and false detections,
respectively, according to~\eqref{eq:true_noise_regions}.}
    \label{fig:duck_fish_detection_examples}
\end{figure}

Using \(\mathcal S^\star\) gives higher estimated
BH power and lower estimated FDR than using
\(\mathcal S\) at every tested SNR level.
The difference is large: although \(\mathcal S\) retains
more than \(90\%\) of each template's energy
(Figure~\ref{fig:duck_fish_projections}), its estimated
BH FDR stays above \(0.45\) and its estimated BH power
below \(0.55\) at every SNR level. This is consistent
with \(\mathcal S\) failing the numerical localization
check.
Among the three candidate detection subspaces
that pass the numerical localization check,
\(\mathcal S^\star\) attains the highest estimated
BH power at the six lowest SNR levels and the lowest
estimated FDR at the three lowest SNR levels.
Figure~\ref{fig:duck_fish_detection_examples}
illustrates the BH reports obtained using
\(\mathcal S^\star\) and \(\mathcal S\) for one
of the 500 images, at all seven SNR levels.
From \(\mathrm{SNR}=0.05\) upward, all reports obtained using
\(\mathcal S^\star\) are true detections, whereas many
false reports obtained using \(\mathcal S\) lie next to
object instances but outside their detection cubes.

\section{Discussion}
\label{sec:discussion}

We give a procedure for detecting translated and
rotated copies of known templates in stationary
Gaussian noise with known law, using a steerable
detection subspace.
Under the sufficient conditions established in
Section~\ref{sec:theoretical_foundations}
for a fixed detection subspace,
the simultaneous recovery probability tends to one at
the Bonferroni and BH thresholds, BH has power tending
to one under a weaker condition, and Bonferroni and BH
asymptotically control FWER and FDR, respectively.
When \(N/M_L\to0\),
Proposition~\ref{prop:retained_energy_lower_bounds}
shows that the growth rates that these conditions impose
on the smallest squared projection norm cannot be
improved in general.
These conditions are stated in terms of the ratio
\(a_L^2/\eta_{\mathrm{peak}}\), and both quantities
depend on the detection subspace: a larger subspace can
give the templates larger projection norms but also
raise the expected maximum noise score.
We therefore use the guarantees themselves to choose
the detection subspace. Among the subspaces in
\(\mathfrak S_{\rm loc}\), which satisfy the
localization conditions
in~\eqref{eq:retained_detection_subspaces}, we maximize
the pSNR criterion in~\eqref{eq:psnr_V}, given by the ratio of
the smallest squared template projection norm to the
expected maximum noise score over the comparison cube.
By Corollary~\ref{cor:asymptotic_selected}, this choice
gives the weakest sufficient conditions on the smallest
squared amplitude \(b_L^2\)
for BH power and simultaneous recovery, so the detection
subspace is chosen for detection rather than for
approximating the templates.
In the simulations of Section~\ref{sec:simulations},
the selected subspace \(\mathcal S^\star\) has smaller
squared template projection norms than the ambient space
\(\mathcal S\) (Figure~\ref{fig:duck_fish_projections}),
but gives higher
estimated BH power and lower estimated FDR at every tested
SNR level, and the highest estimated BH power among the
candidate detection subspaces (that pass the numerical
localization check) at the six lowest SNR levels
(Figure~\ref{fig:duck_fish_subspace_comparison}).

The choice of the detection subspace does not arise
in~\cite{wakninEtAl2026object},
where the detection subspace is given and each object
instance, translated to the origin, belongs to it
(the linear subspace model).
However, the score map, the selection of candidate locations,
and the multiple-testing step in the current work follow~\cite{wakninEtAl2026object}.
For the linear subspace model, our results improve
on~\cite{wakninEtAl2026object} in two respects.
First, \cite{wakninEtAl2026object} proposed calibrating
the test values from the maximum noise score over the
comparison cube without rigorous performance proofs.
Under Assumption~\ref{ass:noise_reg_main},
Theorem~\ref{thm:continuity_test} shows that this maximum
has an absolutely continuous distribution, so
\(\Prob(Z_{\max}\ge u_c)=p(u_c)=c\)
(Corollary~\ref{S-cor:generalized_inverse_attainment}),
and Section~\ref{sec:theoretical_foundations} gives the
missing proofs.
Second, our guarantees apply without requiring the
detection subspace to be steerable, and the sufficient
condition for simultaneous recovery
in~\eqref{eq:object_strength_condition} requires
\(a_L^2=\min_{1\le i\le N}\|x_i\|_{L^2(\mathbb R^d)}^2\)
to eventually exceed a multiple of the logarithm of the
image volume (Corollary~\ref{S-cor:exact_linear_subspace}
of the supplementary material, which also covers \(d=1\)),
whereas the sufficient condition
in~\cite{wakninEtAl2026object} requires \(a_L^2\) to grow
faster than the image volume.

Within the setting of this paper, with known templates
and noise law, questions remain about performance for a
fixed side length \(L\) and the accuracy of the numerical
implementation.
The simulations in Section~\ref{sec:simulations}
show that Bonferroni and BH detect all object instances
in every image at the three highest SNR levels, but
the estimated BH FDR exceeds \(\alpha\) at the two
lowest SNR levels and the estimated Bonferroni FWER
at the four lowest.
These exceedances do not contradict
Section~\ref{sec:theoretical_foundations}:
Proposition~\ref{prop:asymptotic_error_power}
bounds \(\FWER(u_{\rm Bon})\) by
\(\alpha+\varepsilon_L\) rather than \(\alpha\) for
every \(L\), and the FDR guarantee in
Proposition~\ref{prop:bh_fdr_control} is asymptotic.
Sharper bounds for a fixed \(L\) than those in
Propositions~\ref{prop:asymptotic_error_power}
and~\ref{prop:bh_density_strength} could describe more
precisely how power and error rates depend on
\(a_L^2/\eta_{\mathrm{peak}}\).
Guarantees for the numerical implementation also
require accounting for spatial discretization,
the accuracy of the numerical localization check,
and uncertainty in the Monte Carlo estimates of
\(\eta_{\mathrm{peak}}(\mathcal V)\) and \(p(u)\).

When the noise is Gaussian but its law is unknown,
the expected maximum noise score
\(\eta_{\mathrm{peak}}(\mathcal V)\)
in~\eqref{eq:eta_peak_V}, used by the pSNR criterion,
and the tail probability \(p(u)\)
in~\eqref{eq:test_value_intro}, used for threshold
calibration, can no longer be estimated from simulated
noise images as in Section~\ref{sec:simulations}.
The question is how to estimate these quantities
from images containing object instances accurately
enough to preserve the power and error-control
guarantees, particularly at the levels \(\alpha/M_L\)
and \(k\alpha/M_L\) used by Bonferroni and BH.

When clean templates are unavailable, the detection
subspace must also be learned from observed images.
The ambient space in
Proposition~\ref{prop:ambient_uniform_coverage}, the
localization conditions
in~\eqref{eq:retained_detection_subspaces}, and the
pSNR criterion in~\eqref{eq:psnr_V} are all computed
from the known templates.
An open problem is to learn a steerable detection
subspace while ensuring that
\(a_L^2/\eta_{\mathrm{peak}}\) is large enough for the
conditions of Section~\ref{sec:theoretical_foundations}
and that the localization property in
Assumption~\ref{def:localization_property} holds.
The guarantees of Section~\ref{sec:theoretical_foundations}
take the detection subspace and \(p(u)\) to be fixed
before the image is observed.
If the observed image \(y\) is also used to learn
the subspace or to estimate
\(\eta_{\mathrm{peak}}(\mathcal V)\) and \(p(u)\),
the analysis must account for the dependence of the
learned subspace and of these estimated quantities on \(y\).

Finally, the guarantees of
Section~\ref{sec:theoretical_foundations} assume
Gaussian noise. In particular, Gaussianity is used to
prove that the maximum noise score over the comparison
cube has an absolutely continuous distribution
(Theorem~\ref{thm:continuity_test}), to obtain the
high-probability upper bound
in~\eqref{eq:search_domain_noise_bound} on the maximum
noise score over the search domain, and to
show that the variance bound in
Assumption~\ref{ass:noise_exceedance_variance} holds
if \(\Cov(z(0),z(h))=0\) whenever \(\|h\|_\infty\)
exceeds a fixed constant
(Lemma~\ref{S-lem:noise_exceedance_variance_finite_range}).
The first assertion of
Proposition~\ref{prop:bh_fdr_control} uses the variance
bound, together with \(N\to\infty\), for BH FDR control.
Extending the analysis to non-Gaussian noise therefore
requires non-Gaussian replacements for each of these
results.
The pSNR criterion in~\eqref{eq:psnr_V} also relies on
the power bound in~\eqref{eq:bh_average_power_rate}
and on the condition in~\eqref{eq:object_strength_condition},
both stated in terms of \(\eta_{\mathrm{peak}}\).
Extending the pSNR criterion to non-Gaussian noise
would require analogues of~\eqref{eq:bh_average_power_rate}
and~\eqref{eq:object_strength_condition}
expressed through \(\eta_{\mathrm{peak}}\).

\begin{acks}[Acknowledgments]
The authors thank Amitay Eldar for helpful and valuable discussions.
\end{acks}

\section*{Statement on AI use}
The authors used ChatGPT and Codex (OpenAI) to assist
with wording and exposition and to obtain feedback
on the clarity and consistency of the manuscript,
mathematical arguments, notation, and references.
These tools were also used for developing, debugging,
and running code for the numerical experiments and
the preparation of figures and tables. The authors
reviewed the suggestions and outputs and take full
responsibility for the final content of the article
and its supplementary material.

\section*{Supplementary material}
\addcontentsline{toc}{section}{Supplementary material}
The supplement contains the proofs of the results in
		Sections~\ref{sec:theoretical_foundations} and~\ref{sec:repr_framework},
		the proof of the lower bound in
		Proposition~\ref{prop:retained_energy_lower_bounds}, and sufficient conditions for the noise assumptions.

\setcounter{section}{0}
\setcounter{equation}{0}
\renewcommand{\thesection}{S\arabic{section}}
\renewcommand{\theequation}{S\arabic{equation}}
\renewcommand{\theHsection}{supp.\arabic{section}}
\renewcommand{\theHsubsection}{\theHsection.\arabic{subsection}}
\renewcommand{\theHsubsubsection}{\theHsubsection.\arabic{subsubsection}}
\renewcommand{\theHequation}{supp.\arabic{equation}}
\renewcommand{\theHtheorem}{\theHsection.\arabic{theorem}}

	\section{Absolute continuity of the maximum noise score}
	\label{S-app:continuity_noise_test}
	
	\begin{lemma}[Continuity of the projected noise field]
		\label{S-lem:projected_noise_continuity}
		Let \(z\) be a jointly measurable centered stationary Gaussian field on
		\(\R^d\), let \(\psi_1,\dots,\psi_M\in L^2(\R^d)\) be real-valued
		and compactly supported, and define
		\[
		X_j(t):=(z*\tilde\psi_j)(t),
		\qquad
		\mathbf X(t):=(X_1(t),\dots,X_M(t))^\top,
		\qquad
		\tilde\psi_j(s):=\psi_j(-s).
		\]
		Then \(\mathbf X\) is a centered Gaussian field and admits a
		modification with continuous sample paths on \(\R^d\).
	\end{lemma}
	
	\begin{proof}
		For \(D_n=[-n,n]^d\), Tonelli's theorem
		and stationarity give
		\[
		\E\int_{D_n}|z(u)|^2\,\dd u
		=|D_n|\,\E|z(0)|^2<\infty.
		\]
		Thus the measurable event
		\(\Omega_0:=\bigcap_{n\ge1}\{\int_{D_n}|z(u)|^2\,\dd u<\infty\}\)
		has probability one, and \(z(\cdot,\omega)\in L^2_{\mathrm{loc}}(\R^d)\)
		for every \(\omega\in\Omega_0\). On this event, define
		\[
		\widehat X_j(t,\omega)
		:=\int_{\R^d}z(u,\omega)\psi_j(u-t)\,\dd u,
		\qquad j=1,\dots,M,
		\]
		and set \(\widehat X_j(t,\omega)=0\) outside \(\Omega_0\).
		By the Cauchy-Schwarz inequality and compact support, these integrals
		are finite for every \(j\) and \(t\) on \(\Omega_0\). Integrating the
		jointly measurable function
		\(\mathbf1_{\Omega_0}(\omega)z(u,\omega)\psi_j(u-t)\) over \(u\)
		gives joint measurability of \(\widehat X_j\). Since
		\(\tilde\psi_j(t-u)=\psi_j(u-t)\), the vector
		\(\widehat{\mathbf X}:=(\widehat X_1,\ldots,\widehat X_M)^\top\)
		agrees with \(\mathbf X\) almost surely at every fixed \(t\) and is
		therefore a modification of \(\mathbf X\).
		Fix \(\omega\in\Omega_0\) and \(t_0\in\R^d\). For all \(j\) and all
		\(t\) near \(t_0\), the supports of \(\psi_j(\cdot-t)\) are contained
		in a common compact set \(C\). By the Cauchy-Schwarz inequality,
		\[
		|\widehat X_j(t,\omega)-\widehat X_j(t_0,\omega)|
		\le \|z(\cdot,\omega)\mathbf1_C\|_{L^2}
		\|\psi_j(\cdot-t)-\psi_j(\cdot-t_0)\|_{L^2}.
		\]
		The last factor tends to zero as \(t\to t_0\) by continuity of
		translations in \(L^2(\R^d)\). Thus \(\widehat{\mathbf X}\) has
		continuous sample paths.
		To prove Gaussianity, for real compactly supported \(g\in L^2(\R^d)\),
		which also belongs to \(L^1(\R^d)\), let
		\(I_g:=\int_{\R^d}g(u)z(u)\,\dd u\) on \(\Omega_0\), and set
		\(I_g=0\) outside this event. Write \(\sigma_z^2=\Var(z(0))\).
		Minkowski's inequality and stationarity give
		\[
		\|I_g\|_{L^2(\Omega)}
		\le\int_{\R^d}|g(u)|\,\|z(u)\|_{L^2(\Omega)}\,\dd u
		=\sigma_z\|g\|_{L^1(\R^d)}<\infty.
		\]
		Let
		\[
		H:=\overline{\spn\{z(u):u\in\R^d\}}^{\,L^2(\Omega)}.
		\]
		Every element of \(H\) is centered Gaussian, since it is an
		\(L^2(\Omega)\)-limit of finite linear combinations of the jointly
		Gaussian variables \(z(u)\), and both distributions and variances
		converge under such a limit. For \(Y\in H^\perp\), Cauchy-Schwarz inequality gives
		\[
		\int_{\R^d}|g(u)|\,\E|Yz(u)|\,\dd u
		\le\sigma_z\|g\|_{L^1(\R^d)}\|Y\|_{L^2(\Omega)}<\infty.
		\]
		Fubini's theorem therefore yields
		\[
		\E[YI_g]=\int_{\R^d}g(u)\E[Yz(u)]\,\dd u=0,
		\]
		because \(z(u)\in H\) for every \(u\). Hence
		\(I_g\in(H^\perp)^\perp=H\), so \(I_g\) is centered Gaussian.
		Finally, every finite linear combination
		\(\sum_{k=1}^n a_k\widehat X_{j_k}(t_k)\) equals \(I_g\) with
		\(g=\sum_{k=1}^n a_k\psi_{j_k}(\cdot-t_k)\), which is compactly
		supported and belongs to \(L^2(\R^d)\).
		Thus \(\widehat{\mathbf X}\), and therefore its modification
		\(\mathbf X\), is a centered Gaussian field.
	\end{proof}
	
	Under Assumption~\ref{ass:noise_reg_main}, use the continuous
	modification from Lemma~\ref{S-lem:projected_noise_continuity} and
	denote it by \(\mathbf X\). Define
	\begin{equation}\label{S-eq:noise_maximum_space}
		K:=\overline{\mathcal B}(0,r/2),
		\qquad E:=C(K;\R^M),
		\qquad \|f\|_E:=\max_{t\in K}\|f(t)\|,
	\end{equation}
	and set
	\begin{equation}\label{S-eq:noise_maximum_random_element}
		G:=\mathbf X|_K,
		\qquad Q:=\|G\|_E,
		\qquad Z_{\max}=Q^2.
	\end{equation}
	To see that \(G\) is a centered Gaussian random element of \(E\),
	take regular grids of \(K\) with spacing tending to zero. On the
	\(n\)th grid, let \(G_n\) be the function obtained by interpolating
	the values of the projected noise field linearly in each spatial
	coordinate. Each \(G_n\) depends continuously on finitely many field
	values and is therefore measurable as a map into \(E\). Since the
	sample paths are continuous on the compact set \(K\), \(G_n\) converges
	uniformly to \(G\), so \(G\) is also measurable. For every continuous
	linear functional \(\ell\) on \(E\), the variable \(\ell(G_n)\) is a
	finite linear combination of the variables \(X_j(t)\) at grid points.
	By Lemma~\ref{S-lem:projected_noise_continuity}, these variables are
	jointly Gaussian with mean zero. Thus \(\ell(G_n)\) is centered
	Gaussian, and so is its almost sure limit \(\ell(G)\).
	Since \(E\) is a separable Banach space, the distribution of \(G\)
	on \(E\) is a Radon probability measure. Continuity of the sample paths
	on \(K\) also gives \(Q=\|G\|_E<\infty\) almost surely.
	Recall that a seminorm on \(\R^M\) is a function
	\(\nu:\R^M\to[0,\infty)\) satisfying
	\[
	\nu(cx)=|c|\nu(x),
	\qquad
	\nu(x+y)\le \nu(x)+\nu(y)
	\]
	for all \(c\in\R\) and \(x,y\in\R^M\).
	Its kernel \(\ker\nu=\{x:\nu(x)=0\}\) may contain nonzero vectors.
	We call \(\nu\) nonzero if \(\nu(x)>0\) for some \(x\in\R^M\).
	
	\begin{corollary}[Maxima of seminorm scores]
		\label{S-cor:seminorm_score_continuity}
		Under Assumption~\ref{ass:noise_reg_main}, let
		\(K=\overline{\mathcal B}(0,r/2)\), let \(\nu\) be a nonzero
		seminorm on \(\R^M\), and let \(\beta>0\). Then
		\[
		Z_{\nu,\beta}
		:=\max_{t\in K}\nu\bigl(\mathbf X(t)\bigr)^\beta
		\]
		has an absolutely continuous distribution. In particular, this
		applies to
		\[
		\max_{t\in K}\sum_{j=1}^M w_j|X_j(t)|^q,
		\qquad 1\le q<\infty,\quad w_1,\dots,w_M>0,
		\]
		and to
		\[
		\max_{t\in K}\mathbf X(t)^\top H\mathbf X(t),
		\qquad H\in\R^{M\times M},\quad H\succeq0,\quad H\ne0.
		\]
	\end{corollary}
	
	\begin{proof}
		Let \(e_1,\dots,e_M\) be the standard basis of \(\R^M\).
		The seminorm properties and the Cauchy-Schwarz inequality give
		\[
		\nu(x)\le\sum_{j=1}^M |x_j|\nu(e_j)
		\le C_\nu\|x\|,
		\qquad
		C_\nu:=\left(\sum_{j=1}^M\nu(e_j)^2\right)^{1/2}.
		\]
		Since \(\nu\) is nonzero, \(C_\nu>0\). Also,
		\(|\nu(x)-\nu(y)|\le\nu(x-y)\le C_\nu\|x-y\|\), so \(\nu\)
		is continuous. Thus, for each \(f\in E\), the function
		\(t\mapsto\nu(f(t))\) attains its maximum on \(K\). Define
		\[
		\Phi_\nu(f):=\max_{t\in K}\nu(f(t)),
		\qquad Q_\nu:=\Phi_\nu(G).
		\]
		The triangle inequality and absolute homogeneity for \(\nu\)
		give the same properties for \(\Phi_\nu\). Moreover,
		\[
		|\Phi_\nu(f)-\Phi_\nu(g)|
		\le \max_{t\in K}\nu(f(t)-g(t))
		\le C_\nu\|f-g\|_E.
		\]
		Hence \(\Phi_\nu\) is a continuous seminorm on \(E\), and
		\(Q_\nu\le C_\nu Q<\infty\) almost surely.
		
		Let
		\[
		b_\nu:=\inf\{a>0:\Prob(Q_\nu\le a)>0\}.
		\]
		By Theorem~1.2 of \cite{hoffmannJorgensenSheppDudley1979},
		the distribution of \(Q_\nu\) is absolutely continuous on
		\((b_\nu,\infty)\) and may have an atom at \(b_\nu\).
		The authors explain on pp.~319-320 how the theorem applies
		to seminorms under a centered Gaussian Radon measure.
		It therefore suffices to show that \(b_\nu=0\) and
		\(\Prob(Q_\nu=0)=0\). For the first claim, \(Q_\nu\le C_\nu Q\),
		so it suffices to prove that \(\Prob(Q<a)>0\) for every \(a>0\).
		Fix \(a>0\) and let \(G'\) be an independent copy of \(G\).
		By separability, \(E\) has a countable cover by open balls of
		radius \(a/\sqrt{2}\). At least one ball \(B\) in this cover
		satisfies \(\Prob(G\in B)>0\). If \(G\) and \(G'\) both lie in
		\(B\), then \(\|(G-G')/\sqrt{2}\|_E<a\).
		The fields \(G\) and \((G-G')/\sqrt{2}\) are centered Gaussian
		with the same covariance, so their finite-dimensional
		distributions agree. Since both have continuous paths, their
		values on a countable dense subset of \(K\) determine their
		distributions on \(E\). Therefore, by independence,
		\begin{equation}\label{S-eq:noise_norm_small_ball}
			\Prob(Q<a)
			=\Prob\!\left(\left\|(G-G')/\sqrt{2}\right\|_E<a\right)
			\ge \Prob(G\in B)^2>0.
		\end{equation}
		Applying \eqref{S-eq:noise_norm_small_ball} with \(a/C_\nu\) gives
		\(\Prob(Q_\nu<a)\ge\Prob(Q<a/C_\nu)>0\) for every \(a>0\).
		This proves \(b_\nu=0\).
		It remains to show that \(\Prob(Q_\nu=0)=0\). Since \(0\in K\),
		the event \(Q_\nu=0\) implies \(\nu(\mathbf X(0))=0\).
		The kernel of \(\nu\) is a proper linear subspace of \(\R^M\)
		because \(\nu\) is nonzero. By
		Assumption~\ref{ass:noise_reg_main}, the Gaussian vector
		\(\mathbf X(0)\) has positive definite covariance and hence a
		density on \(\R^M\). It follows that
		\[
		\Prob(Q_\nu=0)
		\le \Prob\{\mathbf X(0)\in\ker\nu\}=0.
		\]
		Together with absolute continuity on \((0,\infty)\), this proves
		that \(Q_\nu\) has an absolutely continuous distribution.
		Finally, \(Z_{\nu,\beta}=Q_\nu^\beta\). If \(h_\nu\) denotes a
		density of \(Q_\nu\), the change of variables \(u=s^\beta\) on
		\((0,\infty)\) gives the density
		\[
		\frac{1}{\beta}u^{1/\beta-1}h_\nu(u^{1/\beta}),
		\qquad u>0,
		\]
		for \(Z_{\nu,\beta}\), and \(\Prob(Z_{\nu,\beta}=0)=0\).
		The first example follows by taking
		\(\nu(x)=(\sum_{j=1}^M w_j|x_j|^q)^{1/q}\) and \(\beta=q\).
		For the second, take \(\nu(x)=\|H^{1/2}x\|\) and \(\beta=2\).
	\end{proof}
	
	\begin{proof}[Proof of Theorem~\ref{thm:continuity_test}]
		Apply Corollary~\ref{S-cor:seminorm_score_continuity} with
		\(\nu(x)=\|x\|\) and \(\beta=2\), for which
		\(Z_{\nu,\beta}=Z_{\max}\).
	\end{proof}
	
	\begin{corollary}[Threshold attainment]
		\label{S-cor:generalized_inverse_attainment}
		Under Assumption~\ref{ass:noise_reg_main}, for \(c\in(0,1)\), let
		\[
		p(u):=\Prob(Z_{\max}>u),
		\qquad
		u_c:=\inf\{u\in\R:p(u)\le c\}.
		\]
		Then
		\[
		0<u_c<\infty,
		\qquad
		p(u_c)=\Prob(Z_{\max}\ge u_c)=c,
		\]
		and, for every \(s\in\R\),
		\[
		p(s)\le c
		\quad\Longleftrightarrow\quad
		s\ge u_c.
		\]
	\end{corollary}
	
	\begin{proof}
		By Theorem~\ref{thm:continuity_test}, \(p\) is continuous and
		nonincreasing, with \(p(0)=1\) and \(p(u)\to0\) as \(u\to\infty\).
		Hence \(A_c:=\{u:p(u)\le c\}\) is a nonempty closed upper ray whose left
		endpoint \(u_c\) lies in \((0,\infty)\). If \(p(u_c)<c\), continuity
		would give some \(v<u_c\) with \(p(v)\le c\), contradicting the
		definition of \(u_c\). Thus \(p(u_c)=c\), and
		\(A_c=[u_c,\infty)\) gives the stated equivalence. Since the law of
		\(Z_{\max}\) has no atoms, \(\Prob(Z_{\max}\ge u_c)=p(u_c)=c\).
	\end{proof}

	\section{Proofs for Section~\ref{sec:theoretical_foundations}}
	\label{S-app:proofs_theoretical_foundations}
	
	\subsection{Preliminaries}
	
	Throughout this section, the orthonormal family
	\(\{\psi_j\}_{j=1}^M\) from
	Subsection~\ref{subsec:algorithm_outline} is fixed and does not depend on \(L\).
Write \(\mathcal V=\spn\{\psi_1,\dots,\psi_M\}\), and let
\(P_{\mathcal V}\) denote the orthogonal projection onto
\(\mathcal V\) in \(L^2(D_{B/2})\).
	Recall the definitions of the coefficient vectors \(a_i\) and their minimum norm \(a_L\)~\eqref{eq:x}, \eqref{eq:aL_def}.
	We set \(r=2B+\delta\), with \(0<\delta\le2B\).
	For \(f\in\{x,y,z\}\), set
	\[
	\mathbf X^f(t):=\bigl((f*\tilde\psi_j)(t)\bigr)_{j=1}^M,
	\qquad
	S^f(t):=\|\mathbf X^f(t)\|^2,
	\qquad
	\tilde\psi_j(s):=\psi_j(-s).
	\]
    Write \(X_j^f\) for the \(j\)th component of
\(\mathbf X^f\).
In particular, \(\mathbf X^z=\mathbf X\) and
\(\Sigma_{\mathbf X}=\Cov(\mathbf X^z(0))\),
as in Assumption~\ref{ass:noise_reg_main}.
	The model \(y=x+z\) in~\eqref{eq:the_model} and linearity
	of convolution give \(\mathbf X^y=\mathbf X^x+\mathbf X^z\).
	Then, by the reverse triangle inequality
	\begin{equation}\label{S-eq:score_triangle_bounds}
		\bigl|\sqrt{S^y(t)}-\sqrt{S^x(t)}\bigr|
		\le\sqrt{S^z(t)}.
	\end{equation}
The tail probability in~\eqref{eq:test_value_intro} is
\[
	p(u)=\Prob\!\left\{
	\max_{t\in\overline{\mathcal B}(0,r/2)}S^z(t)>u
	\right\}.
	\]
	Recall from~\eqref{eq:Qm}-\eqref{eq:ML_def} that the grid
	consists of the half-open cubes
	\[
	Q_m=\frac r2m+[0,r/2)^d,
	\qquad m\in\mathbb Z^d,
	\]
	and that
	\[
	\mathcal I_L
	=
	\left\{m\in\mathbb Z^d:
	Q_m\cap\overline{\mathcal B}(0,L)\neq\emptyset\right\},
	\qquad
	M_L=|\mathcal I_L|.
	\]
	Each \(\overline Q_m\) is a translate of the comparison cube
	\(\overline{\mathcal B}(0,r/2)\).
	Let \(n_L\) be the number of grid intervals of length \(r/2\)
	that intersect \([-L/2,L/2]\). Then
	\[
	\frac{2L}{r}\le n_L\le\frac{2L}{r}+2,
	\qquad M_L=n_L^d,
	\]
	which yields
	\begin{equation}\label{S-eq:grid_count_asymptotic}
		\frac{M_L}{L^d}\longrightarrow\left(\frac2r\right)^d.
	\end{equation} 
	
	The following lemma describes where an object instance can
	contribute to the score and how selecting a candidate affects
	the detection cubes.
	
	\begin{lemma}[Geometry of candidate selection]
		\label{S-lem:support_and_exclusion}
		Suppose that Assumption~\ref{assumptions} holds.
		Let \(T\) be the candidate set from
		Step~\ref{item:candidate_peaks}, and recall that
		\(r=2B+\delta\), with \(0<\delta\le2B\).
		Then the following statements hold.
		
		\begin{enumerate}
			\item[(i)]
			For every \(i=1,\dots,N\) and \(m=1,\dots,M\),
			\begin{equation}\label{S-eq:candidate_support_cutoff}
				(x_i*\tilde\psi_m)(t)=0
				\quad\text{if }\|t-\tau_i\|_\infty\ge B.
			\end{equation}
			For \(q\in\mathcal B(\tau_i,\delta)\),
			\((x*\tilde\psi_m)(q)=(x_i*\tilde\psi_m)(q)\).
			At the center,
			\begin{equation}\label{S-eq:candidate_center_score}
				S^x(\tau_i)=\|a_i\|^2.
			\end{equation}
			
			\item[(ii)]
			The detection cubes \(\mathcal B(\tau_i,\delta)\)
			are pairwise disjoint.
			For every \(s_i\in\mathcal B(\tau_i,\delta)\),
			\begin{equation}\label{S-eq:candidate_geometric_facts}
				\mathcal B(\tau_i,2B)\subset\mathcal B(s_i,r),
				\qquad
				\mathcal B(s_i,r)\cap\mathcal B(\tau_j,\delta)
				=\emptyset
				\quad(j\ne i).
			\end{equation}
			For every \(t\in\R^d\),
			\begin{equation}\label{S-eq:candidate_local_exclusion}
				\|t-\tau_i\|_\infty<B
				\quad\Longrightarrow\quad
				\mathcal B(\tau_i,\delta)\subset\mathcal B(t,r).
			\end{equation}
			
			\item[(iii)]
			Each detection cube contains at most one candidate.
			If \(s_i\in T\cap\mathcal B(\tau_i,\delta)\), then
			every \(t\in T\setminus\{s_i\}\) satisfies
			\begin{equation}\label{S-eq:candidate_remaining_distance}
				\|t-\tau_i\|_\infty>B.
			\end{equation}
		\end{enumerate}
	\end{lemma}

	\begin{proof}
		\textup{(i)}
		We have \((x_i*\tilde\psi_m)(t)=\int_{\mathbb{R}^d} x_i(u)\psi_m(u-t)\,\dd u\), where
		\(x_i\) is supported in \(\tau_i+D_{B/2}\) and \(\psi_m(\cdot-t)\) in
		\(t+D_{B/2}\). If \(\|t-\tau_i\|_\infty\ge B\), then
		\(\|t-\tau_i\|\ge B\), so these supports (balls) meet in at most one point and
		the integral vanishes which gives~\eqref{S-eq:candidate_support_cutoff}.
		Let \(q\in\mathcal B(\tau_i,\delta)\) and \(j\ne i\). By the separation
		condition~\eqref{eq:separation},
		\[
		\|q-\tau_j\|_\infty
		\ge\|\tau_i-\tau_j\|_\infty-\|q-\tau_i\|_\infty
		>B+\delta,
		\]
		so \((x_j*\tilde\psi_m)(q)=0\) by~\eqref{S-eq:candidate_support_cutoff}.
		Hence \((x*\tilde\psi_m)(q)=(x_i*\tilde\psi_m)(q)\). At \(q=\tau_i\), by the substitution \(u=\tau_i+s\), the model
		\eqref{eq:template_instance_model} and the support of \(\psi_m\) give
		\[
		(x_i*\tilde\psi_m)(\tau_i)
		=\int_{\R^d}x_i(u)\,\psi_m(u-\tau_i)\,\dd u
		=\int_{D_{B/2}}b_i(R_iu_{\ell_i})(s)\,\psi_m(s)\,\dd s
		=a_{i,m},
		\]
		where the last equality is definition~\eqref{eq:x}. Therefore
		\[
		S^x(\tau_i)=\sum_{m=1}^M(x*\tilde\psi_m)^2(\tau_i)
		=\sum_{m=1}^M a_{i,m}^2=\|a_i\|^2,
		\]
		which proves~\eqref{S-eq:candidate_center_score}.
		
		\textup{(ii)}
		For \(i\ne j\), the separation condition gives
		\(\|\tau_i-\tau_j\|_\infty>B+3\delta/2>\delta\), so the detection cubes
		are pairwise disjoint. If \(v\in\mathcal B(\tau_i,2B)\) and
		\(s\in\mathcal B(\tau_i,\delta)\), then
		\[
		\|v-s\|_\infty\le\|v-\tau_i\|_\infty+\|\tau_i-s\|_\infty
		<B+\frac\delta2=\frac r2 .
		\]
		Taking \(s=s_i\) gives the first inclusion
		in~\eqref{S-eq:candidate_geometric_facts}. Since
		\(\|t-\tau_i\|_\infty<B\) means \(t\in\mathcal B(\tau_i,2B)\), taking
		\(v=t\) gives~\eqref{S-eq:candidate_local_exclusion}.
		For \(j\ne i\) and \(q_j\in\mathcal B(\tau_j,\delta)\),
		\[
		\|q_j-s_i\|_\infty
		\ge\|\tau_j-\tau_i\|_\infty-\|q_j-\tau_j\|_\infty-\|s_i-\tau_i\|_\infty
		>\Bigl(B+\frac{3\delta}2\Bigr)-\delta=\frac r2,
		\]
		so \(q_j\notin\mathcal B(s_i,r)\), which is the second assertion
		in~\eqref{S-eq:candidate_geometric_facts}.
		
		\textup{(iii)}
		Distinct candidates have \(\ell_\infty\)-distance at least \(r/2\)
		from each other; see Step~\ref{item:candidate_peaks}. Let
		\(s_i\in T\cap\mathcal B(\tau_i,\delta)\) and
		\(t\in T\setminus\{s_i\}\). Then
		\[
		\|t-\tau_i\|_\infty
		\ge\|t-s_i\|_\infty-\|s_i-\tau_i\|_\infty
		>\frac r2-\frac\delta2=B,
		\]
		which is~\eqref{S-eq:candidate_remaining_distance}. Since
		\(B\ge\delta/2\), the point \(t\) does not lie in
		\(\mathcal B(\tau_i,\delta)\), so each detection cube contains at most
		one candidate.
	\end{proof}
	
	Recall the true and false detection regions
	\(\mathbb B_1^\delta,\mathbb B_0^\delta\) from
	\eqref{eq:true_noise_regions}, and the set of reported locations \(T(u)\) from
	Step~\ref{step:multiple_testing}.
	As in Subsection~\ref{subsec:performance_criteria}, let
	\[
	W(u)=\#(T(u)\cap\mathbb B_1^\delta),
	\qquad
	V(u)=\#(T(u)\cap\mathbb B_0^\delta)
	\]
	be the numbers of true and false detections, respectively.
	For \(u\in\R\cup\{+\infty\}\), let \(D(u)\) count the object
	instances whose detection cubes contain no reported location:
	\begin{equation}\label{S-eq:bh_missed_neighborhoods}
		D(u):=\sum_{i=1}^N
		\mathbf1_{\{T(u)\cap\mathcal B(\tau_i,\delta)=\emptyset\}}.
	\end{equation}
	Under Assumption~\ref{assumptions},
	Lemma~\ref{S-lem:support_and_exclusion}\textup{(ii)-(iii)}
	shows that the detection cubes are pairwise disjoint and that
	each contains at most one candidate, and hence at most one
	reported location. Therefore,
	\begin{equation}\label{S-eq:bh_W_equals_N_minus_D}
		W(u)
		=
		\sum_{i=1}^N
		\mathbf1_{\{T(u)\cap\mathcal B(\tau_i,\delta)\ne\emptyset\}}
		=
		N-D(u).
	\end{equation}
	The definition of power in~\eqref{def:power} gives
	\begin{equation}\label{S-eq:bh_power_missed_identity}
		\Power(u)=1-\E[D(u)/N].
	\end{equation}
	For \(u\in\R\), let \(U_S(u)\) count the grid cubes whose
	maximum noise score is at least \(u\):
	\begin{equation}\label{S-eq:noise_exceedance_count}
		U_S(u):=\sum_{m\in\mathcal I_L}
		\mathbf1_{\{\max_{t\in\overline Q_m}S^z(t)\ge u\}},
		\qquad U_S(+\infty):=0.
	\end{equation}
	Under Assumption~\ref{ass:noise_reg_main}, stationarity and
	Theorem~\ref{thm:continuity_test} give
	\(\Prob(\max_{\overline Q_m}S^z\ge u)=p(u)\) for all \(m\in\mathcal I_L\), so
	\begin{equation}\label{S-eq:noise_box_expectation_identity}
		\E[U_S(u)]=M_Lp(u),
		\qquad u\in\R.
	\end{equation}
	
	\begin{lemma}[False detections and missed object instances]
		\label{S-lem:bh_false_candidates_missed}
		Suppose that Assumption~\ref{assumptions} holds.
		For every realization of the noise field \(z\) and every
		\(u\in\R\cup\{+\infty\}\),
		\begin{equation}\label{S-eq:bh_false_candidates_missed}
			V(u)\le U_S(u)+2^dD(u).
		\end{equation}
		If in addition every detection cube contains a candidate from Step~\ref{item:candidate_peaks},
		then, for every \(u\in\R\cup\{+\infty\}\),
		\begin{equation}\label{S-eq:false_detections_full_candidates}
			V(u)\le U_S(u).
		\end{equation}
	\end{lemma}
	
	\begin{proof}
		Fix a realization of \(z\). For \(u=+\infty\), the report set \(T(u)\) is
		empty and \(V(u)=U_S(u)=0\). Let \(u\in\R\). Every false report
		\(t\in T(u)\cap\mathbb B_0^\delta\) satisfies \(t\in T\) and
		\(S^y(t)\ge u\). We split the false reports according to whether the
		clean score map contributes at \(t\): reports with
		\(\mathbf X^x(t)=0\) are bounded by noise exceedances, and reports
		with \(\mathbf X^x(t)\ne0\) by missed instances.
		We start with false reports with \(\mathbf X^x(t)=0\). By linearity, \(\mathbf X^y=\mathbf X^x+\mathbf X^z\), so
		\(S^z(t)=S^y(t)\ge u\). The half-open grid
		cubes form a partition of \(\R^d\), so \(t\) lies in exactly one of
		them, \(Q_m\), and \(m\in\mathcal I_L\) because \(t\) lies in the
		search domain. Since
		\(\max_{s\in\overline Q_m}S^z(s)\ge S^z(t)\ge u\), the cube \(Q_m\) is
		counted by \(U_S(u)\). Two points in the same half-open grid cube
		have \(\ell_\infty\)-distance strictly less than \(r/2\), whereas
		distinct candidates have \(\ell_\infty\)-distance at least \(r/2\);
		see Step~\ref{item:candidate_peaks}. Hence, distinct false reports with
		\(\mathbf X^x(t)=0\) lie in distinct grid cubes counted by
		\(U_S(u)\), so there are at most \(U_S(u)\) of them.
		Next, we consider false reports with \(\mathbf X^x(t)\ne0\). We assign
		to each such report a missed instance, and then show that at most
		\(2^d\) reports are assigned to the same instance.
		Since \(\mathbf X^x(t)\ne0\), we have \((x*\tilde\psi_j)(t)\ne0\) for
		some \(j\). Since \(x=\sum_{i=1}^Nx_i\), linearity gives an index \(i\)
		with \((x_i*\tilde\psi_j)(t)\ne0\). If there are several, we choose one
		and denote it by \(i(t)\). Since \((x_{i(t)}*\tilde\psi_j)(t)\ne0\),
		\eqref{S-eq:candidate_support_cutoff} forces
		\(\|t-\tau_{i(t)}\|_\infty<B\), and
		\eqref{S-eq:candidate_local_exclusion} gives
		\(\mathcal B(\tau_{i(t)},\delta)\subset\mathcal B(t,r)\). A candidate \(s\in\mathcal B(\tau_{i(t)},\delta)\) would belong to
		\(\mathcal B(t,r)\) and hence satisfy \(\|s-t\|_\infty<r/2\).
		Since \(t\in\mathbb B_0^\delta\) lies outside every detection cube,
		\(s\) and \(t\) would be distinct candidates, contradicting their
		separation by at least \(r/2\). Since \(T(u)\subseteq T\), the cube
		\(\mathcal B(\tau_{i(t)},\delta)\) also contains no reported location,
		so instance \(i(t)\) is missed and is counted by \(D(u)\).
		Now fix a missed instance \(i\). Every false report \(t\) with
		\(i(t)=i\) satisfies \(\|t-\tau_i\|_\infty<B\), and therefore lies in
		the cube \(\mathcal B(\tau_i,2B)\). Splitting each coordinate interval
		of this cube at the corresponding coordinate of \(\tau_i\) covers it
		by \(2^d\) closed cubes of side length \(B\). Two points in the same
		closed cube have \(\ell_\infty\)-distance at most \(B\), whereas
		distinct candidates have \(\ell_\infty\)-distance at least
		\(r/2=B+\delta/2>B\). Each of the \(2^d\) closed cubes therefore
		contains at most one false report \(t\) with \(i(t)=i\). Summing over
		the \(D(u)\) missed instances, there are at most \(2^dD(u)\) false
		reports with \(\mathbf X^x(t)\ne0\).
		Finally, suppose that every detection cube contains a candidate. We have shown that a false report with
		\(\mathbf X^x(t)\ne0\) would leave the detection cube
		\(\mathcal B(\tau_{i(t)},\delta)\) without a candidate. Hence every
		false report satisfies \(\mathbf X^x(t)=0\), and the first bound
		gives~\eqref{S-eq:false_detections_full_candidates}.
	\end{proof}
	
	\subsection{Noise bounds}
	
	\begin{lemma}[Bounds on noise maxima]
		\label{S-lem:noise_maxima_fixed_growing}
		Suppose that Assumption~\ref{ass:noise_reg_main} holds.
		Recall \(Z_{\max}\) and \(p\) from
		Subsection~\ref{subsec:atomlessness} and
		\(\eta_{\mathrm{peak}}\) from~\eqref{eq:eta_peak_def}.
		Set \(R_L:=\max_{t\in\overline{\mathcal B}(0,L)}S^z(t)\),
		and define
		\begin{equation}\label{S-eq:projected_noise_max_variance}
			\sigma_{\mathcal V}^2
			:=\lambda_{\max}(\Sigma_{\mathbf X})
			=\max_{\substack{v\in\R^M\\ \|v\|=1}}
			v^\top\Sigma_{\mathbf X}v.
		\end{equation}
		Let \(\sigma_{\mathcal V}\) denote the nonnegative square root
		of \(\sigma_{\mathcal V}^2\).
		Then \(0<\sigma_{\mathcal V}^2\le\eta_{\mathrm{peak}}<\infty\) and 
		for every \(u\ge\eta_{\mathrm{peak}}\),
		\begin{equation}\label{S-eq:local_noise_sigma_bound}
			p(u)\le\exp\!\left\{
			-\frac{\bigl(\sqrt u-\sqrt{\eta_{\mathrm{peak}}}\bigr)^2}
			{2\sigma_{\mathcal V}^2}
			\right\},
		\end{equation}
		and hence
		\begin{equation}\label{S-eq:local_noise_eta_bound}
			p(u)\le\exp\!\left\{
			-\frac12\left(
			\sqrt{\frac{u}{\eta_{\mathrm{peak}}}}-1
			\right)^2
			\right\}.
		\end{equation}
		For every \(L>0\) and \(u\in\R\),
		\begin{equation}\label{S-eq:growing_noise_union_bound}
			\Prob\{R_L>u\}\le M_Lp(u).
		\end{equation}
		Moreover, for every \(L>0\),
		\begin{equation}\label{S-eq:growing_noise_log_bound}
			\Prob\!\left\{
			R_L>\eta_{\mathrm{peak}}
			\bigl(1+2\sqrt{\log M_L}\bigr)^2
			\right\}
			\le\frac1{M_L}.
		\end{equation}
	\end{lemma}
	
	\begin{proof}
		Write \(K=\overline{\mathcal B}(0,r/2)\). Since
		\(\max_{\|v\|=1}\langle v,x\rangle=\|x\|\) for every \(x\in\R^M\),
		\[
		\sqrt{Z_{\max}}
		=\sup_{\substack{t\in K\\v\in\R^M,\ \|v\|=1}}
		\langle v,\mathbf X^z(t)\rangle.
		\]
		By Lemma~\ref{S-lem:projected_noise_continuity}, the process
		\((t,v)\mapsto\langle v,\mathbf X^z(t)\rangle\) is a centered Gaussian
		process with continuous sample paths on the compact set
		\(K\times\{v\in\R^M:\|v\|=1\}\). It is therefore separable, and its
		supremum \(\sqrt{Z_{\max}}\) is finite almost surely. By
		stationarity, \(\Var\langle v,\mathbf X^z(t)\rangle
		=v^\top\Sigma_{\mathbf X}v\), so the maximal variance of the process
		is \(\sigma_{\mathcal V}^2\), which is positive by
		Assumption~\ref{ass:noise_reg_main}. By the Borell-TIS inequality \cite[Theorem~2.1.1]{adlerTaylor2007},
		we have \(\E\sqrt{Z_{\max}}<\infty\) and
		\begin{equation}\label{S-eq:noise_maximum_borell}
			\Prob\bigl\{\sqrt{Z_{\max}}>\E\sqrt{Z_{\max}}+h\bigr\}
			\le\exp\!\left\{-\frac{h^2}{2\sigma_{\mathcal V}^2}\right\},
			\qquad h\ge0.
		\end{equation}
		Integrating this inequality gives
		\(\eta_{\mathrm{peak}}=\E Z_{\max}<\infty\), and Jensen's inequality
		gives \(\E\sqrt{Z_{\max}}\le\sqrt{\eta_{\mathrm{peak}}}\). Since
		\(\mathbf X^z(0)\) is centered,
		\[
		\sigma_{\mathcal V}^2
		=\lambda_{\max}(\Sigma_{\mathbf X})
		\le\operatorname{tr}\Sigma_{\mathbf X}
		=\sum_{j=1}^M\Var X_j^z(0)
		=\E\|\mathbf X^z(0)\|^2.
		\]
		Since \(0\in K\), \(\|\mathbf X^z(0)\|^2\le Z_{\max}\), and taking
		expectations gives \(\E\|\mathbf X^z(0)\|^2\le\eta_{\mathrm{peak}}\).
		We then get
		\begin{equation}\label{S-eq:noise_maximum_parameters}
			0<\sigma_{\mathcal V}^2\le\E\|\mathbf X^z(0)\|^2
			\le\eta_{\mathrm{peak}}.
		\end{equation}
		Let \(u\ge\eta_{\mathrm{peak}}\). Since \(Z_{\max}\ge0\),
		\(p(u)=\Prob\{\sqrt{Z_{\max}}>\sqrt u\}\), and
		\(\sqrt u-\E\sqrt{Z_{\max}}\ge\sqrt u-\sqrt{\eta_{\mathrm{peak}}}\ge0\).
		Applying~\eqref{S-eq:noise_maximum_borell} with
		\(h=\sqrt u-\E\sqrt{Z_{\max}}\) gives
		\[
		p(u)
		\le\exp\!\left\{-\frac{\bigl(\sqrt u-\E\sqrt{Z_{\max}}\bigr)^2}
		{2\sigma_{\mathcal V}^2}\right\}
		\le\exp\!\left\{
		-\frac{\bigl(\sqrt u-\sqrt{\eta_{\mathrm{peak}}}\bigr)^2}
		{2\sigma_{\mathcal V}^2}\right\},
		\]
		which is~\eqref{S-eq:local_noise_sigma_bound}. Since
		\(\sigma_{\mathcal V}^2\le\eta_{\mathrm{peak}}\),
		\[
		\frac{\bigl(\sqrt u-\sqrt{\eta_{\mathrm{peak}}}\bigr)^2}
		{2\sigma_{\mathcal V}^2}
		\ge\frac{\bigl(\sqrt u-\sqrt{\eta_{\mathrm{peak}}}\bigr)^2}
		{2\eta_{\mathrm{peak}}}
		=\frac12\left(\sqrt{\frac u{\eta_{\mathrm{peak}}}}-1\right)^2,
		\]
		which gives~\eqref{S-eq:local_noise_eta_bound}.
		The closed grid cubes \(\overline Q_m\), \(m\in\mathcal I_L\), cover
		the search domain, so
		\(R_L\le\max_{m\in\mathcal I_L}\max_{t\in\overline Q_m}S^z(t)\). Each
		\(\overline Q_m\) is a translate of \(K\), so by stationarity
		\(\Prob\{\max_{t\in\overline Q_m}S^z(t)>u\}=p(u)\). A union bound
		over the \(M_L\) cubes gives~\eqref{S-eq:growing_noise_union_bound}.
		For \(u_L:=\eta_{\mathrm{peak}}(1+2\sqrt{\log M_L})^2\ge
		\eta_{\mathrm{peak}}\), we have
		\(\sqrt{u_L/\eta_{\mathrm{peak}}}-1=2\sqrt{\log M_L}\),
		so~\eqref{S-eq:local_noise_eta_bound} gives
		\(p(u_L)\le e^{-2\log M_L}=M_L^{-2}\).
		Combining with~\eqref{S-eq:growing_noise_union_bound}, we get
		\(\Prob\{R_L>u_L\}\le M_L\cdot M_L^{-2}=M_L^{-1}\), which
		proves~\eqref{S-eq:growing_noise_log_bound} as desired.
	\end{proof}
	
	\begin{corollary}[Bounds on detection thresholds]
		\label{S-cor:calibrated_thresholds}
		Suppose that Assumption~\ref{ass:noise_reg_main} holds.
		Let \(\eta_{\mathrm{peak}}\) and \(\sigma_{\mathcal V}>0\)
		be as in Lemma~\ref{S-lem:noise_maxima_fixed_growing},
		and let \(u_c\) be defined by~\eqref{eq:uc_def_inline}.
		For \(c\in(0,1)\), let \(q_c>0\) satisfy
		\(\Prob(G^2>q_c)=c\), where \(G\sim N(0,1)\).
		Then
		\begin{equation}\label{S-eq:calibrated_threshold_quantile_bounds}
			\sigma_{\mathcal V}^2q_c
			\le u_c
			\le
			\left(
			\sqrt{\eta_{\mathrm{peak}}}
			+\sigma_{\mathcal V}\sqrt{2\log(1/c)}
			\right)^2.
		\end{equation}
		If \(c\le\Prob(G^2>1)\), then
		\begin{equation}\label{S-eq:calibrated_threshold_explicit_bounds}
			u_c
			\ge
			2\sigma_{\mathcal V}^2\log(1/c)
			-\sigma_{\mathcal V}^2
			\log\!\bigl(4\pi\log(1/c)\bigr).
		\end{equation}
		Consequently, for fixed \(0<\alpha<1\), every
		\(0<\varepsilon<1\), and all sufficiently large \(L\),
		\begin{equation}\label{S-eq:bonferroni_threshold_order}
			2(1-\varepsilon)\sigma_{\mathcal V}^2
			\log\!\left(\frac{M_L}{\alpha}\right)
			\le u_{\rm Bon}
			\le
			2(1+\varepsilon)\sigma_{\mathcal V}^2
			\log\!\left(\frac{M_L}{\alpha}\right).
		\end{equation}
	\end{corollary}
	
	\begin{proof}
		Fix \(c\in(0,1)\) and set
		\(U_c:=\bigl(\sqrt{\eta_{\mathrm{peak}}}
		+\sigma_{\mathcal V}\sqrt{2\log(1/c)}\bigr)^2\).
		Since \(\sigma_{\mathcal V}>0\), \(\log(1/c)>0\), we have
		\(U_c>\eta_{\mathrm{peak}}\), so~\eqref{S-eq:local_noise_sigma_bound}
		gives \(p(U_c)\le e^{-\log(1/c)}=c\).
		Since \(u_c=\inf\{u\in\R:p(u)\le c\}\), this gives \(u_c\le U_c\),
		which is the upper bound
		in~\eqref{S-eq:calibrated_threshold_quantile_bounds}.
		For the lower bound, let \(v_0\) be a unit eigenvector of
		\(\Sigma_{\mathbf X}\) with eigenvalue
		\(\sigma_{\mathcal V}^2\).
		It holds that
		\[
		Z_{\max}
		\ge \|\mathbf X^z(0)\|^2
		\ge \langle v_0,\mathbf X^z(0)\rangle^2,
		\]
		and \(\langle v_0,\mathbf X^z(0)\rangle/\sigma_{\mathcal V}\)
		has the same law as \(G\sim N(0,1)\).
		By Corollary~\ref{S-cor:generalized_inverse_attainment},
		\(u_c>0\) and \(p(u_c)=c\). If \(u_c<\sigma_{\mathcal V}^2q_c\), then,
		since \(x\mapsto\Prob(G^2>x)\) is strictly decreasing on
		\([0,\infty)\),
		\[
		p(u_c)\ge\Prob\!\left(G^2>\frac{u_c}{\sigma_{\mathcal V}^2}\right)
		>\Prob(G^2>q_c)=c,
		\]
		a contradiction. Hence \(u_c\ge\sigma_{\mathcal V}^2q_c\), which is
		the lower bound in~\eqref{S-eq:calibrated_threshold_quantile_bounds}.
		To prove~\eqref{S-eq:calibrated_threshold_explicit_bounds},
		use the Gaussian tail inequalities
		\[
		\frac{x}{1+x^2}\cdot\frac{e^{-x^2/2}}{\sqrt{2\pi}}
		\le \Prob(G>x)
		\le \frac12e^{-x^2/2},
		\qquad x>0.
		\]
		At \(x=\sqrt{q_c}\), symmetry and the upper bound give
		\(c=2\Prob(G>\sqrt{q_c})\le e^{-q_c/2}\), so
		\(q_c\le 2\log(1/c)\).
		If \(c\le\Prob(G^2>1)\), then \(q_c\ge1\) since  \(x\mapsto\Prob(G^2>x)\) is strictly
		decreasing on \([0,\infty)\). Applying the lower
		Gaussian tail bound gives
		\[
		c
		\ge
		\frac{2\sqrt{q_c}}{1+q_c}
		\frac{e^{-q_c/2}}{\sqrt{2\pi}}
		\ge
		\frac{e^{-q_c/2}}{\sqrt{2\pi q_c}}.
		\]
		Taking logarithms and using \(q_c\le2\log(1/c)\), we obtain
		\[
		\begin{aligned}
			q_c
			&\ge 2\log(1/c)-\log(2\pi q_c)\\
			&\ge 2\log(1/c)-\log\!\bigl(4\pi\log(1/c)\bigr).
		\end{aligned}
		\]
		Combining this with \(u_c\ge\sigma_{\mathcal V}^2q_c\)
		proves~\eqref{S-eq:calibrated_threshold_explicit_bounds}.
		Finally, for \(c\le\Prob(G^2>1)\), the upper bound
		in~\eqref{S-eq:calibrated_threshold_quantile_bounds}
		and the lower bound
		in~\eqref{S-eq:calibrated_threshold_explicit_bounds} give
		\[
		\begin{aligned}
			1-\frac{\log\!\bigl(4\pi\log(1/c)\bigr)}
			{2\log(1/c)}
			\le
			\frac{u_c}{2\sigma_{\mathcal V}^2\log(1/c)}
			\le
			\left(
			1+\sqrt{\frac{\eta_{\mathrm{peak}}}
				{2\sigma_{\mathcal V}^2\log(1/c)}}
			\right)^2.
		\end{aligned}
		\]
		Since \(\eta_{\mathrm{peak}}\) and \(\sigma_{\mathcal V}>0\)
		are fixed, both bounds tend to one as \(c\downarrow0\). Hence
		\begin{equation}\label{S-eq:noise_quantile_asymptotic}
			\lim_{c\downarrow0}
			\frac{u_c}{2\sigma_{\mathcal V}^2\log(1/c)}=1.
		\end{equation}
		By~\eqref{eq:u_bon_def} and~\eqref{eq:uc_def_inline},
		\(u_{\rm Bon}=u_{\alpha/M_L}\).
		Since \(M_L\to\infty\) by~\eqref{S-eq:grid_count_asymptotic},
		\(c=\alpha/M_L\to0\) and \(\log(1/c)=\log(M_L/\alpha)\).
		Hence, for every \(0<\varepsilon<1\) and all sufficiently
		large \(L\),
		\[
		\frac{u_{\rm Bon}}{2\sigma_{\mathcal V}^2\log(M_L/\alpha)}
		\in[1-\varepsilon,1+\varepsilon],
		\]
		which proves~\eqref{S-eq:bonferroni_threshold_order}.
	\end{proof}
	
	\subsection{Candidate selection}
	\label{S-app:candidate_selection}
	
	\begin{lemma}[Bounds on the clean score at representative points]
		\label{S-lem:representative_scores}
Suppose that Assumption~\ref{assumptions} holds,
that \(a_L>0\), and that
Assumption~\ref{def:localization_property} holds
with localization constant \(\rho\).
		In the template model \eqref{eq:template_instance_model}, assume
		that \(\mathcal V\) is steerable and define
		\begin{equation}\label{S-eq:kappa_V_def}
			\kappa_{\mathcal V}
			:=\max_{\substack{1\le\ell\le n\\P_{\mathcal V}u_\ell\ne0}}
			\frac{\|u_\ell\|_{L^2(D_{B/2})}}
			{\|P_{\mathcal V}u_\ell\|_{L^2(D_{B/2})}}.
		\end{equation}
		In the linear subspace model of
		Corollary~\ref{S-cor:exact_linear_subspace}, set
		\(\kappa_{\mathcal V}:=1\).
		
		In either model, \(1\le\kappa_{\mathcal V}<\infty\), and the
		representative points \(q_1,\ldots,q_N\) can be chosen so that
		the localization inequality holds with the same constant \(\rho\)
		and
		\begin{equation}\label{S-eq:bh_representative_clean_score_bounds}
			\|a_i\|^2\le S^x(q_i)\le\kappa_{\mathcal V}^2\|a_i\|^2,
			\qquad i=1,\ldots,N.
		\end{equation}
	\end{lemma}
	
	\begin{proof}
		Fix \(i\) and \(t\in\R^d\), and write \(x_{i,t}(s):=x_i(t+s)\)
		for \(s\in D_{B/2}\). Since
		\((x_i*\tilde\psi_m)(t)=\langle x_{i,t},\psi_m\rangle_{L^2(D_{B/2})}\),
		the convolution values are the coefficients of \(P_{\mathcal V}x_{i,t}\)
		in the orthonormal basis of \(\mathcal V\). Hence
		\begin{equation}\label{S-eq:single_instance_bessel}
			\left(\sum_{m=1}^M(x_i*\tilde\psi_m)^2(t)\right)^{1/2}
			=\|P_{\mathcal V}x_{i,t}\|_{L^2(D_{B/2})}
			\le\|x_{i,t}\|_{L^2(D_{B/2})}
			\le\|x_i\|_{L^2(\R^d)}.
		\end{equation}
		
		In the template model,
		\(\|x_i\|_{L^2(\R^d)}=b_i\|u_{\ell_i}\|_{L^2(D_{B/2})}\).
		By \eqref{eq:x} and Lemma~\ref{S-lem:proj_rot_appendix},
		\[
		\|a_i\|=b_i\|P_{\mathcal V}u_{\ell_i}\|_{L^2(D_{B/2})}
		\ge a_L>0.
		\]
		Thus \(u_{\ell_i}\) is included in the maximum in
		\eqref{S-eq:kappa_V_def}, and
		\(\|x_i\|_{L^2(\R^d)}\le\kappa_{\mathcal V}\|a_i\|\).
		The maximum is finite because the template family is finite and
		every denominator in the maximum is positive. Each ratio is at
		least one because orthogonal projection does not increase the norm.
		In the linear subspace model, the centered instance
		\(x_i(\tau_i+\cdot)\) belongs to \(\mathcal V\) and has coefficient
		vector \(a_i\) in its orthonormal basis. Hence
		\(\|x_i\|_{L^2(\R^d)}=\|a_i\|\), so the same bound holds with
		\(\kappa_{\mathcal V}=1\).
		Fix \(i\). By Lemma~\ref{S-lem:support_and_exclusion}\textup{(i)},
		only instance \(i\) contributes to the clean score in
		\(\mathcal B(\tau_i,\delta)\), and \(S^x(\tau_i)=\|a_i\|^2\).
		The same conclusions hold in the linear subspace model by separation
		of the supports and the definition of \(a_i\). Thus, for every
		\(q\in\mathcal B(\tau_i,\delta)\), \eqref{S-eq:single_instance_bessel}
		gives
		\[
		\sqrt{S^x(q)}
		=\|P_{\mathcal V}x_{i,q}\|_{L^2(D_{B/2})}
		\le \|x_i\|_{L^2(\R^d)}
		\le \kappa_{\mathcal V}\|a_i\|.
		\]
		Choose \(q_i\) as in Assumption~\ref{def:localization_property}.
		If \(S^x(q_i)<S^x(\tau_i)\), replace \(q_i\) by \(\tau_i\).
		The center satisfies \eqref{eq:qi}, and the comparison set
		\eqref{eq:comparison set} depends only on the instance locations
		and their ordering. The replacement therefore increases the score
		difference at every point in the comparison set, so the localization
		inequality still holds with the same constant \(\rho\).
		With this choice, \(S^x(q_i)\ge S^x(\tau_i)=\|a_i\|^2\).
		Since \(q_i\in\mathcal B(\tau_i,\delta)\), the upper bound proved
		for every point in the cube also applies to \(q_i\), giving
		\eqref{S-eq:bh_representative_clean_score_bounds}.
	\end{proof}
	
	\begin{lemma}[Candidate selection in the absence of noise]
		\label{S-lem:clean_candidate_selection}
		Suppose that Assumptions~\ref{def:localization_property}
and~\ref{assumptions} hold, and that \(a_L>0\).
		Then the first \(N\) candidates selected from \(S^x\) by
		Step~\ref{item:candidate_peaks} consist of exactly one point
		in each detection cube \(\mathcal B(\tau_i,\delta)\).
	\end{lemma}
	
	\begin{proof}
		Order the instances and choose representatives \(q_1,\ldots,q_N\)
		as in Assumption~\ref{def:localization_property}. We prove by
		induction that the first \(k\) candidates exist and lie in \(k\)
		distinct detection cubes, for \(k=0,\ldots,N\). The claim holds
		for \(k=0\). Suppose it holds for some \(k<N\), and let \(A\) be
		the set of indices of the cubes containing the first \(k\) candidates.
		By \eqref{S-eq:candidate_geometric_facts}, a candidate selected in
		\(\mathcal B(\tau_j,\delta)\) removes \(\mathcal B(\tau_j,2B)\)
		and leaves every other detection cube untouched. Since
		\(\delta\le 2B\), no point of a detection cube indexed by \(A\)
		remains. Let \(i:=\min(\{1,\ldots,N\}\setminus A)\). Then \(q_i\)
		remains in the search domain, and \(\mathcal B(\tau_j,2B)\) has
		been removed for every \(j<i\).
		Any remaining point \(t\) outside all detection cubes therefore
		belongs to the comparison set \eqref{eq:comparison set} of
		instance \(i\). By localization,
		\[
		S^x(q_i)-S^x(t)\ge \rho\|a_i\|^2\ge \rho a_L^2>0.
		\]
		The remaining domain is nonempty because it contains \(q_i\),
		and it is compact because the exclusion cubes are open.
		Continuity of translations in \(L^2(\R^d)\) gives continuity of
		each \(x*\tilde\psi_m\), and hence of \(S^x\). Thus \(S^x\)
		attains its maximum on the remaining domain. Every maximizer has
		score at least \(S^x(q_i)\), so it lies in a detection cube whose
		index is outside \(A\). The \((k+1)\)st candidate therefore lies
		in a new detection cube, completing the induction.
	\end{proof}
	
	By Lemma~\ref{S-lem:clean_candidate_selection}, the first \(N\)
	locations selected from the clean score \(S^x\) include one point
	in each detection cube. For the observed score \(S^y\), we give
	conditions under which the cube \(\mathcal B(\tau_i,\delta)\)
	contains a selected location \(t\) with
	\(S^y(t)\ge \|a_i\|^2/4\). 
	Throughout the next four lemmas, assume the hypotheses of
	Lemma~\ref{S-lem:representative_scores} and
	Assumption~\ref{ass:noise_reg_main}. Order the instances as in
	Assumption~\ref{def:localization_property} and take the
	representatives \(q_i\) and the constant \(\kappa_{\mathcal V}\)
	from Lemma~\ref{S-lem:representative_scores}. We use the continuous
	version of the projected noise from
	Lemma~\ref{S-lem:projected_noise_continuity}. As shown in
	Subsection~\ref{subsec:algorithm_outline},
	Step~\ref{item:candidate_peaks} exhausts the search domain
	after at most \(M_L\) selections, so every \(q_i\) is eventually
	removed.
	
	\begin{definition}[Blocked instances and local noise events]
		\label{S-def:blocked_instances_local_noise}
		For each instance \(i\), let \(\widehat t_i\) be the first
		selected location whose exclusion cube contains \(q_i\).
		We call \(i\) \emph{blocked} if
		\(\widehat t_i\notin\mathcal B(\tau_i,\delta)\).
		Any point in \(\mathcal B(\tau_i,\delta)\) is at
		\(\ell_\infty\)-distance less than \(\delta\le r/2\) from
		\(q_i\), so for a blocked instance \(q_i\) is removed before
		any location in its detection cube is selected.
		We choose \(\gamma_{\rm BH}>0\) so that
		\[
		\sqrt{\gamma_{\rm BH}}\le\frac12
		\qquad\text{and}\qquad
		4\kappa_{\mathcal V}\sqrt{\gamma_{\rm BH}}\le\frac{\rho}{2}.
		\]
		Equivalently, \(\gamma_{\rm BH}\) must be at most \(1/4\) and
		at most \(\rho^2/(64\kappa_{\mathcal V}^2)\). We therefore set
		\begin{equation}\label{S-eq:bh_gamma_choice}
			\gamma_{\rm BH}
			:=\min\left\{\frac14,\frac{\rho^2}{64\kappa_{\mathcal V}^2}\right\}.
		\end{equation}
		For each instance \(i\), define the local noise event
		\begin{equation}\label{S-eq:bh_local_noise_event}
			E_i^{\rm loc}
			:=\left\{
			\sup_{\|t-\tau_i\|_\infty\le B+\delta} S^z(t)
			>\gamma_{\rm BH}\|a_i\|^2
			\right\},
		\end{equation}
		where the supremum is restricted to points \(t\) in the search
		domain \(\overline{\mathcal B}(0,L)\).
	\end{definition}
	
	\begin{lemma}[Candidate scores near unblocked instances]
		\label{S-lem:bh_high_score_true_candidates}
		Let \(G_L\) be the number of unblocked instances \(i\) for
		which \(E_i^{\rm loc}\) does not occur. For each such \(i\),
		\(\widehat t_i\in\mathcal B(\tau_i,\delta)\) and
		\begin{equation}\label{S-eq:bh_candidate_score_level}
			S^y(\widehat t_i)\ge\frac{\|a_i\|^2}{4}
			\ge\frac{a_L^2}{4}.
		\end{equation}
		These \(G_L\) locations are distinct.
	\end{lemma}
	
	\begin{proof}
		Fix any instance \(i\). When \(\widehat t_i\) is selected,
		\(q_i\) still belongs to the domain over which
		\(\widehat t_i\) maximizes \(S^y\). Hence
		\begin{equation}\label{S-eq:first_removal_score_comparison}
			S^y(\widehat t_i)\ge S^y(q_i).
		\end{equation}
		Moreover, \(q_i\in\mathcal B(\tau_i,\delta)\) and
		\(q_i\in\mathcal B(\widehat t_i,r)\), so
		\[
		\|\widehat t_i-\tau_i\|_\infty
		\le \|\widehat t_i-q_i\|_\infty+\|q_i-\tau_i\|_\infty
		< r/2+\delta/2=B+\delta.
		\]
		Thus both \(q_i\) and \(\widehat t_i\) lie in the region
		in \eqref{S-eq:bh_local_noise_event}. On \((E_i^{\rm loc})^c\),
		\eqref{S-eq:score_triangle_bounds} gives
		\begin{equation}\label{S-eq:local_score_perturbation}
			\bigl|\sqrt{S^y(t)}-\sqrt{S^x(t)}\bigr|
			\le\sqrt{\gamma_{\rm BH}}\|a_i\|,
			\qquad t\in\{q_i,\widehat t_i\}.
		\end{equation}
		Suppose \(i\) is counted by \(G_L\). Then
		\(\widehat t_i\in\mathcal B(\tau_i,\delta)\) and
		\(E_i^{\rm loc}\) does not occur. By
		\eqref{S-eq:local_score_perturbation} at \(t=q_i\) and
		\eqref{S-eq:bh_representative_clean_score_bounds},
		\[
		\sqrt{S^y(q_i)}
		\ge \sqrt{S^x(q_i)}-\sqrt{\gamma_{\rm BH}}\|a_i\|
		\ge (1-\sqrt{\gamma_{\rm BH}})\|a_i\|
		\ge \frac{\|a_i\|}{2},
		\]
		where the last inequality follows from \eqref{S-eq:bh_gamma_choice}.
		Hence \eqref{S-eq:first_removal_score_comparison} gives
		\[
		S^y(\widehat t_i)
		\ge S^y(q_i)
		\ge \frac{\|a_i\|^2}{4}
		\ge \frac{a_L^2}{4}.
		\]
		For any two such indices \(i\ne j\),
		Lemma~\ref{S-lem:support_and_exclusion}\textup{(ii)} gives
		\[
		\widehat t_i\in\mathcal B(\tau_i,\delta),
		\qquad
		\widehat t_j\in\mathcal B(\tau_j,\delta),
		\qquad
		\mathcal B(\tau_i,\delta)\cap\mathcal B(\tau_j,\delta)
		=\emptyset,
		\]
		so \(\widehat t_i\ne\widehat t_j\).
	\end{proof}

	Even when the noise score near \(\tau_i\) is small, a selected
	location outside \(\mathcal B(\tau_i,\delta)\) may remove \(q_i\)
	before any location in that cube is selected. 
	
	\begin{lemma}[Blocking when the local noise score is small]
		\label{S-lem:blocking_candidate}
		Suppose that instance \(i\) is blocked and that \(E_i^{\rm loc}\)
		does not occur. Then \(\widehat t_i\) lies outside every detection
		cube \(\mathcal B(\tau_j,\delta)\), and
		\(\widehat t_i\in\mathcal B(\tau_j,2B)\) for at least one \(j<i\).
		Let \(h\) be the smallest such index. Then,
		\begin{enumerate}
			\item[(i)] Instance \(h\) is blocked.
			\item[(ii)] \(S^x(q_i)\le\theta S^x(q_h)\), where
			\begin{equation}\label{S-eq:bh_multiplicative_clean_score_drop}
				\theta
				:=1-\frac{\min\{\rho,\kappa_{\mathcal V}^2\}}
				{2\kappa_{\mathcal V}^2}
				\in[1/2,1).
			\end{equation}
			\item[(iii)] \(\|\tau_i-\tau_h\|_\infty<r=2B+\delta\).
		\end{enumerate}
	\end{lemma}
	
	\begin{proof}
		Since \(i\) is blocked, \(\widehat t_i\notin\mathcal B(\tau_i,\delta)\)
		by Definition~\ref{S-def:blocked_instances_local_noise}.
		By definition of \(\widehat t_i\), its exclusion cube contains the
		representative \(q_i\). Since \(q_i\in\mathcal B(\tau_i,\delta)\),
		\[
		q_i\in\mathcal B(\widehat t_i,r)\cap\mathcal B(\tau_i,\delta).
		\]
		By \eqref{S-eq:candidate_geometric_facts}, this intersection is empty
		whenever \(\widehat t_i\in\mathcal B(\tau_j,\delta)\) for some
		\(j\ne i\). Thus \(\widehat t_i\) lies outside every detection cube.
		Set
		\[
		h:=\min\bigl(\{i\}\cup
		\{j<i:\widehat t_i\in\mathcal B(\tau_j,2B)\}\bigr).
		\]
		By the choice of \(h\) and since \(\widehat t_i\) lies outside
		every detection cube,
		\[
		\widehat t_i\in\overline{\mathcal B}(0,L)\setminus
		\left(
		\bigcup_{j<h}\mathcal B(\tau_j,2B)
		\;\cup\;
		\bigcup_{j=h}^{N}\mathcal B(\tau_j,\delta)
		\right).
		\]
		Since \(\widehat t_i\) belongs to the comparison set
		\eqref{eq:comparison set} for instance \(h\), the localization
		property of Assumption~\ref{def:localization_property} and the
		representative score bound \eqref{S-eq:bh_representative_clean_score_bounds}
		give
		\begin{equation}\label{S-eq:bh_additive_clean_score_drop}
			S^x(q_h)-S^x(\widehat t_i)\ge\rho\|a_h\|^2,
			\qquad
			\sqrt{S^x(\widehat t_i)}<\sqrt{S^x(q_h)}
			\le\kappa_{\mathcal V}\|a_h\|.
		\end{equation}
		
		On \((E_i^{\rm loc})^c\), \eqref{S-eq:local_score_perturbation} and
		\eqref{S-eq:first_removal_score_comparison} give
		\[
		\begin{aligned}
			\sqrt{S^x(q_i)}-\sqrt{S^x(\widehat t_i)}
			&\le \sqrt{S^y(q_i)}-\sqrt{S^y(\widehat t_i)}
			+2\sqrt{\gamma_{\rm BH}}\|a_i\|\\
			&\le 2\sqrt{\gamma_{\rm BH}}\|a_i\|.
		\end{aligned}
		\]
		Since \(h\le i\), the ordering in
		Assumption~\ref{def:localization_property} gives
		\(\|a_i\|\le\|a_h\|\). Multiplying by
		\(\sqrt{S^x(q_i)}+\sqrt{S^x(\widehat t_i)}\) and using
		\eqref{S-eq:bh_representative_clean_score_bounds},
		\eqref{S-eq:bh_additive_clean_score_drop}, and \eqref{S-eq:bh_gamma_choice},
		we obtain
		\[
		\begin{aligned}
			S^x(q_i)-S^x(\widehat t_i)
			&\le 2\sqrt{\gamma_{\rm BH}}\|a_i\|\,
			\kappa_{\mathcal V}(\|a_i\|+\|a_h\|)\\
			&\le 4\kappa_{\mathcal V}\sqrt{\gamma_{\rm BH}}\|a_h\|^2
			\le \frac{\rho}{2}\|a_h\|^2.
		\end{aligned}
		\]
		Combining with \eqref{S-eq:bh_additive_clean_score_drop} gives
		\begin{equation}\label{S-eq:bh_theta_drop}
			\begin{aligned}
				S^x(q_h)-S^x(q_i)
				&=\bigl[S^x(q_h)-S^x(\widehat t_i)\bigr]
				-\bigl[S^x(q_i)-S^x(\widehat t_i)\bigr]\\
				&\ge \rho\|a_h\|^2-\frac{\rho}{2}\|a_h\|^2
				=\frac{\rho}{2}\|a_h\|^2>0.
			\end{aligned}
		\end{equation}
		Thus \(h<i\), and the definition of \(h\) yields
		\[
		h=\min\{j<i:\widehat t_i\in\mathcal B(\tau_j,2B)\}.
		\]
		\begin{enumerate}
			\item[(i)] Suppose that \(h\) is unblocked. By
			Definition~\ref{S-def:blocked_instances_local_noise},
			\(\widehat t_h\in T\cap\mathcal B(\tau_h,\delta)\).
			Since \(\widehat t_i\notin\mathcal B(\tau_h,\delta)\), we have
			\(\widehat t_i\in T\setminus\{\widehat t_h\}\). Applying
			\eqref{S-eq:candidate_remaining_distance} with
			\(s_h=\widehat t_h\) and \(t=\widehat t_i\) gives
			\[
			\|\widehat t_i-\tau_h\|_\infty>B,
			\]
			contradicting \(\widehat t_i\in\mathcal B(\tau_h,2B)\).
			Hence \(h\) is blocked.
			
			\item[(ii)] By \eqref{S-eq:bh_multiplicative_clean_score_drop},
			\eqref{S-eq:bh_representative_clean_score_bounds}, and
			\eqref{S-eq:bh_theta_drop},
			\[
			\begin{aligned}
				(1-\theta)S^x(q_h)
				&=\frac{\min\{\rho,\kappa_{\mathcal V}^2\}}
				{2\kappa_{\mathcal V}^2}\,S^x(q_h)\\
				&\le \frac{\rho}{2}\|a_h\|^2
				\le S^x(q_h)-S^x(q_i).
			\end{aligned}
			\]
			Rearranging gives \(S^x(q_i)\le\theta S^x(q_h)\).
			
			\item[(iii)] Since \(q_i\in\mathcal B(\tau_i,\delta)
			\cap\mathcal B(\widehat t_i,r)\) and
			\(\widehat t_i\in\mathcal B(\tau_h,2B)\),
			\[
			\begin{aligned}
				\|\tau_i-\tau_h\|_\infty
				&\le \|\tau_i-q_i\|_\infty
				+\|q_i-\widehat t_i\|_\infty
				+\|\widehat t_i-\tau_h\|_\infty\\
				&< \frac{\delta}{2}+\frac{r}{2}+B
				=2B+\delta=r.
			\end{aligned}
			\]
		\end{enumerate}
	\end{proof}
	
	We now use Lemma~\ref{S-lem:blocking_candidate} to bound the total
	number of blocked instances by a sum over the indices \(i\) for
	which \(E_i^{\rm loc}\) occurs.
	
	\begin{lemma}[Counting blocked instances]
		\label{S-lem:bh_displacement_charging}
		Let \(D_L^{\rm block}\) be the number of blocked instances
		in the sense of Definition~\ref{S-def:blocked_instances_local_noise},
		and let \(E_l^{\rm loc}\) be the local noise event for instance
		\(l\) defined in \eqref{S-eq:bh_local_noise_event}.
		There exists a finite constant \(C_{\rm group}\), depending only
		on \(d\), \(\kappa_{\mathcal V}\), and \(\rho\), such that
		\begin{equation}\label{S-eq:bh_displacement_deterministic_bound}
			D_L^{\rm block}
			\le C_{\rm group}\sum_{l=1}^N
			\mathbf 1_{E_l^{\rm loc}}
			\left(1+\frac{\|a_l\|^2}{a_L^2}\right)^2.
		\end{equation}
		In particular, if none of the events \(E_l^{\rm loc}\) occurs,
		no instance is blocked.
	\end{lemma}
	
	\begin{proof}
		To bound \(D_L^{\rm block}\), we first show that each blocked
		instance \(i\) can be assigned to an instance \(l\le i\) for
		which \(E_l^{\rm loc}\) occurs. We then bound the number of
		blocked instances assigned to each \(l\).
		If no instance is blocked, the bound is immediate. Otherwise,
		let \(i\) be the smallest index of a blocked instance in the
		ordering fixed in Assumption~\ref{def:localization_property}.
		Then \(E_i^{\rm loc}\) must occur, since otherwise
		Lemma~\ref{S-lem:blocking_candidate} would give a blocked instance
		with index \(h<i\). We therefore take the sequence
		\((i)\), whose last index is \(l=i\).
		Now fix a blocked instance \(i\) and suppose that, for every
		blocked instance with smaller index, we have constructed a
		sequence ending at an instance where the local noise event
		occurs. If \(E_i^{\rm loc}\) occurs, take the sequence \((i)\).
		Otherwise, let
		\[
		h:=\min\{j<i:\widehat t_i\in\mathcal B(\tau_j,2B)\}.
		\]
		By Lemma~\ref{S-lem:blocking_candidate}, the minimum exists and
		instance \(h\) is blocked. By the induction hypothesis, the
		sequence for \(h\) ends at an index \(l\) for which
		\(E_l^{\rm loc}\) occurs. Placing \(i\) before the sequence
		for \(h\) gives a sequence beginning at \(i\) and ending at
		the same \(l\).
		
		For \(1\le l\le N\), let \(\mathcal G_l\) be the set of starting
		indices of the sequences that end at \(l\). Each blocked instance
		starts exactly one sequence, and a sequence can end at \(l\)
		only if \(E_l^{\rm loc}\) occurs. Hence
		\begin{equation}\label{S-eq:blocked_count_by_endpoint}
			D_L^{\rm block}
			=\sum_{l=1}^N|\mathcal G_l|
			=\sum_{l=1}^N\mathbf 1_{E_l^{\rm loc}}|\mathcal G_l|.
		\end{equation}
		Fix \(l\) with \(\mathcal G_l\ne\emptyset\) and set
		\(s=\|a_l\|^2/a_L^2\ge1\), where \(a_L\) is defined in
		\eqref{eq:aL_def}. For \(i\in\mathcal G_l\), write its
		sequence as \(j_0,\ldots,j_m\), where \(j_0=i\), \(j_m=l\),
		and \(m\) is the number of steps. By
		Lemma~\ref{S-lem:blocking_candidate}\textup{(ii)},
		\[
		S^x(q_{j_k})\le\theta S^x(q_{j_{k+1}})
		\qquad (0\le k<m).
		\]
		Iterating and applying \eqref{S-eq:bh_representative_clean_score_bounds}
		at \(i\) and \(l\) gives
		\[
		a_L^2\le\|a_i\|^2\le S^x(q_i)
		\le\theta^m S^x(q_l)
		\le\theta^m\kappa_{\mathcal V}^2\|a_l\|^2.
		\]
		Dividing by \(\theta^m a_L^2\) gives
		\((1/\theta)^m\le\kappa_{\mathcal V}^2s\). Since \(0<\theta<1\)
		by \eqref{S-eq:bh_multiplicative_clean_score_drop}, taking logarithms
		yields
		\begin{equation}\label{S-eq:blocked_sequence_length}
			m\le\frac{\log(\kappa_{\mathcal V}^2s)}{\log(1/\theta)}.
		\end{equation}
		By Lemma~\ref{S-lem:blocking_candidate}\textup{(iii)}, consecutive
		centers are at \(\ell_\infty\)-distance less than \(r\). Applying
		the triangle inequality along the sequence from \(j_0=i\) to
		\(j_m=l\) and using \eqref{S-eq:blocked_sequence_length} gives
		\begin{equation}\label{S-eq:blocked_group_radius}
			\|\tau_i-\tau_l\|_\infty
			\le \sum_{k=0}^{m-1}\|\tau_{j_k}-\tau_{j_{k+1}}\|_\infty
			\le rm
			\le r\,\frac{\log(\kappa_{\mathcal V}^2s)}{\log(1/\theta)}
			=:R_l.
		\end{equation}
		By \eqref{S-eq:blocked_group_radius}, every center \(\tau_i\) with
		\(i\in\mathcal G_l\) lies in the closed cube centered at \(\tau_l\)
		of side length \(2R_l\).
		To count these centers, place a cube of side length \(r/2\)
		around each one. By \eqref{eq:separation}, any two distinct
		centers differ by more than \(r/2\) in at least one coordinate,
		so the cubes are disjoint. For any point \(x\) in the cube
		centered at \(\tau_i\), we have \(\|x-\tau_i\|_\infty\le r/4\).
		Hence, by \eqref{S-eq:blocked_group_radius},
		\[
		\|x-\tau_l\|_\infty
		\le \|x-\tau_i\|_\infty+\|\tau_i-\tau_l\|_\infty
		\le r/4+R_l.
		\]
		The \(|\mathcal G_l|\) disjoint cubes of side length \(r/2\)
		all lie in the cube centered at \(\tau_l\) of side length
		\(2R_l+r/2\). Hence
		\[
		|\mathcal G_l|\left(\frac r2\right)^d
		\le \left(2R_l+\frac r2\right)^d.
		\]
		Dividing by \((r/2)^d\) and substituting \(R_l\) from
		\eqref{S-eq:blocked_group_radius} yields
		\begin{equation}\label{S-eq:blocked_group_size}
			|\mathcal G_l|
			\le \left(1+\frac{4R_l}{r}\right)^d
			= \left(1+4\,\frac{\log(\kappa_{\mathcal V}^2s)}
			{\log(1/\theta)}\right)^d.
		\end{equation}
		To obtain a bound that is valid for
		every \(s\ge1\), set
		\begin{equation}\label{S-eq:blocked_group_constant}
			C_{\rm group}:=\sup_{u\ge1}
			\frac{\bigl[1+4\log(\kappa_{\mathcal V}^2u)/
				\log(1/\theta)\bigr]^d}{(1+u)^2}.
		\end{equation}
		The ratio is continuous on \([1,\infty)\) and tends to zero as
		\(u\to\infty\), since \((1+\log u)^d/u^2\to0\). Thus
		\(C_{\rm group}<\infty\). By \eqref{S-eq:bh_multiplicative_clean_score_drop},
		\(\theta\) depends only on \(\kappa_{\mathcal V}\) and \(\rho\),
		so \(C_{\rm group}\) depends only on \(d\), \(\kappa_{\mathcal V}\),
		and \(\rho\).
		For each \(l\) with \(\mathcal G_l\ne\emptyset\), applying
		\eqref{S-eq:blocked_group_constant} at \(u=s=\|a_l\|^2/a_L^2\)
		in \eqref{S-eq:blocked_group_size} gives
		\[
		|\mathcal G_l|
		\le C_{\rm group}
		\left(1+\frac{\|a_l\|^2}{a_L^2}\right)^2.
		\]
		Substituting into \eqref{S-eq:blocked_count_by_endpoint} proves
		\eqref{S-eq:bh_displacement_deterministic_bound}.
	\end{proof}
	
	Lemma~\ref{S-lem:bh_high_score_true_candidates} gives a selected
	location with observed score at least \(\|a_i\|^2/4\) in the
	detection cube of each instance \(i\) counted by \(G_L\).
	For the BH power bound, we now combine
	\eqref{S-eq:bh_displacement_deterministic_bound} with bounds on the
	probabilities of the local noise events to control \(\E[1-G_L/N]\).
	
	\begin{lemma}[{An upper bound on \(\E[1-G_L/N]\)}]
		\label{S-lem:bh_good_fraction}
		Let \(A_{\rm BH}:=4\eta_{\mathrm{peak}}/\gamma_{\rm BH}\),
		where \(\eta_{\mathrm{peak}}\) and \(\gamma_{\rm BH}\) are
		defined in \eqref{eq:eta_peak_def} and
		\eqref{S-eq:bh_gamma_choice}, respectively.
		If \(a_L^2\ge A_{\rm BH}\), then, for every \(i\),
		\begin{equation}\label{S-eq:bh_local_noise_probability}
			\Prob(E_i^{\rm loc})
			\le 3^d\exp\!\left\{
			-\frac{\gamma_{\rm BH}\|a_i\|^2}{8\eta_{\mathrm{peak}}}
			\right\}.
		\end{equation}
		Moreover,
		\begin{equation}\label{S-eq:bh_good_fraction_L1}
			\E\!\left[1-\frac{G_L}{N}\right]
			\le 3^d\left(1+\frac{16}{e}C_{\rm group}\right)
			\exp\!\left\{
			-\frac{\gamma_{\rm BH}a_L^2}{8\eta_{\mathrm{peak}}}
			\right\},
		\end{equation}
		where \(C_{\rm group}\) is defined in
		\eqref{S-eq:blocked_group_constant}.
	\end{lemma}
	
	\begin{proof}
		The region in \eqref{S-eq:bh_local_noise_event} is contained in
		a cube of side length \(2(B+\delta)\). Since \(\delta\le2B\)
		and \(r=2B+\delta\),
		\[
		2(B+\delta)\le 3(B+\delta/2)=3r/2.
		\]
		Thus the region is contained in the cube centered at \(\tau_i\)
		of side length \(3r/2\). Dividing each side into three intervals
		of length \(r/2\) gives a cover by \(3^d\) translates of the
		comparison cube \(\overline{\mathcal B}(0,r/2)\).
		By stationarity, the definition of \(p\) in
		\eqref{eq:test_value_intro}, and a union bound,
		\[
		\Prob(E_i^{\rm loc})
		\le 3^d p(\gamma_{\rm BH}\|a_i\|^2).
		\]
		Set \(v:=\gamma_{\rm BH}\|a_i\|^2/\eta_{\mathrm{peak}}\).
		By \eqref{eq:aL_def} and the assumption \(a_L^2\ge A_{\rm BH}\),
		\[
		v
		\ge \frac{\gamma_{\rm BH}a_L^2}{\eta_{\mathrm{peak}}}
		\ge \frac{\gamma_{\rm BH}A_{\rm BH}}{\eta_{\mathrm{peak}}}
		=4.
		\] Applying \eqref{S-eq:local_noise_eta_bound} gives
		\[
		p(\gamma_{\rm BH}\|a_i\|^2)
		\le \exp\!\left\{-\tfrac12(\sqrt v-1)^2\right\}
		\le e^{-v/8},
		\]
		where the last inequality uses \(\sqrt v-1\ge\sqrt v/2\).
		This proves \eqref{S-eq:bh_local_noise_probability}.
		By the definition of \(G_L\) in
		Lemma~\ref{S-lem:bh_high_score_true_candidates}, an instance \(i\)
		is counted in \(N-G_L\) if and only if it is blocked or
		\(E_i^{\rm loc}\) occurs. Hence
		\begin{equation}\label{S-eq:bh_good_count_complement}
			N-G_L
			\le D_L^{\rm block}
			+\sum_{i=1}^N\mathbf 1_{E_i^{\rm loc}}.
		\end{equation}
		Set
		\[
		\lambda:=\frac{\gamma_{\rm BH}a_L^2}
		{8\eta_{\mathrm{peak}}}\ge\frac12,
		\qquad
		s_l:=\frac{\|a_l\|^2}{a_L^2}\ge1.
		\]
		Taking expectations in \eqref{S-eq:bh_displacement_deterministic_bound}
		and applying \eqref{S-eq:bh_local_noise_probability} gives
		\[
		\E[D_L^{\rm block}]
		\le C_{\rm group}\sum_{l=1}^N(1+s_l)^2\Prob(E_l^{\rm loc})
		\le 3^d C_{\rm group}\sum_{l=1}^N(1+s_l)^2e^{-\lambda s_l}.
		\]
		For \(s\ge1\) and \(\lambda\ge1/2\),
		\[
		(1+s)^2e^{-\lambda s}
		\le e^{-\lambda}(1+s)^2e^{-(s-1)/2}
		\le \frac{16}{e}e^{-\lambda},
		\]
		since \((1+s)^2e^{-(s-1)/2}\) attains its maximum at \(s=3\).
		Therefore
		\begin{equation}\label{S-eq:bh_expected_blocked_count}
			\E[D_L^{\rm block}]
			\le 3^d C_{\rm group}\frac{16}{e}Ne^{-\lambda}.
		\end{equation}
		Also, \eqref{S-eq:bh_local_noise_probability} and \(s_l\ge1\) give
		\[
		\sum_{l=1}^N\Prob(E_l^{\rm loc})
		\le 3^d\sum_{l=1}^N e^{-\lambda s_l}
		\le 3^dNe^{-\lambda}.
		\]
		Taking expectations in \eqref{S-eq:bh_good_count_complement}
		and using \eqref{S-eq:bh_expected_blocked_count} yields
		\[
		\E\!\left[1-\frac{G_L}{N}\right]
		\le \frac{\E[D_L^{\rm block}]
			+\sum_{l=1}^N\Prob(E_l^{\rm loc})}{N}
		\le 3^d\left(1+\frac{16}{e}C_{\rm group}\right)e^{-\lambda}.
		\]
		Substituting the definition of \(\lambda\) proves
		\eqref{S-eq:bh_good_fraction_L1}.
	\end{proof}
	
	\subsection{Simultaneous recovery and Bonferroni control}
	\label{S-app:simultaneous_recovery}
	
	\begin{proof}[Proof of Proposition~\ref{prop:asymptotic_error_power}
		and Corollary~\ref{cor:simultaneous_recovery_asymptotic}]
		For \(a_L>0\), we identify an event on which every detection
		cube contains a candidate with observed score at least \(a_L^2/4\).
		We then show that, if \(\varepsilon_L\le\alpha\), the Bonferroni
		and BH procedures in Step~\ref{step:multiple_testing} both
		report these candidates on that event.
		Since \(Z_{\max}\ge0\) and Theorem~\ref{thm:continuity_test}
		gives \(\Prob(Z_{\max}=0)=0\), the definition of \(p\) in
		\eqref{eq:test_value_intro} gives \(p(0)=1\). For any \(L\)
		with \(a_L=0\) and any \(\Gamma_{\mathcal V}>0\), it follows that
		\[
		\varepsilon_L
		=M_Lp\!\left(\frac{2a_L^2}{\Gamma_{\mathcal V}}\right)
		=M_Lp(0)=M_L\ge1.
		\]
		The stated FWER bound then reduces to \(\FWER(u_{\rm Bon})\le1\).
		Since \(\alpha<1\), every \(L\) satisfying the condition
		\(\varepsilon_L\le\alpha\) in the simultaneous recovery bound
		\eqref{eq:simultaneous_recovery_finite_bound} also satisfies
		\(a_L>0\).
		Fix \(L\) with \(a_L>0\). By Lemma~\ref{S-lem:representative_scores},
		the set of templates in \eqref{S-eq:kappa_V_def} is nonempty, so
		\(\kappa_{\mathcal V}\) is well defined. With \(\gamma_{\rm BH}\)
		as in \eqref{S-eq:bh_gamma_choice}, take
		\(\Gamma_{\mathcal V}=2/\gamma_{\rm BH}\). The definition of
		\(\varepsilon_L\) in Proposition~\ref{prop:asymptotic_error_power}
		then gives
		\begin{equation}\label{S-eq:global_recovery_remainder}
			\varepsilon_L
			=M_Lp\!\left(\frac{2a_L^2}{\Gamma_{\mathcal V}}\right)
			=M_Lp(\gamma_{\rm BH}a_L^2).
		\end{equation}
		Define
		\[
		\Omega_L:=
		\left\{
		\max_{t\in\overline{\mathcal B}(0,L)} S^z(t)
		\le \gamma_{\rm BH}a_L^2
		\right\}.
		\]
		Applying \eqref{S-eq:growing_noise_union_bound} at
		\(u=\gamma_{\rm BH}a_L^2\) and using
		\eqref{S-eq:global_recovery_remainder} gives
		\begin{equation}\label{S-eq:candidate_noise_event_bound}
			\Prob(\Omega_L^c)
			\le M_Lp(\gamma_{\rm BH}a_L^2)
			=\varepsilon_L.
		\end{equation}
		On \(\Omega_L\), the definition of \(a_L\) in \eqref{eq:aL_def}
		gives, for every \(t\in\overline{\mathcal B}(0,L)\) and every \(i\),
		\[
		S^z(t)\le\gamma_{\rm BH}a_L^2
		\le\gamma_{\rm BH}\|a_i\|^2.
		\]
		Hence none of the local noise events defined in
		\eqref{S-eq:bh_local_noise_event} occurs and using Lemma~\ref{S-lem:bh_displacement_charging}, no instance
		is blocked. Thus \(G_L=N\), and
		Lemma~\ref{S-lem:bh_high_score_true_candidates} gives \(N\) distinct
		candidates satisfying
		\[
		\widehat t_i\in\mathcal B(\tau_i,\delta),
		\qquad
		S^y(\widehat t_i)\ge\frac{\|a_i\|^2}{4}
		\ge\frac{a_L^2}{4},
		\qquad i=1,\ldots,N.
		\]
		Since every detection cube contains a candidate on \(\Omega_L\),
		\eqref{S-eq:false_detections_full_candidates} gives
		\(V(u)\le U_S(u)\) on \(\Omega_L\) for every threshold \(u\).
		At the Bonferroni threshold, \eqref{eq:u_bon_def} and
		Corollary~\ref{S-cor:generalized_inverse_attainment} give
		\(p(u_{\rm Bon})=\alpha/M_L\). Hence
		\eqref{S-eq:noise_box_expectation_identity} yields
		\begin{equation}\label{S-eq:bon_noise_exceedance_mean}
			\E[U_S(u_{\rm Bon})]=M_Lp(u_{\rm Bon})=\alpha.
		\end{equation}
		Using \eqref{S-eq:candidate_noise_event_bound}, we obtain
		\begin{equation}\label{S-eq:bon_fwer_finite}
			\begin{aligned}
				\FWER(u_{\rm Bon})
				&=\Prob\{V(u_{\rm Bon})\ge1\}\\
				&\le\Prob(\Omega_L^c)+\Prob\{U_S(u_{\rm Bon})\ge1\}\\
				&\le\varepsilon_L+\E[U_S(u_{\rm Bon})]
				=\varepsilon_L+\alpha.
			\end{aligned}
		\end{equation}
		Taking the minimum with one proves the stated FWER bound.
		
		Now suppose that \(\varepsilon_L\le\alpha\). On \(\Omega_L\),
		the score bound for each \(\widehat t_i\), the monotonicity
		of \(p\), and \(\gamma_{\rm BH}\le1/4\) give
		\[
		p(S^y(\widehat t_i))
		\le p(a_L^2/4)
		\le p(\gamma_{\rm BH}a_L^2)
		=\frac{\varepsilon_L}{M_L}
		\le\frac{\alpha}{M_L}.
		\]
		By \eqref{eq:u_bon_def} and
		Corollary~\ref{S-cor:generalized_inverse_attainment}, all \(N\)
		candidates are reported by Bonferroni. Since they are distinct,
		the \(N\)th ordered test value satisfies
		\[
		p_{(N)}\le\frac{\alpha}{M_L}\le\frac{N\alpha}{M_L}.
		\]
		Hence \eqref{eq:k_BH_def} gives \(k\ge N\), and the threshold
		definitions \eqref{eq:u_bon_def} and \eqref{eq:u_bh_def}
		give \(u_{\rm BH}\le u_{\rm Bon}\). It follows that
		\begin{equation}\label{S-eq:candidate_reports_on_noise_event}
			\Omega_L\subseteq
			\bigcap_{i=1}^N
			\{T(u)\cap\mathcal B(\tau_i,\delta)\ne\emptyset\},
			\qquad u\in\{u_{\rm Bon},u_{\rm BH}\}.
		\end{equation}
		Taking probabilities and using \eqref{S-eq:candidate_noise_event_bound}
		proves \eqref{eq:simultaneous_recovery_finite_bound}.
		
		To prove the asymptotic conclusions of
		Corollary~\ref{cor:simultaneous_recovery_asymptotic}, it suffices
		to show that \(\varepsilon_L\to0\). Indeed, \eqref{S-eq:bon_fwer_finite}
		then gives \(\limsup_{L\to\infty}\FWER(u_{\rm Bon})\le\alpha\).
		Moreover, \(\varepsilon_L\le\alpha\) for all sufficiently large \(L\),
		so \eqref{eq:simultaneous_recovery_finite_bound} gives simultaneous
		recovery with probability tending to one at both thresholds.
		We prove that \(\varepsilon_L\to0\) under the following condition,
		which is implied by the strength condition
		\eqref{eq:object_strength_condition} in
		Corollary~\ref{cor:simultaneous_recovery_asymptotic}:
		\begin{equation}\label{S-eq:global_recovery_sigma_condition}
			\liminf_{L\to\infty}
			\frac{a_L^2}{\sigma_{\mathcal V}^2\log M_L}
			>\Gamma_{\mathcal V},
		\end{equation}
		where \(\sigma_{\mathcal V}^2\) is defined in
		\eqref{S-eq:projected_noise_max_variance}.
		Since
		\(0<\sigma_{\mathcal V}^2\le\eta_{\mathrm{peak}}\) by
		\eqref{S-eq:noise_maximum_parameters}, for all sufficiently large \(L\),
		\[
		\frac{a_L^2}{\sigma_{\mathcal V}^2\log M_L}
		\ge
		\frac{a_L^2}{\eta_{\mathrm{peak}}\log M_L}.
		\]
		Thus the strength condition \eqref{eq:object_strength_condition}
		in Corollary~\ref{cor:simultaneous_recovery_asymptotic} implies
		\eqref{S-eq:global_recovery_sigma_condition}.
		Since \(\Gamma_{\mathcal V}=2/\gamma_{\rm BH}\),
		condition~\eqref{S-eq:global_recovery_sigma_condition} allows us
		to choose \(c>0\) such that
		\[
		1+c<
		\liminf_{L\to\infty}
		\frac{\gamma_{\rm BH}a_L^2}{2\sigma_{\mathcal V}^2\log M_L}.
		\]
		By \eqref{S-eq:grid_count_asymptotic}, \(M_L\to\infty\). Applying
		\eqref{S-eq:noise_quantile_asymptotic} at the level \(M_L^{-(1+c)}\)
		gives
		\[
		\frac{u_{M_L^{-(1+c)}}}{2\sigma_{\mathcal V}^2\log M_L}
		=(1+c)\frac{u_{M_L^{-(1+c)}}}
		{2\sigma_{\mathcal V}^2\log(M_L^{1+c})}
		\longrightarrow 1+c.
		\]
		Thus \(u_{M_L^{-(1+c)}}\le\gamma_{\rm BH}a_L^2\) for all
		sufficiently large \(L\). By \eqref{S-eq:global_recovery_remainder},
		monotonicity of \(p\), and
		Corollary~\ref{S-cor:generalized_inverse_attainment},
		\begin{equation}\label{S-eq:recovery_remainder_vanishes}
			\varepsilon_L
			=M_Lp(\gamma_{\rm BH}a_L^2)
			\le M_Lp\bigl(u_{M_L^{-(1+c)}}\bigr)
			=M_L^{-c}\longrightarrow0.
		\end{equation}	Since \(\varepsilon_L\le\alpha\) for all sufficiently large \(L\),
		the asymptotic conclusions follow from
		\eqref{eq:simultaneous_recovery_finite_bound} and
		\eqref{S-eq:bon_fwer_finite}.
		Finally, if \(P_{\mathcal V}u_\ell=0\) for every template
		\(u_\ell\), then \eqref{eq:x} and
		Lemma~\ref{S-lem:proj_rot_appendix} give \(a_L=0\) for every \(L\).
	\end{proof}
	\begin{corollary}[Bonferroni FDR under simultaneous recovery]
		\label{S-cor:bon_fdr_recovery}
		Under the assumptions of
		Proposition~\ref{prop:asymptotic_error_power}, if
		\(\varepsilon_L\le\alpha\), then
		\begin{equation}\label{S-eq:bon_fdr_recovery_bound}
			\FDR(u_{\rm Bon})
			\le\min\left\{1,\frac{\alpha}{N+1}+\varepsilon_L\right\}.
		\end{equation}
		If \(N\to\infty\) and either \eqref{eq:object_strength_condition}
		or \eqref{S-eq:global_recovery_sigma_condition} holds, then
		\(\FDR(u_{\rm Bon})\to0\).
	\end{corollary}
	\begin{proof}
		Fix \(L\) with \(\varepsilon_L\le\alpha\) and set \(u=u_{\rm Bon}\).
		On \(\Omega_L\), every detection cube contains a report by
		\eqref{S-eq:candidate_reports_on_noise_event}. Hence \(D(u)=0\),
		and \eqref{S-eq:bh_W_equals_N_minus_D} gives \(W(u)=N\).
		Moreover, \eqref{S-eq:false_detections_full_candidates} gives
		\(V(u)\le U_S(u)\) on \(\Omega_L\).
		If \(V(u)\ge1\), then \(N+V(u)\ge N+1\). If \(V(u)=0\),
		the false discovery proportion is zero. Thus, on \(\Omega_L\),
		\[
		\frac{V(u)}{\max\{W(u)+V(u),1\}}
		=\frac{V(u)}{N+V(u)}
		\le\frac{U_S(u)}{N+1}.
		\]
		By \eqref{def:fdr}, \eqref{S-eq:bon_noise_exceedance_mean}, and
		\eqref{S-eq:candidate_noise_event_bound},
		\[
		\begin{aligned}
			\FDR(u)
			&=\E\!\left[
			\frac{V(u)\mathbf 1_{\Omega_L}}{\max\{W(u)+V(u),1\}}
			\right]
			+\E\!\left[
			\frac{V(u)\mathbf 1_{\Omega_L^c}}{\max\{W(u)+V(u),1\}}
			\right]\\
			&\le \min\left\{1,\,
			\frac{\E[U_S(u)\mathbf 1_{\Omega_L}]}{N+1}
			+\Prob(\Omega_L^c)\right\}\\
			&\le \min\left\{1,\frac{\alpha}{N+1}+\varepsilon_L\right\}.
		\end{aligned}
		\]
		Suppose \(N\to\infty\) and either \eqref{eq:object_strength_condition}
		or \eqref{S-eq:global_recovery_sigma_condition} holds.
		By \eqref{S-eq:recovery_remainder_vanishes}, \(\varepsilon_L\to0\),
		so \(\varepsilon_L\le\alpha\) for all sufficiently large \(L\).
		For such \(L\),
		\[
		0\le\FDR(u_{\rm Bon})
		\le\frac{\alpha}{N+1}+\varepsilon_L
		\longrightarrow0.
		\]
	\end{proof}
	
	\subsection{Power and FDR at the BH threshold}
	\label{S-app:bh_density_strength}
	
	\subsubsection{BH power}
	\label{S-subsubsec:bh_power_bounds}
	
	By Lemma~\ref{S-lem:bh_high_score_true_candidates}, at least \(G_L\)
	distinct detection cubes contain a candidate with score at least
	\(a_L^2/4\). We give a sufficient condition for the BH procedure in
	\eqref{eq:k_BH_def} and \eqref{eq:u_bh_def} to report
	all \(G_L\) candidates identified by the lemma.
	
	\begin{lemma}[BH power for a fixed image size]
		\label{S-lem:bh_finite_power_bound}
		Under the assumptions of Proposition~\ref{prop:bh_density_strength},
		suppose \(a_L>0\) and set
		\begin{equation}\label{S-eq:bh_tail_density_ratio}
			q_L:=\frac{M_L}{\alpha N}p(a_L^2/4).
		\end{equation}
		If \(q_L<1\), then
		\begin{equation}\label{S-eq:bh_power_q_bound}
			1-\Power(u_{\rm BH})\le\frac{\E[1-G_L/N]}{1-q_L}.
		\end{equation}
		If also \(a_L^2\ge A_{\rm BH}\) and \(q_L\le1/2\), where \(A_{\rm BH}\)
		is defined in Lemma~\ref{S-lem:bh_good_fraction}, then
		\begin{equation}\label{S-eq:bh_power_exponential_bound}
			1-\Power(u_{\rm BH})
			\le 2\cdot3^d\left(1+\frac{16}{e}C_{\rm group}\right)
			\exp\!\left\{-\frac{\gamma_{\rm BH}a_L^2}{8\eta_{\mathrm{peak}}}\right\}.
		\end{equation}
	\end{lemma}
	
	\begin{proof}
		Suppose \(q_L<1\). Recall that \(D(u_{\rm BH})\) counts the
		instances whose detection cubes contain no report at the BH
		threshold.
		On the event \(\{G_L/N>q_L\}\), we have \(G_L\ge1\), since
		\(q_L\ge0\) and \(G_L\) is an integer. By
		Lemma~\ref{S-lem:bh_high_score_true_candidates}, there are \(G_L\)
		candidates in distinct detection cubes, each with observed score
		at least \(a_L^2/4\). Since \(p\) is nonincreasing, each of
		these candidates has test value at most \(p(a_L^2/4)\). Thus,
		by the ordering of the test values and \eqref{S-eq:bh_tail_density_ratio},
		\[
		p_{(G_L)}
		\le p(a_L^2/4)
		= \frac{\alpha Nq_L}{M_L}
		< \frac{\alpha G_L}{M_L}.
		\]
		By \eqref{eq:k_BH_def}, \(k\ge G_L\), so
		\[
		p(a_L^2/4)<\frac{\alpha G_L}{M_L}
		\le\frac{\alpha k}{M_L}.
		\]
		Thus all \(G_L\) candidates meet the BH reporting criterion by
		\eqref{eq:u_bh_def} and
		Corollary~\ref{S-cor:generalized_inverse_attainment}. Since they
		lie in distinct detection cubes, \(D(u_{\rm BH})\le N-G_L\).
		Using \(0\le q_L<1\), we obtain
		\[
		(1-q_L)\frac{D(u_{\rm BH})}{N}
		\le \frac{D(u_{\rm BH})}{N}
		\le 1-\frac{G_L}{N}.
		\]
		If \(G_L/N\le q_L\), then \(D(u_{\rm BH})\le N\) gives
		\[
		(1-q_L)\frac{D(u_{\rm BH})}{N}
		\le 1-q_L
		\le 1-\frac{G_L}{N}.
		\]
		The same inequality therefore holds in both cases. Taking
		expectations, dividing by \(1-q_L>0\), and using
		\eqref{S-eq:bh_power_missed_identity} gives
		\[
		1-\Power(u_{\rm BH})
		=\E\!\left[\frac{D(u_{\rm BH})}{N}\right]
		\le\frac{\E[1-G_L/N]}{1-q_L},
		\]
		which proves \eqref{S-eq:bh_power_q_bound}.
		If \(q_L\le1/2\), then \((1-q_L)^{-1}\le2\). For
		\(a_L^2\ge A_{\rm BH}\), substituting the bound on
		\(\E[1-G_L/N]\) from \eqref{S-eq:bh_good_fraction_L1} into
		\eqref{S-eq:bh_power_q_bound} gives \eqref{S-eq:bh_power_exponential_bound}.
	\end{proof}
	
	\begin{proof}[Proof of Proposition~\ref{prop:bh_density_strength}]
		Recall from \eqref{S-eq:bh_tail_density_ratio} that
		\[
		q_L=\frac{p(a_L^2/4)}{\alpha N/M_L},
		\]
		where \(a_L^2/4\) is the lower bound on the candidate scores in
		\eqref{S-eq:bh_candidate_score_level}.
		By Lemma~\ref{S-lem:bh_finite_power_bound}, the power bound
		\eqref{S-eq:bh_power_exponential_bound} holds when
		\(a_L^2\ge A_{\rm BH}\) and \(q_L\le1/2\).
		The first inequality is assumed in \eqref{eq:bh_finite_power_conditions}.
		By the formula for \(q_L\), the second is equivalent to
		\[
		p(a_L^2/4)\le\frac{\alpha N}{2M_L}.
		\]
		We derive this bound from the two conditions in
		\eqref{eq:bh_finite_power_conditions} with \(C_{\rm BH}=32\).
		Suppose first that at least one template has nonzero projection
		onto \(\mathcal V\), so that \(\kappa_{\mathcal V}\) and
		\(\gamma_{\rm BH}\) are defined by \eqref{S-eq:kappa_V_def}
		and \eqref{S-eq:bh_gamma_choice}. Take \(A_{\rm BH}\) as in
		Lemma~\ref{S-lem:bh_good_fraction}. Since \(\gamma_{\rm BH}\le1/4\),
		the first condition in \eqref{eq:bh_finite_power_conditions}
		gives
		\[
		a_L^2\ge A_{\rm BH}
		=\frac{4\eta_{\mathrm{peak}}}{\gamma_{\rm BH}}
		\ge16\eta_{\mathrm{peak}}.
		\]
		In particular, \(a_L^2/4\ge\eta_{\mathrm{peak}}\), so the
		tail bound \eqref{S-eq:local_noise_sigma_bound} applies at
		\(u=a_L^2/4\). Moreover,
		\(\sqrt{\eta_{\mathrm{peak}}}\le a_L/4\), and hence
		\(a_L/2-\sqrt{\eta_{\mathrm{peak}}}\ge a_L/4\).
		Using also \(\sigma_{\mathcal V}^2\le\eta_{\mathrm{peak}}\)
		from Lemma~\ref{S-lem:noise_maxima_fixed_growing}, we obtain
		\begin{equation}\label{S-eq:bh_candidate_tail_bound}
			\begin{aligned}
				p(a_L^2/4)
				&\le \exp\!\left\{
				-\frac{(a_L/2-\sqrt{\eta_{\mathrm{peak}}})^2}
				{2\sigma_{\mathcal V}^2}
				\right\}\\
				&\le \exp\!\left\{-\frac{a_L^2}{32\sigma_{\mathcal V}^2}\right\}
				\le \exp\!\left\{-\frac{a_L^2}{32\eta_{\mathrm{peak}}}\right\}.
			\end{aligned}
		\end{equation}
		Now use the second condition in
		\eqref{eq:bh_finite_power_conditions}. With \(C_{\rm BH}=32\),
		it gives
		\[
		\frac{a_L^2}{32\eta_{\mathrm{peak}}}
		\ge \log\!\left(\frac{2M_L}{\alpha N}\right).
		\]
		Substituting into \eqref{S-eq:bh_candidate_tail_bound} yields
		\[
		p(a_L^2/4)
		\le \exp\!\left\{-\log\!\left(\frac{2M_L}{\alpha N}\right)\right\}
		=\frac{\alpha N}{2M_L}.
		\]
		Thus \(q_L\le1/2\), and \eqref{S-eq:bh_power_exponential_bound}
		gives \eqref{eq:bh_average_power_rate} with
		\begin{equation}\label{S-eq:bh_power_explicit_constants}
			A_{\rm BH}=\frac{4\eta_{\mathrm{peak}}}{\gamma_{\rm BH}},
			\qquad
			C=2\cdot3^d\left(1+\frac{16}{e}C_{\rm group}\right),
			\qquad
			c=\frac{\gamma_{\rm BH}}{8}.
		\end{equation}
		By \eqref{eq:eta_peak_def}, \eqref{S-eq:kappa_V_def},
		\eqref{S-eq:bh_gamma_choice}, and \eqref{S-eq:blocked_group_constant},
		these constants depend only on the parameters listed in the
		proposition.
		If all template projections vanish, \eqref{eq:x} and
		Lemma~\ref{S-lem:proj_rot_appendix} give \(a_L=0\) for every \(L\).
	\end{proof}
	
	\begin{proof}[Proof of Corollary~\ref{cor:bh_power_asymptotic}]
		Assume \eqref{eq:bh_log_density_balance}. Since
		\(a_L^2\to\infty\) and \(A_{\rm BH}\) is fixed,
		\(a_L^2\ge A_{\rm BH}\) for all sufficiently large \(L\).
		The logarithmic inequality in \eqref{eq:bh_log_density_balance}
		is precisely the second condition in
		\eqref{eq:bh_finite_power_conditions}. Thus
		Proposition~\ref{prop:bh_density_strength} gives, for all
		sufficiently large \(L\),
		\[
		0\le 1-\Power(u_{\rm BH})
		\le C\exp\!\left\{-c\frac{a_L^2}{\eta_{\mathrm{peak}}}\right\}
		\longrightarrow0,
		\]
		since \(0<\eta_{\mathrm{peak}}<\infty\) is independent of \(L\).
		Finally, suppose \eqref{eq:positive_density_weak_strength}
		holds. Choose \(0<\beta<\liminf_{L\to\infty}N/M_L\).
		Since \(N/M_L\ge\beta\) for all sufficiently large \(L\) and
		\(a_L^2/\eta_{\mathrm{peak}}\to\infty\), we have, for all
		sufficiently large \(L\),
		\[
		C_{\rm BH}\log\!\left(\frac{2M_L}{\alpha N}\right)
		\le C_{\rm BH}\log\!\left(\frac{2}{\alpha\beta}\right)
		\le \frac{a_L^2}{\eta_{\mathrm{peak}}}.
		\]
		Together with \(a_L^2\to\infty\), this verifies
		\eqref{eq:bh_log_density_balance}. Hence
		\(\Power(u_{\rm BH})\to1\).
	\end{proof}

	\begin{remark}[BH power under a variance-based strength condition]
		\label{S-rem:bh_sigma_refinement}
		The first inequality in \eqref{S-eq:bh_candidate_tail_bound} gives
		\(q_L\le1/2\) whenever \(a_L^2\ge A_{\rm BH}\) and
		\[
		\frac{a_L^2}{\sigma_{\mathcal V}^2}
		\ge32\log\!\left(\frac{2M_L}{\alpha N}\right).
		\]
		Thus the logarithmic conditions in \eqref{eq:bh_finite_power_conditions}
		and \eqref{eq:bh_log_density_balance} also suffice with
		\(\sigma_{\mathcal V}^2\) in the denominator. Equation
		\eqref{S-eq:bh_power_exponential_bound} gives the same power bound
		in terms of \(a_L^2/\eta_{\mathrm{peak}}\).
	\end{remark}
	
	\subsubsection{BH false discovery rate}
	\label{S-subsubsec:bh_false_discovery_rate}
	
	By \eqref{S-eq:bh_false_candidates_missed}, the number of false reports
	at the BH threshold satisfies
	\[
	V(u_{\rm BH})\le U_S(u_{\rm BH})+2^dD(u_{\rm BH}).
	\]
	Recall that the search domain is covered by \(M_L\) grid cubes
	\(Q_m\) of side length \(r/2\), as in \eqref{eq:Qm}.
	By \eqref{S-eq:noise_exceedance_count}, \(U_S(u_{\rm BH})\) counts
	those cubes for which \(\max_{t\in\overline{Q_m}}S^z(t)\ge u_{\rm BH}\).
	The count \(D(u_{\rm BH})\) is the number of instances whose
	detection cubes contain no BH report. Proposition~\ref{prop:bh_fdr_control} assumes
	that \(\Power(u_{\rm BH})\to1\), so \eqref{S-eq:bh_power_missed_identity}
	gives \(\E[D(u_{\rm BH})/N]\to0\). To bound the false discovery
	proportion, we compare \(U_S(u_{\rm BH})\) with the number of
	BH reports.
	
Fix \(0<\alpha<1\). Recall that \(m_L=|T|\) is
the number of candidate locations from
Step~\ref{item:candidate_peaks}.
Let \(k\) be the BH index in~\eqref{eq:k_BH_def}.
If \(k\ge1\), then
	\(p_{(k)}\le\alpha k/M_L\). If also \(k<m_L\), maximality of \(k\)
	gives
	\[
	p_{(k+1)}>\frac{\alpha(k+1)}{M_L}>
	\frac{\alpha k}{M_L}.
	\]
	Thus the candidates with the \(k\) smallest test values are exactly
	those with test value at most \(\alpha k/M_L\). By
	\eqref{eq:u_bh_def} and Corollary~\ref{S-cor:generalized_inverse_attainment},
	BH reports precisely these \(k\) candidates.
	Using
	\eqref{S-eq:bh_W_equals_N_minus_D}, we obtain
	\begin{equation}\label{S-eq:bh_threshold_report_identity}
		p(u_{\rm BH})=\frac{\alpha k}{M_L}\quad(k\ge1),
		\qquad
		k=|T(u_{\rm BH})|=N-D(u_{\rm BH})+V(u_{\rm BH}).
	\end{equation}
	The count identity also holds when \(k=0\), since the procedure then
	reports no candidates. By candidate separation, each grid cube in
	\eqref{eq:Qm} contains at most one candidate, so \(k\le m_L\le M_L\).
	Likewise, \eqref{eq:separation} places distinct instance centers
	in distinct grid cubes, so \(N\le M_L\).
	Write \(D=D(u_{\rm BH})\) and \(V=V(u_{\rm BH})\). Suppose that
	\[
	U_S(u_{\rm BH})\le Ak+E
	\]
	for some \(0\le A<1\) and \(E\ge0\). Then
	\eqref{S-eq:bh_false_candidates_missed} gives \(V\le Ak+E+2^dD\).
	Substituting into \(k=N-D+V\) from
	\eqref{S-eq:bh_threshold_report_identity} gives
	\begin{equation}\label{S-eq:report_count_algebra}
		(1-A)k\le N+E+(2^d-1)D.
	\end{equation}
	For \(k\ge1\), we have \(V/k=1-(N-D)/k\). Since \(N-D\ge0\),
	\eqref{S-eq:report_count_algebra} gives
	\[
	\frac{V}{k}
	\le 1-\frac{(1-A)(N-D)}{N+E+(2^d-1)D}
	= A+(1-A)\frac{E+2^dD}{N+E+(2^d-1)D}.
	\]
	For \(k=0\), the false discovery proportion is zero. Thus, for every
	\(k\),
	\begin{equation}\label{S-eq:report_proportion_algebra}
		\begin{aligned}
			\frac{V}{\max\{k,1\}}
			&\le A+(1-A)\frac{E+2^dD}{N+E+(2^d-1)D}\\
			&\le A+(1-A)\frac{E}{N+E}
			+2^d(1-A)\frac{D}{N},
		\end{aligned}
	\end{equation}
	where the last inequality uses \(N+E+(2^d-1)D\ge N+E\ge N\).
	
	To bound \(U_S(u_{\rm BH})\) in terms of \(k\), consider the
	fraction of grid cubes whose maximum noise score is at least \(u\),
	\begin{equation}\label{S-eq:empirical_noise_tail}
		\widehat p_L(u):=\frac{U_S(u)}{M_L},
		\qquad u\in\R.
	\end{equation}
	For each fixed \(u\), \eqref{S-eq:noise_box_expectation_identity}
	gives \(\E[\widehat p_L(u)]=p(u)\), and
	Assumption~\ref{ass:noise_exceedance_variance} gives
	\begin{equation}\label{S-eq:normalized_noise_exceedance_variance}
		\Var(\widehat p_L(u))\le\frac{C_{\rm var}}{M_L}p(u).
	\end{equation}
	
	\begin{lemma}[Uniform upper bound on \(\widehat p_L(u)/p(u)\)]
		\label{S-lem:uniform_relative_noise_tail}
		Under Assumptions~\ref{ass:noise_reg_main}
		and~\ref{ass:noise_exceedance_variance}, fix \(L>0\) and
		\(t_0,\epsilon\in(0,1)\). The event
		\begin{equation}\label{S-eq:relative_noise_tail_event}
			\mathcal G_L(t_0,\epsilon)
			:=\bigl\{\widehat p_L(u)\le(1+\epsilon)p(u)
			\text{ for every }u\in\R\text{ with }p(u)\ge t_0\bigr\}
		\end{equation}
		is measurable and satisfies
		\begin{equation}\label{S-eq:bh_relative_tail_probability}
			\Prob\bigl(\mathcal G_L(t_0,\epsilon)^c\bigr)
			\le\frac{27C_{\rm var}}{M_L\epsilon^3t_0}.
		\end{equation}
	\end{lemma}
	
	\begin{proof}
		By \eqref{S-eq:noise_exceedance_count} and
		\eqref{S-eq:empirical_noise_tail}, \(\widehat p_L\) is nonincreasing
		and left-continuous. The tail function \(p\) is nonincreasing
		and is continuous by Theorem~\ref{thm:continuity_test}.
		Suppose that \(\widehat p_L(q)\le(1+\epsilon)p(q)\) for every
		\(q\in\mathbb Q\) with \(p(q)\ge t_0\). For any \(u\in\R\) with
		\(p(u)\ge t_0\), take \(q_n\in\mathbb Q\) with \(q_n\uparrow u\).
		Then \(p(q_n)\ge p(u)\ge t_0\), and left-continuity gives
		\[
		\widehat p_L(u)
		=\lim_{n\to\infty}\widehat p_L(q_n)
		\le (1+\epsilon)\lim_{n\to\infty}p(q_n)
		=(1+\epsilon)p(u).
		\]
		Hence
		\[
		\mathcal G_L(t_0,\epsilon)
		=\bigcap_{\substack{q\in\mathbb Q\\p(q)\ge t_0}}
		\{\widehat p_L(q)\le(1+\epsilon)p(q)\},
		\]
		a countable intersection of measurable events.
		
		Set
		\[
		t_j:=t_0(1+\epsilon/2)^j,\qquad j=0,1,\ldots.
		\]
		For each \(j\) with \(t_j<1\), the corresponding score threshold is
		\[
		u_{t_j}=\inf\{v\in\R:p(v)\le t_j\}
		\]
		by \eqref{eq:uc_def_inline}, and
		\(p(u_{t_j})=t_j\) by Corollary~\ref{S-cor:generalized_inverse_attainment}.
		Since \(u_{t_j}\) is deterministic, \eqref{S-eq:noise_box_expectation_identity}
		and \eqref{S-eq:normalized_noise_exceedance_variance} give
		\[
		\E[\widehat p_L(u_{t_j})]=t_j,
		\qquad
		\Var(\widehat p_L(u_{t_j}))\le\frac{C_{\rm var}t_j}{M_L}.
		\]
		Chebyshev's inequality therefore yields
		\begin{equation}\label{S-eq:bh_noise_tail_at_fixed_level}
			\Prob\!\left\{
			\widehat p_L(u_{t_j})-t_j>\frac{\epsilon t_j}{3}
			\right\}
			\le \frac{9\Var(\widehat p_L(u_{t_j}))}{\epsilon^2t_j^2}
			\le \frac{9C_{\rm var}}{M_L\epsilon^2t_j}.
		\end{equation}
		On \(\mathcal G_L(t_0,\epsilon)^c\), there is a threshold \(u\)
		such that
		\[
		p(u)\ge t_0,
		\qquad
		(1+\epsilon)p(u)<\widehat p_L(u)\le1.
		\]
		Choose the smallest \(j\ge0\) with \(t_j\ge p(u)\). If \(j=0\),
		then \(t_j=p(u)=t_0\). If \(j\ge1\), then \(t_{j-1}<p(u)\) and
		\(t_j=(1+\epsilon/2)t_{j-1}\). In both cases,
		\[
		p(u)\le t_j\le(1+\epsilon/2)p(u)
		<\frac{1+\epsilon/2}{1+\epsilon}<1.
		\]
		Since \(p(u)\le t_j\), the definition of \(u_{t_j}\) gives
		\(u_{t_j}\le u\). Monotonicity of \(\widehat p_L\) now gives
		\[
		\widehat p_L(u_{t_j})
		\ge\widehat p_L(u)
		>(1+\epsilon)p(u)
		\ge\frac{1+\epsilon}{1+\epsilon/2}\,t_j.
		\]
		Subtracting \(t_j\) and using \(0<\epsilon<1\), we obtain
		\[
		\widehat p_L(u_{t_j})-t_j
		>\frac{\epsilon}{2+\epsilon}\,t_j
		\ge\frac{\epsilon t_j}{3}.
		\]
		Thus
		\[
		\mathcal G_L(t_0,\epsilon)^c
		\subseteq
		\bigcup_{j:t_j<1}
		\left\{\widehat p_L(u_{t_j})-t_j>\frac{\epsilon t_j}{3}\right\}.
		\]
		A union bound and \eqref{S-eq:bh_noise_tail_at_fixed_level} give
		\[
		\begin{aligned}
			\Prob\bigl(\mathcal G_L(t_0,\epsilon)^c\bigr)
			&\le \frac{9C_{\rm var}}{M_L\epsilon^2}
			\sum_{j:t_j<1}\frac1{t_j}\\
			&\le \frac{9C_{\rm var}}{M_L\epsilon^2t_0}
			\sum_{j=0}^{\infty}(1+\epsilon/2)^{-j}\\
			&= \frac{9C_{\rm var}(2+\epsilon)}{M_L\epsilon^3t_0}
			\le \frac{27C_{\rm var}}{M_L\epsilon^3t_0},
		\end{aligned}
		\]
		which proves \eqref{S-eq:bh_relative_tail_probability}.
	\end{proof}
	
	If \(k\ge1\) and \(\alpha k/M_L\ge t_0\), then
	\(p(u_{\rm BH})=\alpha k/M_L\ge t_0\) by
	\eqref{S-eq:bh_threshold_report_identity}. On
	\(\mathcal G_L(t_0,\epsilon)\), applying
	\eqref{S-eq:relative_noise_tail_event} at \(u=u_{\rm BH}\) gives
	\begin{equation}\label{S-eq:bh_noise_at_random_threshold}
		U_S(u_{\rm BH})
		=M_L\widehat p_L(u_{\rm BH})
		\le M_L(1+\epsilon)p(u_{\rm BH})
		=\alpha(1+\epsilon)k.
	\end{equation}
	
	\begin{lemma}[BH FDR under the variance bound]
		\label{S-lem:bh_fdr_relative_noise_bound}
		Fix \(0<\alpha<1\). Under Assumptions~\ref{assumptions},
		\ref{ass:noise_reg_main}, and~\ref{ass:noise_exceedance_variance},
		for every \(L>0\) and \(0<\epsilon<1\),
		\begin{equation}\label{S-eq:bh_fdr_relative_noise_bound}
			\FDR(u_{\rm BH})
			\le \alpha(1+\epsilon)
			+2^d\bigl[1-\Power(u_{\rm BH})\bigr]
			+\frac{54C_{\rm var}}{\alpha\epsilon^3N},
		\end{equation}
		where \(C_{\rm var}\) is the constant in
		\eqref{eq:noise_exceedance_variance}.
	\end{lemma}
	
	\begin{proof}
		Write \(D=D(u_{\rm BH})\) and \(V=V(u_{\rm BH})\).
		If \(\alpha(1+\epsilon)\ge1\), the bound follows from
		\(\FDR(u_{\rm BH})\le1\). Suppose \(\alpha(1+\epsilon)<1\).
		When \(D\le N/2\), the report identity
		\eqref{S-eq:bh_threshold_report_identity} gives
		\[
		k=N-D+V\ge N-D\ge N/2.
		\]
		In particular, \(k\ge1\) and \(\alpha k/M_L\ge\alpha N/(2M_L)\).
		Set \(t_0:=\alpha N/(2M_L)\), which belongs to \((0,1)\) since
		\(1\le N\le M_L\). On
		\(\{D\le N/2\}\cap\mathcal G_L(t_0,\epsilon)\),
		\eqref{S-eq:bh_noise_at_random_threshold} and
		\eqref{S-eq:bh_false_candidates_missed} give
		\[
		U_S(u_{\rm BH})\le\alpha(1+\epsilon)k,
		\qquad
		V\le\alpha(1+\epsilon)k+2^dD.
		\]
		Applying \eqref{S-eq:report_proportion_algebra} with
		\(A=\alpha(1+\epsilon)\) and \(E=0\), we obtain
		\[
		\frac Vk
		\le \alpha(1+\epsilon)
		+2^d[1-\alpha(1+\epsilon)]\frac DN
		\le \alpha(1+\epsilon)+2^d\frac DN.
		\]
		When \(D>N/2\), we have \(2^dD/N\ge1\), so the same upper
		bound holds for \(V/\max\{k,1\}\). On
		\(\mathcal G_L(t_0,\epsilon)^c\), the false discovery proportion
		is at most one. Hence
		\[
		\frac{V}{\max\{k,1\}}
		\le \alpha(1+\epsilon)+2^d\frac DN
		+\mathbf 1_{\mathcal G_L(t_0,\epsilon)^c}.
		\]
		Taking expectations and using \eqref{def:fdr} and
		\eqref{S-eq:bh_power_missed_identity} gives
		\[
		\FDR(u_{\rm BH})
		\le \alpha(1+\epsilon)
		+2^d[1-\Power(u_{\rm BH})]
		+\Prob\bigl(\mathcal G_L(t_0,\epsilon)^c\bigr).
		\]
		Finally, \eqref{S-eq:bh_relative_tail_probability} and the choice
		\(t_0=\alpha N/(2M_L)\) give
		\[
		\Prob\bigl(\mathcal G_L(t_0,\epsilon)^c\bigr)
		\le \frac{27C_{\rm var}}{M_L\epsilon^3t_0}
		= \frac{54C_{\rm var}}{\alpha\epsilon^3N},
		\]
		proving \eqref{S-eq:bh_fdr_relative_noise_bound}.
	\end{proof}
	
	For the second assertion of Proposition~\ref{prop:bh_fdr_control},
	we bound the BH FDR in terms of the maximum error of \(\widehat p_L\) defined in \eqref{S-eq:empirical_noise_tail}.
	Set
	\begin{equation}\label{S-eq:bh_uniform_noise_tail_error}
		\Delta_L:=\sup_{u\in\R}|\widehat p_L(u)-p(u)|.
	\end{equation}
	Since \(\widehat p_L\) is left-continuous and \(p\) is continuous
	by Theorem~\ref{thm:continuity_test}, the supremum in
	\eqref{S-eq:bh_uniform_noise_tail_error} equals the supremum over
	\(u\in\mathbb Q\). Hence \(\Delta_L\) is measurable.
	The following lemma gives a bound for each \(L\).
	Lemma~\ref{S-lem:bh_uniform_noise_box_lln} will then show that
	\(\E[\Delta_L]\to0\) when the noise field is ergodic.
	
	\begin{lemma}[BH FDR in terms of the uniform noise tail error]
		\label{S-lem:bh_fdr_general_bound}
		Fix \(0<\alpha<1\). Under Assumptions~\ref{assumptions}
		and~\ref{ass:noise_reg_main}, let \(\Delta_L\) be as in
		\eqref{S-eq:bh_uniform_noise_tail_error} and set
		\begin{equation}\label{S-eq:betaL}
			\beta_L:=\frac{N}{M_L}.
		\end{equation}
		For every \(L>0\),
		\begin{equation}\label{S-eq:bh_fdr_finite_bound}
			\FDR(u_{\rm BH})
			\le \alpha
			+(1-\alpha)\E\!\left[
			\frac{\Delta_L}{\beta_L+\Delta_L}
			\right]
			+2^d(1-\alpha)\bigl[1-\Power(u_{\rm BH})\bigr].
		\end{equation}
		If \(\Power(u_{\rm BH})\to1\) and \(\Delta_L/\beta_L\to0\)
		in probability, then
		\[
		\limsup_{L\to\infty}\FDR(u_{\rm BH})\le\alpha.
		\]
	\end{lemma}
	
	\begin{proof}
		Write \(D=D(u_{\rm BH})\) and \(V=V(u_{\rm BH})\).
		For \(k\ge1\), the definition of \(\Delta_L\) applies at
		\(u_{\rm BH}\), and \eqref{S-eq:bh_threshold_report_identity} gives
		\(M_Lp(u_{\rm BH})=\alpha k\). Hence
		\[
		U_S(u_{\rm BH})
		=M_L\widehat p_L(u_{\rm BH})
		\le M_Lp(u_{\rm BH})+M_L\Delta_L
		=\alpha k+M_L\Delta_L.
		\]
		By \eqref{S-eq:bh_false_candidates_missed},
		\[
		V\le U_S(u_{\rm BH})+2^dD
		\le \alpha k+M_L\Delta_L+2^dD.
		\]
		Applying \eqref{S-eq:report_proportion_algebra} with \(A=\alpha\)
		and \(E=M_L\Delta_L\) gives
		\[
		\begin{aligned}
			\frac Vk
			&\le \alpha
			+(1-\alpha)\frac{M_L\Delta_L}{N+M_L\Delta_L}
			+2^d(1-\alpha)\frac DN\\
			&= \alpha
			+(1-\alpha)\frac{\Delta_L}{\beta_L+\Delta_L}
			+2^d(1-\alpha)\frac DN.
		\end{aligned}
		\]
		When \(k=0\), the false discovery proportion is zero, so the
		same upper bound holds for \(V/\max\{k,1\}\). Taking expectations
		and using \eqref{def:fdr} and \eqref{S-eq:bh_power_missed_identity}
		proves \eqref{S-eq:bh_fdr_finite_bound}.
		
		Now suppose that \(\Delta_L/\beta_L\to0\) in probability.
		Since
		\[
		0\le\frac{\Delta_L}{\beta_L+\Delta_L}
		\le\min\!\left\{1,\frac{\Delta_L}{\beta_L}\right\},
		\]
		for every \(\eta>0\) we have
		\[
		\E\!\left[\frac{\Delta_L}{\beta_L+\Delta_L}\right]
		\le \eta+\Prob\{\Delta_L/\beta_L>\eta\}.
		\]
		Taking the upper limit and then letting \(\eta\downarrow0\)
		gives \(\E[\Delta_L/(\beta_L+\Delta_L)]\to0\).
		Together with \(\Power(u_{\rm BH})\to1\),
		\eqref{S-eq:bh_fdr_finite_bound} yields
		\(\limsup_{L\to\infty}\FDR(u_{\rm BH})\le\alpha\).
	\end{proof}

	\begin{lemma}[Uniform convergence of the empirical noise tail]
		\label{S-lem:bh_uniform_noise_box_lln}
		Under Assumption~\ref{ass:noise_reg_main}, if \(z\) is ergodic
		under translations of \(\R^d\), then
		\begin{equation}\label{S-eq:bh_uniform_noise_box_L1}
			\E[\Delta_L]\longrightarrow0.
		\end{equation}
	\end{lemma}
	
	\begin{proof}
		For \(m\in\mathbb Z^d\), let \(Z_m:=\max_{t\in\overline Q_m}S^z(t)\).
		The stationary Gaussian field \(z\) is ergodic under translations
		of \(\R^d\), hence also under the grid shifts
		\(t\mapsto t+(r/2)m\), \(m\in\mathbb Z^d\)
		\cite{passeggeriVeraart2019,bergelsonRosenblatt1988}.
		Since the cubes in \eqref{eq:Qm} are grid translates of \(Q_0\),
		their maxima \((Z_m)_{m\in\mathbb Z^d}\) inherit stationarity
		and ergodicity.
		For each fixed \(u\in\R\), the ergodic theorem gives
		\[
		\widehat p_L(u)
		=\frac{1}{M_L}\sum_{m\in\mathcal I_L}\mathbf 1_{\{Z_m\ge u\}}
		\longrightarrow \Prob(Z_0\ge u)
		\qquad\text{almost surely}.
		\]
		By stationarity and Theorem~\ref{thm:continuity_test},
		\(\Prob(Z_0\ge u)=p(u)\) and \(p\) is continuous.
		Since \(\mathbb Q\) is countable, the convergence holds almost surely
		simultaneously for all \(u\in\mathbb Q\). Applying P\'olya's theorem
		to the distribution functions \(x\mapsto\widehat p_L(-x)\) gives
		\[
		\Delta_L=\sup_{u\in\R}|\widehat p_L(u)-p(u)|
		\longrightarrow0
		\qquad\text{almost surely}.
		\]
		Since \(0\le\Delta_L\le1\), bounded convergence gives
		\(\E[\Delta_L]\to0\), proving \eqref{S-eq:bh_uniform_noise_box_L1}.
	\end{proof}

	\begin{proof}[Proof of Proposition~\ref{prop:bh_fdr_control}]
		Suppose first that \(N\to\infty\) and
		Assumption~\ref{ass:noise_exceedance_variance} holds.
		Fix \(0<\epsilon<1\). Since \(\Power(u_{\rm BH})\to1\)
		and \(N\to\infty\), taking the upper limit in
		\eqref{S-eq:bh_fdr_relative_noise_bound} gives
		\[
		\limsup_{L\to\infty}\FDR(u_{\rm BH})
		\le \alpha(1+\epsilon).
		\]
		Letting \(\epsilon\downarrow0\) gives the FDR bound in
		Proposition~\ref{prop:bh_fdr_control} when \(N\to\infty\)
		and Assumption~\ref{ass:noise_exceedance_variance} holds.
		Now suppose that \(z\) is ergodic under translations of \(\R^d\)
		and that \(\liminf_{L\to\infty}\beta_L>0\), where
		\(\beta_L=N/M_L\) as in \eqref{S-eq:betaL}.
		By \eqref{S-eq:bh_uniform_noise_box_L1}, \(\E[\Delta_L]\to0\).
		Hence, for every \(\eta>0\), Markov's inequality gives
		\[
		\Prob\!\left\{\frac{\Delta_L}{\beta_L}>\eta\right\}
		\le \frac{\E[\Delta_L]}{\eta\beta_L}\longrightarrow0.
		\]
		Thus \(\Delta_L/\beta_L\to0\) in probability. Together with
		\(\Power(u_{\rm BH})\to1\), Lemma~\ref{S-lem:bh_fdr_general_bound}
		yields
		\[
		\limsup_{L\to\infty}\FDR(u_{\rm BH})\le\alpha.
		\]
	\end{proof}
	
	\subsubsection{The BH threshold}
	\label{S-subsubsec:bh_threshold_order}
	
	We prove the two claims in Remark~\ref{rem:threshold_orders}
	about the BH threshold when \(\Power(u_{\rm BH})\to1\).
	Whenever \(k\ge1\), \eqref{eq:u_bh_def} gives
	\(u_{\rm BH}=u_{\alpha k/M_L}\). We therefore first bound the
	number \(k\) of BH reports.
	
	\begin{lemma}[Bounds on the number of BH reports]
		\label{S-lem:bh_report_count_bounds}
		Fix \(0<\alpha<1\), and let \(k\) be the BH index in
		\eqref{eq:k_BH_def}. Under Assumptions~\ref{assumptions}
		and~\ref{ass:noise_reg_main}, for every \(L>0\),
		\begin{equation}\label{S-eq:bh_report_count_lower}
			\Prob\!\left(k<\frac N2\right)
			\le 2\bigl[1-\Power(u_{\rm BH})\bigr].
		\end{equation}
		If Assumption~\ref{ass:noise_exceedance_variance} also holds,
		then for every \(0<\epsilon<1\) with \(\alpha(1+\epsilon)<1\)
		and every \(H>2^d/[1-\alpha(1+\epsilon)]\),
		\begin{equation}\label{S-eq:bh_report_count_upper}
			\Prob(k\ge HN)
			\le \frac{27C_{\rm var}}{\alpha\epsilon^3HN},
		\end{equation}
		where \(C_{\rm var}\) is the constant in
		\eqref{eq:noise_exceedance_variance}.
	\end{lemma}
	
	\begin{proof}
		Write \(D=D(u_{\rm BH})\). By \eqref{S-eq:bh_threshold_report_identity},
		\(k=N-D+V(u_{\rm BH})\ge N-D\). Hence Markov's inequality and
		\eqref{S-eq:bh_power_missed_identity} give
		\[
		\Prob(k<N/2)
		\le \Prob(D>N/2)
		\le 2\E[D/N]
		= 2[1-\Power(u_{\rm BH})],
		\]
		proving \eqref{S-eq:bh_report_count_lower}.
		For \eqref{S-eq:bh_report_count_upper}, if \(HN>M_L\), then
		\(\Prob(k\ge HN)=0\) since \(k\le M_L\). Suppose \(HN\le M_L\)
		and set \(t_0:=\alpha HN/M_L\in(0,1)\). On \(\{k\ge HN\}\),
		we have \(k\ge1\) and \(\alpha k/M_L\ge t_0\). If the event
		\(\mathcal G_L(t_0,\epsilon)\) in \eqref{S-eq:relative_noise_tail_event}
		also holds, then \eqref{S-eq:bh_noise_at_random_threshold} gives
		\[
		U_S(u_{\rm BH})\le\alpha(1+\epsilon)k.
		\]
		Applying \eqref{S-eq:report_count_algebra} with
		\(A=\alpha(1+\epsilon)<1\) and \(E=0\), we obtain
		\[
		k\le\frac{N+(2^d-1)D}{1-\alpha(1+\epsilon)}
		\le\frac{2^dN}{1-\alpha(1+\epsilon)}
		<HN.
		\]
		Thus \(\{k\ge HN\}\subseteq\mathcal G_L(t_0,\epsilon)^c\).
		By \eqref{S-eq:bh_relative_tail_probability},
		\[
		\Prob(k\ge HN)
		\le\frac{27C_{\rm var}}{M_L\epsilon^3t_0}
		=\frac{27C_{\rm var}}{\alpha\epsilon^3HN}.
		\]
	\end{proof}
	
	\begin{lemma}[BH threshold when \(N/M_L\to0\)]
		\label{S-lem:bh_threshold_sparse}
		Fix \(0<\alpha<1\). Suppose that Assumptions~\ref{assumptions},
		\ref{ass:noise_reg_main}, and~\ref{ass:noise_exceedance_variance}
		hold. If \(N/M_L\to0\) and \(\Power(u_{\rm BH})\to1\), then
		for every \(0<\varepsilon<1\),
		\begin{equation}\label{S-eq:bh_threshold_sparse_order}
			\Prob\!\left\{
			2(1-\varepsilon)\sigma_{\mathcal V}^2\log(M_L/N)
			\le u_{\rm BH}\le
			2(1+\varepsilon)\sigma_{\mathcal V}^2\log(M_L/N)
			\right\}\longrightarrow1.
		\end{equation}
	\end{lemma}
	
	\begin{proof}
		Set \(h_L:=\log(M_L/N)\), which tends to infinity since
		\(N/M_L\to0\). Fix \(0<\zeta<1\) such that \(\alpha(1+\zeta)<1\).
		For all sufficiently large \(L\),
		\(h_L>2^d/[1-\alpha(1+\zeta)]\). Applying
		\eqref{S-eq:bh_report_count_lower} and \eqref{S-eq:bh_report_count_upper}
		with \(H=h_L\) and \(\epsilon=\zeta\) gives
		\[
		\Prob\bigl(\{N/2\le k<Nh_L\}^c\bigr)
		\le 2[1-\Power(u_{\rm BH})]
		+\frac{27C_{\rm var}}{\alpha\zeta^3Nh_L}
		\longrightarrow0,
		\]
		since \(\Power(u_{\rm BH})\to1\) and \(Nh_L\ge h_L\to\infty\).
		
		Define the probability levels
		\[
		c_L^-:=\frac{\alpha N}{2M_L}
		=\frac{\alpha}{2}e^{-h_L},
		\qquad
		c_L^+:=\frac{\alpha Nh_L}{M_L}
		=\alpha h_Le^{-h_L}.
		\]
		Both tend to zero, and \(0<c_L^-\le c_L^+<1\) for all
		sufficiently large \(L\). On \(\{N/2\le k<Nh_L\}\), we have
		\(k\ge1\) and \(c_L^-\le\alpha k/M_L\le c_L^+\).
		Since \(u_{\rm BH}=u_{\alpha k/M_L}\) by \eqref{eq:u_bh_def}
		and \(c\mapsto u_c\) is nonincreasing by \eqref{eq:uc_def_inline},
		\begin{equation}\label{S-eq:bh_threshold_sparse_bracket}
			\Prob\{u_{c_L^+}\le u_{\rm BH}\le u_{c_L^-}\}
			\ge \Prob\{N/2\le k<Nh_L\}\longrightarrow1.
		\end{equation}
		Since
		\[
		\log(1/c_L^-)=h_L+\log(2/\alpha),
		\qquad
		\log(1/c_L^+)=h_L-\log h_L+\log(1/\alpha),
		\]
		we have \(\log(1/c_L^\pm)/h_L\to1\). Applying
		\eqref{S-eq:noise_quantile_asymptotic} at \(c_L^-\) and \(c_L^+\) gives
		\[
		\frac{u_{c_L^\pm}}{2\sigma_{\mathcal V}^2h_L}
		=
		\frac{u_{c_L^\pm}}{2\sigma_{\mathcal V}^2\log(1/c_L^\pm)}
		\frac{\log(1/c_L^\pm)}{h_L}
		\longrightarrow1.
		\]
		Together with \eqref{S-eq:bh_threshold_sparse_bracket}, this proves
		\eqref{S-eq:bh_threshold_sparse_order}.
	\end{proof}
	
	\begin{lemma}[BH threshold when \(\liminf_{L\to\infty} N/M_L>0\)]
		\label{S-lem:bh_threshold_positive_density}
		Fix \(0<\alpha<1\). Suppose that Assumptions~\ref{assumptions}
		and~\ref{ass:noise_reg_main} hold, that \(\Power(u_{\rm BH})\to1\),
		and that \(\liminf_{L\to\infty}N/M_L>0\). Then, for every
		\(0<\beta_0<\liminf_{L\to\infty}N/M_L\),
		\begin{equation}\label{S-eq:bh_threshold_positive_density_bounds}
			\Prob\!\left\{
			u_\alpha\le u_{\rm BH}\le u_{\alpha\beta_0/2}
			\right\}\longrightarrow1,
		\end{equation}
		where \(0<u_\alpha\le u_{\alpha\beta_0/2}<\infty\) and both
		bounds are independent of \(L\).
	\end{lemma}
	
	\begin{proof}
		For all sufficiently large \(L\), \(N/M_L\ge\beta_0\).
		On \(\{k\ge N/2\}\), we have \(k\ge1\), and \(k\le M_L\) gives
		\[
		\frac{\alpha\beta_0}{2}
		\le\frac{\alpha N}{2M_L}
		\le\frac{\alpha k}{M_L}
		\le\alpha.
		\]
		By \eqref{eq:u_bh_def} and the monotonicity of \(u_c\) in
		\eqref{eq:uc_def_inline}, it follows that
		\[
		u_\alpha\le u_{\rm BH}\le u_{\alpha\beta_0/2}.
		\]
		Since \(N\le M_L\), we have \(0<\alpha\beta_0/2<\alpha<1\).
		Corollary~\ref{S-cor:generalized_inverse_attainment} therefore
		gives \(0<u_\alpha\le u_{\alpha\beta_0/2}<\infty\).
		Thus, for all sufficiently large \(L\),
		\[
		\{k\ge N/2\}
		\subseteq \{u_\alpha\le u_{\rm BH}\le u_{\alpha\beta_0/2}\}.
		\]
		By \eqref{S-eq:bh_report_count_lower} and \(\Power(u_{\rm BH})\to1\),
		\[
		\Prob(k\ge N/2)
		=1-\Prob(k<N/2)
		\ge 1-2[1-\Power(u_{\rm BH})]\longrightarrow1.
		\]
		The inclusion gives \eqref{S-eq:bh_threshold_positive_density_bounds}.
	\end{proof}
	
	\subsubsection{The linear subspace model}
	
	\begin{corollary}[The linear subspace model]
		\label{S-cor:exact_linear_subspace}
		Replace the model~\eqref{eq:template_instance_model} by
		\[
		x(t)=\sum_{i=1}^N x_i(t),
		\qquad
		x_i(t)=\sum_{j=1}^M a_{i,j}\psi_j(t-\tau_i),
		\]
		where the coefficient vectors
		\(a_i=(a_{i,1},\dots,a_{i,M})\in\mathbb R^M\setminus\{0\}\)
		are deterministic and may depend on \(L\).
		With these vectors in place of~\eqref{eq:x} in
		Assumption~\ref{def:localization_property}
		and~\eqref{eq:aL_def},
		Propositions~\ref{prop:asymptotic_error_power},
		\ref{prop:bh_density_strength}, and~\ref{prop:bh_fdr_control}
		and Corollaries~\ref{cor:simultaneous_recovery_asymptotic}
		and~\ref{cor:bh_power_asymptotic} hold for this model
		under their respective remaining assumptions, without
		assuming that \(\mathcal V\) is steerable. Moreover,
		\(\Gamma_{\mathcal V}\) may be chosen to depend only on \(\rho\).
        These conclusions also hold for \(d=1\).
	\end{corollary}
	
	\begin{proof}[Proof of Corollary~\ref{S-cor:exact_linear_subspace}]
		In the linear subspace model, Lemma~\ref{S-lem:representative_scores}
		gives \(\kappa_{\mathcal V}=1\). Its proof also verifies
		Lemma~\ref{S-lem:support_and_exclusion}\textup{(i)} for this model.
		The geometry in Lemma~\ref{S-lem:support_and_exclusion}\textup{(ii)-(iii)}
		uses only the detection and exclusion cubes and the separation assumptions.
		The proofs of Proposition~\ref{prop:asymptotic_error_power}
		and Proposition~\ref{prop:bh_density_strength} therefore apply
		with \(\kappa_{\mathcal V}=1\), as does
		Corollary~\ref{S-cor:bon_fdr_recovery}. In particular,
		\eqref{S-eq:bh_gamma_choice} gives
		\[
		\gamma_{\rm BH}=\min\{1/4,\rho^2/64\},
		\qquad \Gamma_{\mathcal V}=2/\gamma_{\rm BH}.
		\]
		Hence \(\Gamma_{\mathcal V}\) depends only on \(\rho\), and
		\eqref{S-eq:bh_power_explicit_constants} gives the claimed
		dependence of the BH power constants with the same absolute
		\(C_{\rm BH}\).
		The false-report bound \eqref{S-eq:bh_false_candidates_missed}
		follows from the same support and separation properties, so the
		proof of Proposition~\ref{prop:bh_fdr_control} applies as well.
		Applying the strength conditions of
		Corollaries~\ref{cor:simultaneous_recovery_asymptotic}
		and~\ref{cor:bh_power_asymptotic} to the respective finite
		bounds gives their asymptotic conclusions in the linear subspace model.
	None of these arguments uses \(d\ge2\): the dimension enters
only through \(M_L=n_L^d\) and the constants \(2^d\), \(3^d\),
and \(C_{\rm group}\). The proof therefore also applies for \(d=1\).
    \end{proof}
	
	\subsection{Sufficient conditions for the noise assumptions}
\label{S-app:noise_exceedance_variance}

We first verify the variance bound in
Assumption~\ref{ass:noise_exceedance_variance}
when the noise covariance vanishes beyond a fixed distance.
We then verify nondegeneracy and ergodicity for the
noise covariance used in the simulations.

Fix the detection subspace \(\mathcal V\) and its
orthonormal basis \(\psi_1,\ldots,\psi_M\) from
Subsection~\ref{subsec:algorithm_outline}.
For the noise field \(z\) in
Subsection~\ref{subsec:background_setting}, write
\(K(w):=\Cov(z(0),z(w))\).
For real-valued, compactly supported
\(f,g\in L^2(\R^d)\), let
\[
    F=\int_{\R^d}f(x)z(x)\,\dd x,
    \qquad
    G=\int_{\R^d}g(y)z(y)\,\dd y.
\]
Stationarity and Fubini's theorem give
\begin{equation}
    \label{S-eq:projected_noise_kernel_covariance}
    \Cov(F,G)
    =\int_{\R^d}\int_{\R^d}
    f(x)g(y)K(y-x)\,\dd y\,\dd x.
\end{equation}
Fubini's theorem applies because compact support and
square integrability imply \(f,g\in L^1(\R^d)\),
and \(\E|z(x)z(y)|\le K(0)\).

\begin{lemma}[Finite-range covariance]
    \label{S-lem:noise_exceedance_variance_finite_range}
    Suppose that Assumption~\ref{ass:noise_reg_main} holds
    and that, for some \(0\le R_0<\infty\),
    \[
        K(w)=0
        \qquad\text{whenever }\|w\|_\infty>R_0.
    \]
    Then Assumption~\ref{ass:noise_exceedance_variance}
    holds with a constant depending only on
    \(d,R_0,B\), and \(r\).
\end{lemma}

\begin{proof}
    Fix \(L>0\) and \(u\in\R\).
    For each grid cube \(Q_m\) in~\eqref{eq:Qm}, let
    \[
        J_m(u)
        :=\mathbf 1_{\{
        \max_{t\in\overline Q_m}S^z(t)\ge u\}}.
    \]
    For \(p(u)\) in~\eqref{eq:test_value_intro},
    stationarity and Theorem~\ref{thm:continuity_test}
    give \(\E[J_m(u)]=p(u)\).

    Recall from Assumption~\ref{ass:noise_reg_main} that
    \[
        X_j(t)
        =(z*\tilde\psi_j)(t)
        =\int_{\R^d}z(x)\psi_j(x-t)\,\dd x.
    \]
    Since \(\supp(\psi_j)\subseteq D_{B/2}
\subseteq\overline{\mathcal B}(0,B)\), the function
\(\psi_j(\cdot-t)\), for \(t\in\overline Q_m\),
is supported in
    \[
        A_m
        :=\overline Q_m+\overline{\mathcal B}(0,B)
        =\frac r2m+[-B/2,r/2+B/2]^d,
    \]
    by the cube definitions
    in~\eqref{eq:Qm} and~\eqref{eq:closed_cube}.

    If
    \[
        \|n-m\|_\infty>1+\frac{2(B+R_0)}r,
    \]
    then every \(x\in A_m\) and \(y\in A_n\) satisfy
    \[
        \|y-x\|_\infty
        \ge \frac r2\|n-m\|_\infty-(r/2+B)>R_0.
    \]
    Equation~\eqref{S-eq:projected_noise_kernel_covariance}
    therefore gives
    \[
        \Cov(X_j(t),X_k(s))=0
    \]
    for all \(t\in\overline Q_m\), \(s\in\overline Q_n\),
    and \(1\le j,k\le M\).
    The projected fields on the two cubes are jointly
Gaussian by Lemma~\ref{S-lem:projected_noise_continuity}
and uncorrelated, hence independent.
    Using the continuous versions from
    Lemma~\ref{S-lem:projected_noise_continuity}, each score
    maximum is determined by values on a countable dense
    subset of its cube. Thus \(J_m(u)\) and \(J_n(u)\)
    are independent.

    For every pair \(m,n\), Cauchy-Schwarz gives
    \[
        |\Cov(J_m(u),J_n(u))|
        \le
        \sqrt{\Var(J_m(u))\Var(J_n(u))}
        =p(u)(1-p(u))
        \le p(u).
    \]
    For each \(m\), at most
    \[
        C_{\rm var}
        :=\#\left\{h\in\mathbb Z^d:
        \|h\|_\infty\le1+\frac{2(B+R_0)}r\right\}
    \]
    indices \(n\) can have a nonzero covariance with
    \(J_m(u)\).
    Since \(U_S(u)=\sum_{m\in\mathcal I_L}J_m(u)\)
    by~\eqref{S-eq:noise_exceedance_count}, and
    \(|\mathcal I_L|=M_L\), we obtain
    \[
        \Var(U_S(u))
        =\sum_{m,n\in\mathcal I_L}
        \Cov(J_m(u),J_n(u))
        \le C_{\rm var}M_Lp(u).
    \]
    The constant \(C_{\rm var}\) is independent of \(L\)
    and \(u\), as required.
\end{proof}

The noise model constructed in
Section~\ref{S-app:intermediate_density_lower_bound}
satisfies Assumption~\ref{ass:noise_reg_main}
and the finite-range condition of
Lemma~\ref{S-lem:noise_exceedance_variance_finite_range}.

\begin{proposition}[Nondegeneracy and ergodicity for an isotropic exponential covariance]
    \label{S-prop:noise_assumptions_isotropic_exponential}
    Suppose that \(d=2\) and that the noise field \(z\) in
    Subsection~\ref{subsec:background_setting} has covariance
    \[
        K(w)=\sigma^2\exp(-\lambda\|w\|),
        \qquad w\in\R^2,\quad \sigma,\lambda>0.
    \]
    Then Assumption~\ref{ass:noise_reg_main} holds,
    and \(z\) is ergodic under translations of \(\R^2\).
    In particular, this applies to the covariance
    in~\eqref{eq:simulation_correlated_noise},
    with \(\lambda=3\log 10\).
\end{proposition}

\begin{proof}
    With the Fourier transform convention
    \[
        \widehat g(\omega)
        =\int_{\R^2}e^{-\mathrm i\omega\cdot x}
        g(x)\,\dd x,
    \]
    the covariance \(K\) satisfies
    \[
        \widehat K(\omega)
        =\frac{2\pi\sigma^2\lambda}
        {(\lambda^2+\|\omega\|^2)^{3/2}}>0,
        \qquad \omega\in\R^2.
    \]
    Fix \(v\in\R^M\setminus\{0\}\) and let
    \(f=\sum_{j=1}^M v_j\psi_j\).
    By orthonormality,
    \(\|f\|_{L^2(\R^2)}^2=\|v\|^2>0\).
    For the covariance matrix \(\Sigma_{\mathbf X}\)
    in Assumption~\ref{ass:noise_reg_main},
    equation~\eqref{S-eq:projected_noise_kernel_covariance}
    and Plancherel's identity give
\begin{equation}
    \label{S-eq:simulation_noise_nondegeneracy}
    \begin{aligned}
        v^\top\Sigma_{\mathbf X}v
        &=\Var\!\left(
        \int_{\R^2}f(x)z(x)\,\dd x\right)\\
        &=(2\pi)^{-2}\int_{\R^2}
        |\widehat f(\omega)|^2\widehat K(\omega)\,\dd\omega
        >0.
    \end{aligned}
\end{equation}
    The integral is positive because \(\widehat K>0\)
    at every frequency and
    \(\|\widehat f\|_{L^2(\R^2)}>0\).
    This proves Assumption~\ref{ass:noise_reg_main}.

    Finally, \(K(w)\to0\) as \(\|w\|\to\infty\).
    Since \(z\) is a stationary Gaussian field, it is
    mixing, and hence ergodic under translations of \(\R^2\)
    \cite[Theorem~2.4]{passeggeriVeraart2019}.
\end{proof}

	\section{Proof of the lower bound under FDR control}
	\label{S-app:intermediate_density_lower_bound}
	
	We first recall the setting of
	Proposition~\ref{prop:retained_energy_lower_bounds}.
	For \(L>B\), let \(\mathcal A_{L,N}\) be the set of configurations
	\(\boldsymbol\tau=(\tau_1,\dots,\tau_N)\) in
	\(\overline{\mathcal B}(0,L-B)\) satisfying the separation
	condition~\eqref{eq:separation}. A detection procedure observes
	\(y\) on \(\overline{\mathcal B}(0,L+B)\) and reports a finite set
	\(\widehat T_L\) of locations in \(\overline{\mathcal B}(0,L)\).
	Its output is a measurable function of the observation and may also
	depend on an auxiliary random variable independent of the observation,
	whose law is the same for every configuration. For fixed
	\(\boldsymbol\tau\in\mathcal A_{L,N}\), write
	\(\FDR_{\boldsymbol\tau}(\widehat T_L)\) and
	\(\Power_{\boldsymbol\tau}(\widehat T_L)\) for the FDR and
	power defined in \eqref{def:fdr} and \eqref{def:power},
	respectively, with \(\widehat T_L\) as the set of reported locations.
	We prove the following statement, which implies
Proposition~\ref{prop:retained_energy_lower_bounds}.
There exist a radial template \(u_1\)
	of unit \(L^2\)-norm, a centered stationary Gaussian noise field,
	and a constant \(c_{\rm low}>0\) with the following properties.
	The template, its one-dimensional steerable span
	\(\mathcal V_0=\spn\{u_1\}\), the noise law, and \(c_{\rm low}\)
	are fixed as \(L\) varies. For the basis \(\{u_1\}\),
	Assumptions~\ref{ass:noise_reg_main}
	and~\ref{ass:noise_exceedance_variance} hold. For every
	\(\boldsymbol\tau\in\mathcal A_{L,N}\) and every choice of positive
	amplitudes, the clean score satisfies the localization property
	of Assumption~\ref{def:localization_property} with \(\rho=1\).
	
	Fix \(0<\alpha<1\), let \(N=N(L)\ge1\) satisfy \(N/M_L\to0\), and
	let each instance be a translated copy of \(a_Lu_1\), where \(a_L>0\).
	Since \(u_1\) has unit norm and spans the detection subspace, each
	instance has squared projection norm \(a_L^2\). Suppose that the
	strength condition~\eqref{eq:intermediate_lower_strength} holds,
	\[
	\limsup_{L\to\infty}
	\frac{a_L^2}{\eta_{\mathrm{peak}}\log(M_L/N)}
	<c_{\rm low}.
	\]
	Consider any sequence of detection procedures \(\{\widehat T_L\}\)
	such that, for all sufficiently large \(L\) and every
	\(\boldsymbol\tau\in\mathcal A_{L,N}\), distinct reported locations
	are almost surely at \(\ell_\infty\)-distance at least \(r/2\).
	If these procedures control the FDR uniformly over the admissible
	configurations, in the sense that
	\begin{equation}\label{S-eq:intermediate_uniform_fdr}
		\limsup_{L\to\infty}
		\sup_{\boldsymbol\tau\in\mathcal A_{L,N}}
		\FDR_{\boldsymbol\tau}(\widehat T_L)
		\le\alpha,
	\end{equation}
	then their power satisfies
	\begin{equation}\label{S-eq:intermediate_lower_power}
		\limsup_{L\to\infty}
		\inf_{\boldsymbol\tau\in\mathcal A_{L,N}}
		\Power_{\boldsymbol\tau}(\widehat T_L)
		\le\alpha.
	\end{equation}
	
	\begin{definition}[Gaussian block model]
		\label{S-def:gaussian_block_model}
		Consider \(D\) locations indexed by \(1,\ldots,D\). Fix an integer
		\(1\le s\le D\) and set \(q=\lfloor D/s\rfloor\). Following
		\cite{abrahamCastilloRoquain2024}, we divide the first \(sq\) locations
		into \(s\) groups of size \(q\), called blocks. We choose one location
		uniformly from each block, independently across blocks, and denote
		the set of chosen locations by \(I\).
		For a fixed \(\mu\ge0\), we observe \(X=(X_1,\ldots,X_D)\), where
		\begin{equation}\label{S-eq:gaussian_block_model}
			X_j=\mu\mathbf 1_{\{j\in I\}}+\varepsilon_j,
			\qquad j=1,\ldots,D.
		\end{equation}
		The noise variables \(\varepsilon_1,\ldots,\varepsilon_D\) are
		independent standard normal variables and are independent of \(I\).
		Thus, given \(I\), the observations are independent, with mean \(\mu\)
		at the chosen locations and zero elsewhere.
	\end{definition}
	The following lemma bounds from below the expected sum of two
	proportions: the fraction of reports at unchosen locations and the
	fraction of chosen locations that the rule misses.
	
	\begin{lemma}[A Gaussian block bound]
		\label{S-lem:gaussian_product_block_lower_bound}
		Under the model of Definition~\ref{S-def:gaussian_block_model},
		let \(A=(A_1,\ldots,A_D)\) be the output of a reporting rule,
		with \(A_j=1\) if the rule reports location \(j\) and \(A_j=0\)
		otherwise. Assume that the rule uses only \(X\), in the sense that
		for every \(a\in\{0,1\}^D\),
		\[
		\Prob(A=a\mid X,I)=\Prob(A=a\mid X)
		\qquad\text{almost surely}.
		\]
		Let \(R\) be the number of reports and \(T\) the number of
		chosen locations that the rule reports,
		\[
		R:=\sum_{j=1}^D A_j,
		\qquad
		T:=\sum_{j\in I}A_j.
		\]
		Since one location is chosen in each of the \(s\) blocks,
		\(|I|=s\), and \((s-T)/s\) is the fraction of chosen locations
		that the rule misses. Then
		\begin{equation}\label{S-eq:gaussian_product_block_lower_bound}
			\E\!\left[\frac{R-T}{\max\{R,1\}}+\frac{s-T}{s}\right]
			\ge 1-\sqrt{\frac{e^{\mu^2}}{q}}.
		\end{equation}
	\end{lemma}
	
	\begin{proof}
		For each location \(j\), let
		\[
		\pi_j(X):=\Prob(j\in I\mid X).
		\]
		By the assumption on the reporting rule, \(A\) and \(I\) are
		conditionally independent given \(X\). Thus, once \(X\) is known,
		knowing the reports leaves the probability that \(j\) was chosen
		unchanged:
		\[
		\Prob(j\in I\mid X,A)=\pi_j(X).
		\]
		
		We first compute the expected number of correct reports given
		\(X\) and \(A\). A report at \(j\) is correct precisely when
		\(A_j=1\) and \(j\in I\), so
		\(T=\sum_{j=1}^D A_j\mathbf 1_{\{j\in I\}}\). By linearity of
		conditional expectation,
		\[
		\begin{aligned}
			m:=\E[T\mid X,A]
			&=\sum_{j=1}^D
			\E[A_j\mathbf 1_{\{j\in I\}}\mid X,A] \\
			&=\sum_{j=1}^D A_j\,
			\E[\mathbf 1_{\{j\in I\}}\mid X,A] \\
			&=\sum_{j=1}^D A_j\,\Prob(j\in I\mid X,A)
			=\sum_{j=1}^D A_j\pi_j(X).
		\end{aligned}
		\]
		In the second line, each \(A_j\) is known because we condition
		on \(A\). In the third, we use the fact that the conditional
		expectation of an indicator is the conditional probability of
		its event.
		
		Since \(T\) counts the locations that are both reported and chosen,
		it is at most the number \(R\) of reports and at most the number \(s\)
		of chosen locations, hence \(0\le T\le\min\{R,s\}\).
		By monotonicity of conditional expectation,
		\[
		0\le m=\E[T\mid X,A]
		\le \E[\min\{R,s\}\mid X,A]
		=\min\{R,s\}.
		\]
		The last equality holds because \(s\) is fixed and
		\(R=\sum_{j=1}^D A_j\) is determined by \(A\).
		On the event \(R>0\), we may therefore take \(1/R\) and \(1/s\)
		outside the conditional expectation. By linearity,
		\begin{equation}\label{S-eq:block_conditional_error}
			\begin{aligned}
				\E\!\left[\left.
				\frac{R-T}{R}+\frac{s-T}{s}
				\right|X,A\right]
				&= 2-\frac{1}{R}\E[T\mid X,A]
				-\frac{1}{s}\E[T\mid X,A] \\
				&= \left(1-\frac{m}{R}\right)
				+\left(1-\frac{m}{s}\right).
			\end{aligned}
		\end{equation}
		
		To obtain a lower bound for \eqref{S-eq:block_conditional_error}
		that depends only on \(X\), recall that \(m\) is the sum of
		\(\pi_j(X)\) over the \(R\) reported locations.
		Cauchy-Schwarz inequality gives
		\begin{equation}\label{S-eq:block_reported_sum_bound}
			\frac{m^2}{Rs}
			= \frac{1}{Rs}\left(\sum_{j:A_j=1}\pi_j(X)\right)^2
			\le \frac{1}{s}\sum_{j:A_j=1}\pi_j(X)^2
			\le \frac{1}{s}\sum_{j=1}^D\pi_j(X)^2.
		\end{equation}
		Since \(1-m/s\ge0\) and \(0\le m/R\le1\), multiplying the second
		term in \eqref{S-eq:block_conditional_error} by \(m/R\) gives a
		lower bound. Applying \eqref{S-eq:block_reported_sum_bound} in the
		last step, we obtain
		\begin{equation}\label{S-eq:block_positive_report_bound}
			\begin{aligned}
				\left(1-\frac{m}{R}\right)+\left(1-\frac{m}{s}\right)
				&\ge \left(1-\frac{m}{R}\right)
				+\frac{m}{R}\left(1-\frac{m}{s}\right) \\
				&= 1-\frac{m^2}{Rs}
				\ge 1-\frac{1}{s}\sum_{j=1}^D\pi_j(X)^2.
			\end{aligned}
		\end{equation}
		If \(R=0\), then \(T=0\). In this case,
		\((R-T)/\max\{R,1\}=0\) and \((s-T)/s=1\).
		The sum of the two proportions is therefore one. Since
		\(\sum_{j=1}^D\pi_j(X)^2\ge0\), we have
		\[
		1\ge 1-\frac{1}{s}\sum_{j=1}^D\pi_j(X)^2.
		\]
		Together with \eqref{S-eq:block_positive_report_bound} for \(R>0\),
		this gives
		\[
		\E\!\left[\left.
		\frac{R-T}{\max\{R,1\}}+\frac{s-T}{s}
		\right|X,A\right]
		\ge 1-\frac{1}{s}\sum_{j=1}^D\pi_j(X)^2
		\qquad\text{almost surely}.
		\]
		Taking expectations with respect to the joint distribution of \(X\)
		and \(A\), the law of total expectation gives the mean over \(I,X,A\)
		on the left. On the right, each \(\pi_j(X)\) depends only on \(X\).
		Thus
		\begin{equation}\label{S-eq:product_block_posterior_risk}
			\E_{I,X,A}\!\left[
			\frac{R-T}{\max\{R,1\}}+\frac{s-T}{s}
			\right]
			\ge 1-\frac{1}{s}\E_X\!\left[
			\sum_{j=1}^D\pi_j(X)^2
			\right].
		\end{equation}
		By Definition~\ref{S-def:gaussian_block_model}, locations outside
		the \(s\) blocks cannot be chosen, so \(\pi_j(X)=0\) at each such
		location. We can therefore bound
		\(\E_X[\sum_{j=1}^D\pi_j(X)^2]\) in
		\eqref{S-eq:product_block_posterior_risk} by summing the bounds for
		the individual blocks.
		Now fix a block \(\mathcal J\). We first derive an upper bound for
		\[
		\E_X\!\left[\sum_{j\in\mathcal J}\pi_j(X)^2\right].
		\]
		Let \(Q_0\) be the law of \(q\) independent standard normal
		observations. If \(j\) is the chosen location in \(\mathcal J\),
		then \(X_j\) has law \(N(\mu,1)\), and the other observations in
		the block are independent standard normal variables, independent
		of \(X_j\). Thus the joint density of the block observations differs
		from the density under \(Q_0\) only in the factor for \(X_j\).
		Their likelihood ratio is
		\[
		L_j
		:=\frac{\exp(-(X_j-\mu)^2/2)}{\exp(-X_j^2/2)}
		=e^{\mu X_j-\mu^2/2}.
		\]
		
		Write \(X_{\mathcal J}:=(X_k)_{k\in\mathcal J}\). Let \(f\) be the
		density of \(X_{\mathcal J}\), let \(f_k\) be its conditional density
		given \(k\in I\), and let \(f_0\) be its density under \(Q_0\).
		Exactly one location in \(\mathcal J\) is chosen, and each location
		has probability \(1/q\) of being chosen. The law of total probability
		therefore gives
		\[
		f(x)=\frac{1}{q}\sum_{k\in\mathcal J}f_k(x).
		\]
		Since \(f_k(X_{\mathcal J})/f_0(X_{\mathcal J})=L_k\), we have
		\begin{equation}\label{S-eq:block_density_ratio}
			\frac{f(X_{\mathcal J})}{f_0(X_{\mathcal J})}
			=\frac{1}{q}\sum_{k\in\mathcal J}
			\frac{f_k(X_{\mathcal J})}{f_0(X_{\mathcal J})}
			=\frac{S}{q},
			\qquad S:=\sum_{k\in\mathcal J}L_k.
		\end{equation}
		The observations outside \(\mathcal J\) are independent of the chosen
		location in \(\mathcal J\) and of \(X_{\mathcal J}\), jointly. Hence
		\(\pi_j(X)=\Prob(j\in I\mid X_{\mathcal J})\) for \(j\in\mathcal J\).
		By Bayes' formula and \eqref{S-eq:block_density_ratio},
		\begin{equation}\label{S-eq:product_block_posterior}
			\pi_j(X)
			=\frac{(1/q)f_j(X_{\mathcal J})}{f(X_{\mathcal J})}
			=\frac{(1/q)L_j}{S/q}
			=\frac{L_j}{S}.
		\end{equation}
		For any nonnegative measurable function \(g\colon\R^q\to[0,\infty)\),
		\eqref{S-eq:block_density_ratio} gives
		\begin{equation}\label{S-eq:block_change_of_measure}
			\begin{aligned}
				\E_X[g(X_{\mathcal J})]
				&= \int_{\R^q} g(x)f(x)\,\dd x \\
				&= \int_{\R^q} g(x)\frac{f(x)}{f_0(x)}f_0(x)\,\dd x \\
				&= \E_{Q_0}\!\left[g(X_{\mathcal J})\frac{S}{q}\right].
			\end{aligned}
		\end{equation}
		By \eqref{S-eq:product_block_posterior}, the sum
		\(\sum_{j\in\mathcal J}\pi_j(X)^2\) is a function of
		\(X_{\mathcal J}\). Applying \eqref{S-eq:block_change_of_measure}
		with \(g(X_{\mathcal J})=\sum_{j\in\mathcal J}(L_j/S)^2\) gives
		\[
		\begin{aligned}
			\E_X\!\left[\sum_{j\in\mathcal J}\pi_j(X)^2\right]
			&= \E_{Q_0}\!\left[
			\frac{S}{q}\sum_{j\in\mathcal J}
			\left(\frac{L_j}{S}\right)^2
			\right] \\
			&= \frac{1}{q}\E_{Q_0}\!\left[
			\frac{\sum_{j\in\mathcal J}L_j^2}{S}
			\right].
		\end{aligned}
		\]
		Since each \(L_j>0\),
		\[
		S^2
		=\sum_{j\in\mathcal J}L_j^2
		+2\sum_{\substack{j,k\in\mathcal J\\j<k}}L_jL_k
		\ge \sum_{j\in\mathcal J}L_j^2.
		\]
		Since \(S\ge\sqrt{\sum_{j\in\mathcal J}L_j^2}>0\), we have
		\begin{equation}\label{S-eq:block_likelihood_ratio_bound}
			\frac{\sum_{j\in\mathcal J}L_j^2}{S}
			\le \sqrt{\sum_{j\in\mathcal J}L_j^2}.
		\end{equation}
		Applying \eqref{S-eq:block_likelihood_ratio_bound} inside the expectation
		and then using Jensen's inequality for the concave square root gives
		\[
		\begin{aligned}
			\E_X\!\left[\sum_{j\in\mathcal J}\pi_j(X)^2\right]
			&= \frac{1}{q}\E_{Q_0}\!\left[
			\frac{\sum_{j\in\mathcal J}L_j^2}{S}
			\right] \\
			&\le \frac{1}{q}\E_{Q_0}\!\left[
			\sqrt{\sum_{j\in\mathcal J}L_j^2}
			\right] \\
			&\le \frac{1}{q}\sqrt{\sum_{j\in\mathcal J}\E_{Q_0}[L_j^2]}.
		\end{aligned}
		\]
		Under \(Q_0\), \(X_j\) is standard normal. Taking \(t=2\mu\) in
		\(\E_{Q_0}[e^{tX_j}]=e^{t^2/2}\) gives
		\[
		\E_{Q_0}[L_j^2]
		=e^{-\mu^2}\E_{Q_0}[e^{2\mu X_j}]
		=e^{-\mu^2}e^{2\mu^2}
		=e^{\mu^2}.
		\]
		Since \(\mathcal J\) contains \(q\) locations, we obtain
		\begin{equation}\label{S-eq:single_block_posterior_bound}
			\E_X\!\left[\sum_{j\in\mathcal J}\pi_j(X)^2\right]
			\le \frac{1}{q}\sqrt{q e^{\mu^2}}
			= \sqrt{\frac{e^{\mu^2}}{q}}.
		\end{equation}
		The block \(\mathcal J\) was arbitrary, so
		\eqref{S-eq:single_block_posterior_bound} holds for each of the
		\(s\) blocks \(\mathcal J_1,\ldots,\mathcal J_s\).
		Since \(\pi_j(X)=0\) outside these blocks, linearity of expectation gives
		\[
		\frac{1}{s}\E_X\!\left[\sum_{j=1}^D\pi_j(X)^2\right]
		= \frac{1}{s}\sum_{\ell=1}^s
		\E_X\!\left[\sum_{j\in\mathcal J_\ell}\pi_j(X)^2\right]
		\le \frac{1}{s}\sum_{\ell=1}^s\sqrt{\frac{e^{\mu^2}}{q}}
		= \sqrt{\frac{e^{\mu^2}}{q}}.
		\]
		Substituting this upper bound for
		\(s^{-1}\E_X[\sum_{j=1}^D\pi_j(X)^2]\) into
		\eqref{S-eq:product_block_posterior_risk} yields
		\eqref{S-eq:gaussian_product_block_lower_bound}.
	\end{proof}
	
	Lemma~\ref{S-lem:gaussian_product_block_lower_bound} gives a lower
	bound for a reporting rule that observes a finite vector of Gaussian
	measurements, with one location chosen in each block. A detection
	procedure in Proposition~\ref{prop:retained_energy_lower_bounds}
	observes the full image. We construct a model in which the image
	values on a grid have the distribution required by the lemma.
	We then show that, given these values, the conditional law of the
	full image is the same for every configuration of instance locations
	on the grid.
	
	Recall that \(r=2B+\delta\), where \(0<\delta\le 2B\).
	Choose a nonzero, nonnegative, radial function
	\(\phi\in C_c^\infty(\R^d)\) supported in \(D_{\delta/32}\),
	the closed Euclidean ball of radius \(\delta/32\) centered at the origin.
	We choose this radius so that after the two convolutions used below,
	the support lies in \(D_{\delta/8}\), inside the detection cube
	\(\mathcal B(0,\delta)\).
	For a function \(f\), write \(\widetilde f(t):=f(-t)\). Define
	\begin{equation}\label{S-eq:lower_bound_template_definition}
		c_\phi:=\|\phi*\widetilde\phi\|_{L^2},
		\qquad
		u_1:=\frac{\phi*\widetilde\phi}{c_\phi}.
	\end{equation}
	The convolution is continuous and
	\[
	(\phi*\widetilde\phi)(0)
	=\int_{\R^d}\phi(v)^2\,\dd v
	=\|\phi\|_{L^2}^2>0,
	\]
	so \(c_\phi>0\). By \eqref{S-eq:lower_bound_template_definition},
	\[
	u_1(t)=\frac{1}{c_\phi}\int_{\R^d}\phi(v)\phi(v-t)\,\dd v.
	\]
	The integrand is nonnegative. For it to be nonzero, both
	\(\|v\|\) and \(\|v-t\|\) must be at most \(\delta/32\).
	The triangle inequality then gives
	\[
	\|t\|\le \|v\|+\|v-t\|\le\delta/16.
	\]
	Thus \(\supp(u_1)\subset D_{\delta/16}\subset D_{B/2}\), where
	the second inclusion uses \(\delta\le 2B\).
	The template \(u_1\) is smooth because \(\phi\) is smooth and
	compactly supported. It is radial because rotating both \(t\) and
	the integration variable \(v\) leaves the integral unchanged.
	By construction, \(\|u_1\|_{L^2}=1\). Since every rotation leaves
	\(u_1\) unchanged, \(\mathcal V_0:=\spn\{u_1\}\) is a one-dimensional
	steerable subspace with orthonormal basis \(\{u_1\}\).
	
	To construct a noise field whose covariance is proportional to \(u_1\),
	let \(\{\mathcal W(f):f\in L^2(\R^d)\}\) be a centered, jointly Gaussian
	family of real random variables, linear in \(f\) as a map into
	\(L^2(\Omega)\), with
	\[
	\E[\mathcal W(f)\mathcal W(g)]
	=\int_{\R^d}f(v)g(v)\,\dd v.
	\]
	For each \(t\in\R^d\), let \(\phi_t(v):=\phi(v-t)\) be the
	translate of \(\phi\) centered at \(t\). Define the noise at \(t\) by
	\begin{equation}\label{S-eq:lower_bound_noise_definition}
		z_0(t):=\mathcal W(\phi_t).
	\end{equation}
	By \eqref{S-eq:lower_bound_noise_definition} and the covariance
	identity for \(\mathcal W\),
	\begin{equation}\label{S-eq:lower_bound_covariance_template}
		\begin{aligned}
			\Cov(z_0(t),z_0(s))
			&=\langle\phi_t,\phi_s\rangle_{L^2}\\
			&=\int_{\R^d}\phi(v-t)\phi(v-s)\,\dd v\\
			&=\int_{\R^d}\phi(w)\phi\bigl(w-(t-s)\bigr)\,\dd w\\
			&=(\phi*\widetilde\phi)(t-s)
			=c_\phi u_1(t-s).
		\end{aligned}
	\end{equation}
	The third line uses \(w=v-s\), and the last line follows from
	\(\widetilde\phi(v)=\phi(-v)\) and
	\eqref{S-eq:lower_bound_template_definition}.
	Since \(z_0\) is a centered Gaussian field and its covariance
	in \eqref{S-eq:lower_bound_covariance_template} depends only on
	\(t-s\), \(z_0\) is stationary.
	
	Each instance uses the radial template \(u_1\). Since rotations
	preserve distance from the origin,
	\[
	(R_i u_1)(v)=u_1(R_i^{-1}v)=u_1(v),
	\qquad v\in\R^d.
	\]
	Thus, for locations
	\(\boldsymbol\tau=(\tau_1,\ldots,\tau_N)\in\mathcal A_{L,N}\)
	and positive amplitudes \(b_1,\ldots,b_N\), the clean image in
	\eqref{eq:template_instance_model} is
	\[
	x(t)=\sum_{i=1}^N b_i u_1(t-\tau_i).
	\]
	With \(z_0\) as the noise field in \eqref{eq:the_model}, write
	\(Y_L\) for the observation \(y\) on \(\overline{\mathcal B}(0,L+B)\).
	Then
	\begin{equation}\label{S-eq:lower_bound_spatial_model}
		Y_L(t)=x(t)+z_0(t),
		\qquad t\in\overline{\mathcal B}(0,L+B).
	\end{equation}
	Note that the template \(u_1\), its span \(\mathcal V_0\), and the law of
	\(z_0\) are fixed as \(L\) varies.
	
	\begin{proof}[Proof of Proposition~\ref{prop:retained_energy_lower_bounds}]
		We verify the noise assumptions and the localization property for
		the image model in \eqref{S-eq:lower_bound_spatial_model}, beginning
		with continuity of the noise field \(z_0\).
		
		Recall that \(\phi_t(v)=\phi(v-t)\) and
		\(z_0(t)=\mathcal W(\phi_t)\). By linearity of \(\mathcal W\),
		\[
		z_0(t)-z_0(s)=\mathcal W(\phi_t-\phi_s).
		\]
		Thus the increment is centered Gaussian with variance
		\(\|\phi_t-\phi_s\|_{L^2}^2\). To bound this variance, integrate
		the derivative of \(\phi\) along the segment from \(v-s\) to \(v-t\):
		\[
		\phi_t(v)-\phi_s(v)
		=-\int_0^1
		\nabla\phi\bigl(v-s-\theta(t-s)\bigr)\cdot(t-s)\,\dd\theta.
		\]
		Taking the \(L^2\)-norm in \(v\) and using the triangle inequality
		for the integral gives
		\[
		\|\phi_t-\phi_s\|_{L^2}
		\le \|t-s\|\int_0^1
		\|\nabla\phi(\,\cdot-s-\theta(t-s))\|_{L^2}\,\dd\theta
		=\|t-s\|\,\|\nabla\phi\|_{L^2},
		\]
		where the equality follows from translation invariance of the
		\(L^2\)-norm.
		Take \(p=2d\) and let \(C_p=\E|G|^p\) for \(G\sim N(0,1)\).
		Since \(z_0(t)-z_0(s)\) has the same distribution as
		\(\|\phi_t-\phi_s\|_{L^2}G\), we obtain
		\[
		\E|z_0(t)-z_0(s)|^p
		=C_p\|\phi_t-\phi_s\|_{L^2}^p
		\le C_p\|\nabla\phi\|_{L^2}^p\|t-s\|^p.
		\]
		Since \(p=2d>d\), the Kolmogorov-Chentsov theorem
		\cite[Theorem~4.23]{kallenberg2021} gives a modification
		of \(z_0\) with continuous sample paths. We use this jointly
		measurable modification throughout the proof.
		
		For the basis \(\{u_1\}\), the projected noise has a single
		coordinate,
		\[
		X(t):=(z_0*\widetilde u_1)(t).
		\]
		Assumption~\ref{ass:noise_reg_main} therefore requires
		\(\Var(X(0))>0\). To verify this, observe that
		\(X(0)=\int_{\R^d}z_0(v)u_1(v)\,\dd v\). By stationarity and
		the Cauchy-Schwarz inequality,
		\[
		\E|z_0(v)z_0(0)|
		\le \Var(z_0(0))=\|\phi\|_{L^2}^2.
		\]
		Since \(u_1\in L^1(\R^d)\), Fubini's theorem and
		\eqref{S-eq:lower_bound_covariance_template} give
		\begin{equation}\label{S-eq:lower_bound_projected_noise}
			\begin{aligned}
				\Cov(X(0),z_0(0))
				&=\int_{\R^d}u_1(v)\Cov(z_0(v),z_0(0))\,\dd v\\
				&=c_\phi\int_{\R^d}u_1(v)^2\,\dd v
				=c_\phi.
			\end{aligned}
		\end{equation}
		By \eqref{S-eq:lower_bound_projected_noise} and Cauchy-Schwarz,
		\[
		0<c_\phi^2
		=\Cov(X(0),z_0(0))^2
		\le \Var(X(0))\Var(z_0(0)).
		\]
		Thus \(\Var(X(0))>0\), as required by
		Assumption~\ref{ass:noise_reg_main}.
		The strength condition \eqref{eq:intermediate_lower_strength}
		uses the expected maximum noise score \(\eta_{\mathrm{peak}}\).
		For the basis \(\{u_1\}\), the noise score is \(X(t)^2\).
		Since Assumption~\ref{ass:noise_reg_main} holds,
		Lemma~\ref{S-lem:noise_maxima_fixed_growing} gives
		\begin{equation}\label{S-eq:lower_bound_eta_fixed}
			0<\eta_{\mathrm{peak}}
			=\E\!\left[
			\max_{t\in\overline{\mathcal B}(0,r/2)}X(t)^2
			\right]<\infty.
		\end{equation}
		Both the law of \(X\) and the comparison cube
		\(\overline{\mathcal B}(0,r/2)\) are fixed as \(L\) varies,
		so \(\eta_{\mathrm{peak}}\) is independent of \(L\).
		
		To verify Assumption~\ref{ass:noise_exceedance_variance}, we use
		Lemma~\ref{S-lem:noise_exceedance_variance_finite_range}. Its remaining
		condition is that the covariance of \(z_0\) vanish beyond a fixed
		distance. If \(\|h\|_\infty>\delta/16\), then
		\(\|h\|\ge\|h\|_\infty>\delta/16\). Since \(u_1\) is supported
		in \(D_{\delta/16}\), it follows that \(u_1(h)=0\). Thus
		\eqref{S-eq:lower_bound_covariance_template} gives
		\[
		\Cov(z_0(h),z_0(0))=c_\phi u_1(h)=0.
		\]
		The lemma therefore applies with \(R_0=\delta/16\).

		We next verify the localization property for the model
		in~\eqref{S-eq:lower_bound_spatial_model}.
		Consider any \(\boldsymbol\tau\in\mathcal A_{L,N}\) and any positive
		amplitudes, ordered so that \(b_1\ge\cdots\ge b_N\).
		The coefficient of the \(i\)th centered instance in the basis
		\(\{u_1\}\) is \(\langle b_i u_1,u_1\rangle=b_i\), so
		\(a_i=(b_i)\). Let \(C_{u_1}:=u_1*\widetilde u_1\).
		Convolving the clean image in \eqref{S-eq:lower_bound_spatial_model}
		with \(\widetilde u_1\) gives
		\[
		S^x(t)
		=\left(\sum_{i=1}^N b_i C_{u_1}(t-\tau_i)\right)^2.
		\]
		Note that \(C_{u_1}(0)=\|u_1\|_{L^2}^2=1\).
		Since \(u_1\) is supported in a ball of radius \(\delta/16\),
		its support and its translate by \(t\) are disjoint whenever
		\(\|t\|>\delta/8\). Hence
		\[
		\supp(C_{u_1})\subset D_{\delta/8}\subset\mathcal B(0,\delta).
		\]
		By the separation condition~\eqref{eq:separation},
		at \(t=\tau_i\) only the \(i\)th term can be nonzero.
		Outside the union of the detection cubes every term vanishes.
		Consequently,
		\[
		S^x(\tau_i)=b_i^2,
		\qquad
		S^x(t)=0
		\quad\text{for }t\notin
		\bigcup_{j=1}^N\mathcal B(\tau_j,\delta).
		\]
		
		Take \(q_i=\tau_i\) in \eqref{eq:qi}. For \(j<i\), the comparison
		set in \eqref{eq:comparison set} excludes
		\(\mathcal B(\tau_j,2B)\), which contains
		\(\mathcal B(\tau_j,\delta)\) because \(\delta\le2B\).
		For \(j\ge i\), it excludes \(\mathcal B(\tau_j,\delta)\) directly.
		Thus every \(t\) in the comparison set satisfies
		\[
		S^x(q_i)-S^x(t)=b_i^2=\|a_i\|^2.
		\]
		This proves the localization property of
		Assumption~\ref{def:localization_property} with \(\rho=1\)
		for every \(\boldsymbol\tau\in\mathcal A_{L,N}\).
		
		Having verified the localization property for every valid configuration,
		we now restrict the instance locations to a separated grid and set all
		amplitudes equal to \(a_L\). Since \(\|u_1\|_{L^2}=1\), each instance
		then has squared projection norm \(a_L^2\).
		
		Set \(h=2B+2\delta\). For all sufficiently large \(L\), define
		\begin{equation}\label{S-eq:lower_bound_grid_definition}
			\Lambda_L:=
			\left\{hk:k\in\mathbb Z^d,\quad
			\|k\|_\infty\le
			\left\lfloor\frac{L-B}{2h}\right\rfloor\right\}.
		\end{equation}
		Every grid point lies in \(\overline{\mathcal B}(0,L-B)\), and
		distinct grid points are at \(\ell_\infty\)-distance at least
		\(h>B+3\delta/2\). Thus any \(N\) distinct grid points satisfy the
		location constraints defining \(\mathcal A_{L,N}\).
		Counting the grid points in each coordinate gives
		\[
		|\Lambda_L|
		=\left(2\left\lfloor\frac{L-B}{2h}\right\rfloor+1\right)^d,
		\qquad
		\frac{|\Lambda_L|}{L^d}\longrightarrow h^{-d}.
		\]
		Since \(M_L/L^d\to(2/r)^d\) by \eqref{S-eq:grid_count_asymptotic},
		\begin{equation}\label{S-eq:lower_bound_log_packing}
			\frac{|\Lambda_L|}{M_L}
			\longrightarrow \left(\frac{r}{2h}\right)^d>0.
		\end{equation}
		In particular, \(N/M_L\to0\) implies \(N\le|\Lambda_L|\) for all
		sufficiently large \(L\).
		
		Fix a set \(I\subset\Lambda_L\) with \(|I|=N\), and let
		\(\boldsymbol\tau_I\) be the vector of its points in any fixed order.
		With an instance of amplitude \(a_L\) at each point of \(I\),
		the observation model \eqref{S-eq:lower_bound_spatial_model} becomes
		\[
		Y_L(t)=a_L\sum_{\tau\in I}u_1(t-\tau)+z_0(t),
		\qquad t\in\overline{\mathcal B}(0,L+B).
		\]
		Let \(\sigma_\phi:=\|\phi\|_{L^2}\) and define
		\begin{equation}\label{S-eq:lower_bound_observed_coordinates}
			\xi_\tau:=\frac{Y_L(\tau)}{\sigma_\phi},
			\quad \tau\in\Lambda_L,
			\qquad
			\lambda_\phi:=\frac{\sigma_\phi}{c_\phi},
			\qquad
			\mu_L:=\lambda_\phi a_L.
		\end{equation}
		Fix \(\tau\in\Lambda_L\). For every other grid point \(\tau'\),
		the grid spacing exceeds the support radius of \(u_1\), so
		\(u_1(\tau-\tau')=0\). Thus \(x(\tau)=a_Lu_1(0)\) if
		\(\tau\in I\), and \(x(\tau)=0\) otherwise. By
		\eqref{S-eq:lower_bound_template_definition},
		\[
		u_1(0)
		=\frac{(\phi*\widetilde\phi)(0)}{c_\phi}
		=\frac{\|\phi\|_{L^2}^2}{c_\phi}
		=\frac{\sigma_\phi^2}{c_\phi}.
		\]
		For \(\tau\in I\), it follows that
		\(x(\tau)/\sigma_\phi=a_L\sigma_\phi/c_\phi=\mu_L\). Hence
		\[
		\xi_\tau
		=\mu_L\mathbf 1_{\{\tau\in I\}}
		+\frac{z_0(\tau)}{\sigma_\phi},
		\qquad \tau\in\Lambda_L.
		\]
		By \eqref{S-eq:lower_bound_noise_definition},
		\(z_0(\tau)=\mathcal W(\phi_\tau)\), where \(\mathcal W\) is a centered,
		jointly Gaussian family. The grid noise values are therefore centered
		and jointly Gaussian. By \eqref{S-eq:lower_bound_covariance_template},
		\[
		\Cov(z_0(\tau),z_0(\tau'))
		=c_\phi u_1(\tau-\tau')
		=
		\begin{cases}
			\sigma_\phi^2, & \tau=\tau',\\
			0,             & \tau\ne\tau'.
		\end{cases}
		\]
		Consequently, the normalized noise values
		\((z_0(\tau)/\sigma_\phi)_{\tau\in\Lambda_L}\) are independent
		standard normal variables. For fixed \(I\), the coordinates of
		\(\boldsymbol\xi:=(\xi_\tau)_{\tau\in\Lambda_L}\) are thus
		independent, with
		\begin{equation}\label{S-eq:intermediate_gaussian_sequence}
			\xi_\tau\sim N\!\left(\mu_L\mathbf 1_{\{\tau\in I\}},1\right),
			\qquad \tau\in\Lambda_L.
		\end{equation}	
		
		So far, we have considered only the image values on the grid.
		The detection procedure takes the full image \(Y_L\) as input and
		reports locations in \(\overline{\mathcal B}(0,L)\).
		To apply the lower bound on the grid, we proceed in two steps.
		First, we show that the conditional law of \(Y_L\) given
		\(\boldsymbol\xi\) is the same for every set \(I\) of instance
		locations. A randomized procedure using only \(\boldsymbol\xi\)
		can therefore reproduce the distribution of the spatial reports
		for every \(I\). We then convert the spatial reports into decisions
		at the grid points while preserving power and without increasing FDR.
		For fixed \(I\), the difference \(Y_L(t)-z_0(t)\) is deterministic
		for every \(t\). Since \(\xi_\tau=Y_L(\tau)/\sigma_\phi\),
		\eqref{S-eq:lower_bound_covariance_template} gives
		\begin{equation}\label{S-eq:lower_bound_image_grid_covariance}
			\Cov(Y_L(t),\xi_\tau)
			=\frac{\Cov(z_0(t),z_0(\tau))}{\sigma_\phi}
			=\frac{u_1(t-\tau)}{\lambda_\phi}.
		\end{equation}
		To determine the conditional law of \(Y_L\) given \(\boldsymbol\xi\),
		we write \(Y_L(t)\) as a weighted sum of the grid values and a
		residual. Given \(\boldsymbol\xi\), the sum is known. It therefore
		suffices to show that the residual is independent of \(\boldsymbol\xi\)
		and has the same law for every \(I\). Using the covariances in
		\eqref{S-eq:lower_bound_image_grid_covariance} as weights, define
		\begin{equation}\label{S-eq:lower_bound_residual_definition}
			\mathcal R_L(t)
			:=Y_L(t)-\frac{1}{\lambda_\phi}
			\sum_{\tau\in\Lambda_L}u_1(t-\tau)\xi_\tau.
		\end{equation}
		Substituting the expressions for \(Y_L(t)\) and \(\xi_\tau\) into
		\eqref{S-eq:lower_bound_residual_definition} gives
		\[
		\mathcal R_L(t)
		=\left(a_L-\frac{\mu_L}{\lambda_\phi}\right)
		\sum_{\tau\in I}u_1(t-\tau)
		+z_0(t)-\frac{1}{\lambda_\phi}
		\sum_{\tau\in\Lambda_L}u_1(t-\tau)
		\frac{z_0(\tau)}{\sigma_\phi}.
		\]
		Since \(\mu_L=\lambda_\phi a_L\), the first term vanishes.
		The remaining terms depend only on \(z_0\), so the law of
		\(\mathcal R_L\) is the same for every \(I\).
		It remains to show that, for fixed \(I\), the residual function
		\(\mathcal R_L\) is independent of the grid vector \(\boldsymbol\xi\).
		By \eqref{S-eq:intermediate_gaussian_sequence}, the coordinates of
		\(\boldsymbol\xi\) are independent and have unit variance. Hence 
		\begin{equation}\label{S-eq:lower_bound_residual_grid_covariance}
			\begin{aligned}
				\Cov(\mathcal R_L(t),\xi_\tau)
				&=\Cov(Y_L(t),\xi_\tau)
				-\frac{1}{\lambda_\phi}\sum_{\tau'\in\Lambda_L}
				u_1(t-\tau')\Cov(\xi_{\tau'},\xi_\tau)\\
				&=\Cov(Y_L(t),\xi_\tau)
				-\frac{u_1(t-\tau)}{\lambda_\phi}
				=0,
			\end{aligned}
		\end{equation}
		where the last equality follows from
		\eqref{S-eq:lower_bound_image_grid_covariance}.
		Let \(K=\overline{\mathcal B}(0,L+B)\), and fix a countable dense
		subset \(D\subset K\). For any \(t_1,\ldots,t_m\in D\), the vectors
		\((\mathcal R_L(t_1),\ldots,\mathcal R_L(t_m))\) and
		\(\boldsymbol\xi\) are jointly Gaussian, since their entries are
		linear combinations of values of \(z_0\) plus deterministic terms.
		By \eqref{S-eq:lower_bound_residual_grid_covariance}, all
		cross-covariances vanish, so the two vectors are independent.
		The residual has continuous paths, and the Borel sigma-field on \(C(K)\)
		is generated by evaluations on \(D\). Since every finite vector of
		residual values on \(D\) is independent of \(\boldsymbol\xi\), it follows
		that \(\mathcal R_L\) and \(\boldsymbol\xi\) are independent.
		The detection procedure
		uses only \(Y_L\) and any independent randomization with a fixed law.
		Thus, for given grid values, the conditional distribution of its
		reported locations is also the same for every \(I\).

		\par\medskip\noindent
		The spatial procedure returns a set of reported locations
		\(\widehat T_L\). To apply Lemma~\ref{S-lem:gaussian_product_block_lower_bound},
		we convert \(\widehat T_L\) into a vector \(\boldsymbol\varphi\) with
		one binary decision at each grid point. A point \(\tau\in\Lambda_L\) is selected if
\(\mathcal B(\tau,\delta)\) contains a reported location:
		\[
		\varphi_\tau
		:=\mathbf 1_{\{\widehat T_L\cap\mathcal B(\tau,\delta)\ne\emptyset\}},
		\qquad
		\boldsymbol\varphi:=(\varphi_\tau)_{\tau\in\Lambda_L}.
		\]
		We compare the false discovery rate and power of the selected grid
		points with those of the original spatial reports.
		Recall that \(I\subset\Lambda_L\) is the set of \(N\) instance
		locations and that
		\(\boldsymbol\xi=(Y_L(\tau)/\sigma_\phi)_{\tau\in\Lambda_L}\)
		contains the normalized image values on the grid.
		By \eqref{S-eq:lower_bound_residual_definition}, the conditional
		law of \(\widehat T_L\) given \(\boldsymbol\xi\) is the same for
		every \(I\), since the residual is independent of \(\boldsymbol\xi\)
		and has the same law for every \(I\). The conversion from
		\(\widehat T_L\) to \(\boldsymbol\varphi\) uses only the fixed grid.
		Hence the conditional law of \(\boldsymbol\varphi\) given
		\(\boldsymbol\xi\) is also the same for every \(I\).
		To compare the counts, note that the cubes
		\(\mathcal B(\tau,\delta)\), \(\tau\in\Lambda_L\), are disjoint
		because their centers are at \(\ell_\infty\)-distance at least
		\(h>\delta\). Each cube contains at most one spatial report almost
		surely. Indeed, two points in the same open cube of side length
		\(\delta\) are at \(\ell_\infty\)-distance less than \(\delta\),
		and distinct reports are separated by at least
		\(r/2=B+\delta/2\ge\delta\).
		For fixed \(I\), let \(C_I\) count the spatial reports in cubes
		centered at points of \(I\), \(F_I\) the spatial reports in cubes
		centered at the remaining grid points, and \(G_I\) the spatial
		reports outside \(\bigcup_{\tau\in\Lambda_L}\mathcal B(\tau,\delta)\). Thus
		\[
		C_I:=\sum_{\tau\in I}\varphi_\tau,
		\qquad
		F_I:=\sum_{\tau\in\Lambda_L\setminus I}\varphi_\tau,
		\qquad
		G_I:=\#\left(\widehat T_L\setminus
		\bigcup_{\tau\in\Lambda_L}\mathcal B(\tau,\delta)\right).
		\]
		The set \(\widehat T_L\) contains \(C_I\) correct reports and
		\(F_I+G_I\) false reports. The vector \(\boldsymbol\varphi\) selects
		\(C_I\) points in \(I\) and \(F_I\) points outside \(I\).
		Thus the conversion preserves every correct detection and discards
		only the \(G_I\) false reports outside all grid cubes. For fixed \(C_I\), the false discovery proportion
		\(x/\max\{C_I+x,1\}\) is nondecreasing in the number \(x\) of false
		reports. Taking \(x=F_I\) for \(\boldsymbol\varphi\) and \(x=F_I+G_I\)
		for \(\widehat T_L\) gives
		\begin{equation}\label{S-eq:intermediate_fdp_comparison}
			\frac{F_I}{\max\{C_I+F_I,1\}}
			\le
			\frac{F_I+G_I}{\max\{C_I+F_I+G_I,1\}}.
		\end{equation}
		For fixed \(I\), let \(\Prob_{I,L}\) and \(\E_{I,L}\) denote
		probability and expectation over the noise and any independent
		randomization used by the procedure. Taking expectations in
		\eqref{S-eq:intermediate_fdp_comparison} and using \eqref{def:fdr}
		gives
		\begin{equation}\label{S-eq:intermediate_fdr_comparison}
			\E_{I,L}\!\left[\frac{F_I}{\max\{C_I+F_I,1\}}\right]
			\le \FDR_{\boldsymbol\tau_I}(\widehat T_L).
		\end{equation}
		For power, the fraction of instance locations selected by
		\(\boldsymbol\varphi\) is \(C_I/N\). By the definition of
		\(\varphi_\tau\) and \eqref{def:power},
		\[
		\frac{\E_{I,L}[C_I]}{N}
		=\frac{1}{N}\sum_{\tau\in I}
		\Prob_{I,L}\{\widehat T_L\cap\mathcal B(\tau,\delta)\ne\emptyset\}
		=\Power_{\boldsymbol\tau_I}(\widehat T_L).
		\]
		Consequently,
		\begin{equation}\label{S-eq:intermediate_fnr_power}
			\E_{I,L}\!\left[\frac{N-C_I}{N}\right]
			=1-\Power_{\boldsymbol\tau_I}(\widehat T_L).
		\end{equation}
		
		\par\medskip\noindent
		The comparisons in \eqref{S-eq:intermediate_fdr_comparison}
		and \eqref{S-eq:intermediate_fnr_power} hold for each fixed set \(I\)
		of instance locations. To bound the infimum of the power over
		\(\mathcal A_{L,N}\), we now average over sets \(I\) chosen according
		to Definition~\ref{S-def:gaussian_block_model}.
		Fix an ordering of \(\Lambda_L\) and let
		\[
		q=\left\lfloor\frac{|\Lambda_L|}{N}\right\rfloor.
		\]
		Partition the first \(Nq\) grid points into \(N\) blocks of size \(q\).
		Choose one point uniformly from each block, independently across
		blocks and independently of the noise and any randomization used
		by the detection procedure. Let \(\pi_L\) denote the uniform distribution
		on the \(q^N\) sets obtained by choosing one grid point from each block.
		By the separation and boundary properties of
		the grid in \eqref{S-eq:lower_bound_grid_definition}, every resulting
		configuration \(\boldsymbol\tau_I\) belongs to \(\mathcal A_{L,N}\).
		Conditional on \(I\), equation~\eqref{S-eq:intermediate_gaussian_sequence}
		gives the Gaussian block model \eqref{S-eq:gaussian_block_model} with
		\(|\Lambda_L|\) coordinates, \(N\) blocks, and mean \(\mu_L\) at the
		chosen locations. The conditional law of \(\boldsymbol\varphi\) given
		\(\boldsymbol\xi\) is the same for every fixed \(I\), as shown using
		\eqref{S-eq:lower_bound_residual_definition}. Thus, under
		\(I\sim\pi_L\), for every \(a\in\{0,1\}^{|\Lambda_L|}\),
		\[
		\Prob(\boldsymbol\varphi=a\mid\boldsymbol\xi,I)
		=\Prob(\boldsymbol\varphi=a\mid\boldsymbol\xi)
		\qquad\text{almost surely}.
		\]
		We therefore apply Lemma~\ref{S-lem:gaussian_product_block_lower_bound}
		to the observations \(\boldsymbol\xi\) and the reports
		\(\boldsymbol\varphi\). In the notation of the lemma, \(R=C_I+F_I\)
		locations are reported, \(T=C_I\) of them belong to \(I\), and
		\(s=N\). Substituting into \eqref{S-eq:gaussian_product_block_lower_bound}
		and writing the expectation by conditioning on \(I\) gives
		\begin{equation}\label{S-eq:intermediate_acr_lower_bound}
			\E_{I\sim\pi_L}\!\left[
			\E_{I,L}\!\left[
			\frac{F_I}{\max\{C_I+F_I,1\}}+\frac{N-C_I}{N}
			\right]\right]
			\ge 1-\sqrt{\frac{e^{\mu_L^2}}{q}}.
		\end{equation}
		By \eqref{S-eq:intermediate_fdr_comparison}, the inner expectation
		of the first fraction is at most the spatial FDR. By
		\eqref{S-eq:intermediate_fnr_power}, the inner expectation of the
		second is one minus the spatial power. Applying both comparisons
		to \eqref{S-eq:intermediate_acr_lower_bound} and using
		\(\boldsymbol\tau_I\in\mathcal A_{L,N}\) gives
		\begin{equation}\label{S-eq:lower_bound_multiple_finite}
			\begin{aligned}
				\inf_{\boldsymbol\tau\in\mathcal A_{L,N}}
				\Power_{\boldsymbol\tau}(\widehat T_L)
				&\le \E_{I\sim\pi_L}
				\Power_{\boldsymbol\tau_I}(\widehat T_L)\\
				&\le \E_{I\sim\pi_L}
				\FDR_{\boldsymbol\tau_I}(\widehat T_L)
				+\sqrt{\frac{e^{\mu_L^2}}{q}}\\
				&\le \sup_{\boldsymbol\tau\in\mathcal A_{L,N}}
				\FDR_{\boldsymbol\tau}(\widehat T_L)
				+\sqrt{\frac{e^{\mu_L^2}}{q}}.
			\end{aligned}
		\end{equation}
		It remains to bound the last term in
		\eqref{S-eq:lower_bound_multiple_finite}. By
		\eqref{S-eq:lower_bound_log_packing},
		\(|\Lambda_L|/M_L\to(r/(2h))^d>0\). Since \(N/M_L\to0\),
		we have \(|\Lambda_L|/N\to\infty\), so
		\(q=\lfloor|\Lambda_L|/N\rfloor\sim|\Lambda_L|/N\). Consequently,
		\[
		\frac{q}{M_L/N}\longrightarrow\left(\frac{r}{2h}\right)^d,
		\qquad
		\log q=\log(M_L/N)+O(1).
		\]
		As \(M_L/N\to\infty\), it follows that
		\begin{equation}\label{S-eq:intermediate_log_equivalence}
			\frac{\log q}{\log(M_L/N)}\longrightarrow 1.
		\end{equation}
		Choose
		\[
		c_{\rm low}:=\frac{1}{\lambda_\phi^2\eta_{\mathrm{peak}}}
		=\frac{c_\phi^2}{\sigma_\phi^2\eta_{\mathrm{peak}}}>0.
		\]
		Since \(\mu_L=\lambda_\phi a_L\) by
		\eqref{S-eq:lower_bound_observed_coordinates}, the strength condition
		\eqref{eq:intermediate_lower_strength} and
		\eqref{S-eq:intermediate_log_equivalence} give
		\[
		\limsup_{L\to\infty}\frac{\mu_L^2}{\log q}
		=
		\lambda_\phi^2\eta_{\mathrm{peak}}\,
		\limsup_{L\to\infty}
		\frac{a_L^2}{\eta_{\mathrm{peak}}\log(M_L/N)}
		<\lambda_\phi^2\eta_{\mathrm{peak}}c_{\rm low}=1.
		\]
		Hence \(\mu_L^2\le(1-\varepsilon)\log q\) for some
		\(\varepsilon>0\) and all sufficiently large \(L\). Since \(q\to\infty\),
		\[
		\sqrt{\frac{e^{\mu_L^2}}{q}}
		\le q^{-\varepsilon/2}\longrightarrow0.
		\]
		Taking upper limits in \eqref{S-eq:lower_bound_multiple_finite} and
		applying the uniform FDR bound \eqref{S-eq:intermediate_uniform_fdr}
		proves \eqref{S-eq:intermediate_lower_power}.
	\end{proof}
	
	\section{Steerable subspaces and angular decomposition}
	\label{S-app:rotation_invariant}
	Recall that rotations act on functions by \((Rf)(t)=f(R^{-1}t)\)
	and that a detection subspace \(\mathcal V\) is steerable if
	\(R\mathcal V=\mathcal V\) for every \(R\in\SO(d)\), as in
	\eqref{def:rotation_invariant}. The following lemma shows that
	rotating a template preserves the norm of its projection onto
	\(\mathcal V\).
	
	\begin{lemma}[Projection commutes with rotations]
		\label{S-lem:proj_rot_appendix}
		Let \(\mathcal V\subseteq L^2(D_{B/2})\) be finite-dimensional and
		steerable. Then, for every \(R\in\SO(d)\) and
		\(f\in L^2(D_{B/2})\),
		\[
		P_{\mathcal V}(R f)=R P_{\mathcal V}f,
		\qquad
		\|P_{\mathcal V}(R f)\|_{L^2(D_{B/2})}
		=\|P_{\mathcal V}f\|_{L^2(D_{B/2})}.
		\]
	\end{lemma}
	
	\begin{proof}
		Rotations preserve \(D_{B/2}\) and Lebesgue measure, so
		\(R\) is unitary. Write \(P=P_{\mathcal V}\).
		Since \(\mathcal V\) is steerable, \(R(Pf)\in\mathcal V\).
		To identify it as the projection of \(Rf\), fix \(v\in\mathcal V\).
		Then \(R^{-1}v\in\mathcal V\), and \(f-Pf\) is orthogonal to
		\(\mathcal V\). Hence
		\[
		\langle Rf-R(Pf),v\rangle
		=\langle f-Pf,R^{-1}v\rangle=0.
		\]
		Thus \(P(Rf)=R(Pf)\). Taking norms gives the second identity,
		since rotations preserve the \(L^2\) norm.
	\end{proof}
	
	To construct steerable subspaces, we write each
	\(v\in D_{B/2}\setminus\{0\}\) as \(v=\varrho\omega\), where
	\(\varrho=\|v\|\) and \(\omega=v/\|v\|\in S^{d-1}\).
	For \(R\in\SO(d)\), we have
	\((Rf)(\varrho\omega)=f(\varrho R^{-1}\omega)\), so a rotation
	acts only on the angular variable. In two dimensions, rotations preserve the span of the Fourier
		basis functions at each angular frequency. The following lemma gives the angular expansion in spherical harmonics
	and shows that rotations preserve the space of harmonics of each degree.
	For each integer \(j\ge0\), let \(h_{d,j}\) be the dimension of the
	space of real spherical harmonics of degree \(j\) on \(S^{d-1}\).
	Choose an orthonormal basis \(Y_{j,1},\ldots,Y_{j,h_{d,j}}\) of
	that space with respect to surface area measure \(d\sigma\).
	For a real-valued \(f\in L^2(D_{B/2})\), write
	\(f(\varrho,\omega)=f(\varrho\omega)\). For almost every
	\(\varrho\in(0,B/2)\), define the angular coefficients of \(f\) by
	\begin{equation}\label{S-eq:angular_coefficients}
		c_{j,m}(\varrho)
		:=\int_{S^{d-1}}f(\varrho,\omega)Y_{j,m}(\omega)\,d\sigma(\omega),
		\qquad j\ge0,\quad 1\le m\le h_{d,j}.
	\end{equation}
For \(a\in L^2([0,B/2],\varrho^{d-1}\,d\varrho)\), write
\[
\|a\|_{L^2([0,B/2],\varrho^{d-1}\,d\varrho)}^2
:=\int_0^{B/2}|a(\varrho)|^2\varrho^{d-1}\,d\varrho.
\]
	
	\begin{lemma}[Real angular expansion and rotations]
		\label{S-lem:real_angular_expansion}
		Every real-valued \(f\in L^2(D_{B/2})\) has the expansion
		\[
		f(\varrho,\omega)
		=\sum_{j=0}^{\infty}\sum_{m=1}^{h_{d,j}}
		c_{j,m}(\varrho)Y_{j,m}(\omega)
		\qquad\text{in }L^2(D_{B/2}),
		\]
		with coefficients given by \eqref{S-eq:angular_coefficients}.
		For each \(j\), the inner sum is the degree-\(j\) component of \(f\).
		Components of distinct degrees are orthogonal in \(L^2(D_{B/2})\),
		and
		\begin{equation}\label{S-eq:angular_parseval}
			\|f\|_{L^2(D_{B/2})}^2
			=\sum_{j=0}^{\infty}\sum_{m=1}^{h_{d,j}}\|c_{j,m}\|_{L^2([0,B/2],\varrho^{d-1}\,d\varrho)}^2.
		\end{equation}
		Every rotation maps the space of spherical harmonics of degree \(j\)
		onto itself. Consequently, each \(\mathcal S_j\) defined in
		Subsection~\ref{subsec:construction_setting} is steerable, as is
		each subspace \(\mathcal S(I)\) defined in
			Subsection~\ref{subsec:candidate_score_maps}.
	\end{lemma}
	
	\begin{proof}
		Integration in spherical coordinates gives
		\[
		\int_0^{B/2}\|f(\varrho,\cdot)\|_{L^2(S^{d-1},d\sigma)}^2
		\,\varrho^{d-1}\,d\varrho
		=\|f\|_{L^2(D_{B/2})}^2<\infty.
		\]
		Thus \(f(\varrho,\cdot)\in L^2(S^{d-1},d\sigma)\) for almost
		every \(\varrho\in(0,B/2)\). At each such radius,
		\eqref{S-eq:angular_coefficients} gives its spherical harmonic
		coefficients. Applying Parseval's identity on \(S^{d-1}\) to
		\(f(\varrho,\cdot)\) for almost every \(\varrho\), we obtain
		\[
		\begin{aligned}
			\|f\|_{L^2(D_{B/2})}^2
			&=\int_0^{B/2}
			\sum_{j=0}^{\infty}\sum_{m=1}^{h_{d,j}}
			|c_{j,m}(\varrho)|^2\,\varrho^{d-1}\,d\varrho\\
			&=
			\sum_{j=0}^{\infty}\sum_{m=1}^{h_{d,j}}
			\int_0^{B/2}|c_{j,m}(\varrho)|^2\,\varrho^{d-1}\,d\varrho\\
			&=\sum_{j=0}^{\infty}\sum_{m=1}^{h_{d,j}}\|c_{j,m}\|_{L^2([0,B/2],\varrho^{d-1}\,d\varrho)}^2.
		\end{aligned}
		\]
		This proves \eqref{S-eq:angular_parseval}.
		Let \(f_K\) be the partial sum through degree \(K\). By
		orthonormality on the sphere, the angular coefficients of
		\(f-f_K\) vanish for \(j\le K\) and equal \(c_{j,m}\) for
		\(j>K\). Applying \eqref{S-eq:angular_parseval} to \(f-f_K\)
		therefore gives
		\[
		\|f-f_K\|_{L^2(D_{B/2})}^2
		=\sum_{j>K}\sum_{m=1}^{h_{d,j}}\|c_{j,m}\|_{L^2([0,B/2],\varrho^{d-1}\,d\varrho)}^2
		\longrightarrow0.
		\]
		Components of distinct degrees have zero angular inner product
		at almost every radius. Integrating this inner product against
		\(\varrho^{d-1}\,d\varrho\) gives their orthogonality in
		\(L^2(D_{B/2})\).
		Fix \(R\in\SO(d)\). Rotations preserve the space of spherical
		harmonics of each degree \cite[Theorem~1.1.7]{daiXu2013spheres}.
		By the definition of \(\mathcal S_j\) in
		Subsection~\ref{subsec:construction_setting}, this space is
		spanned by the functions \(g(\varrho,\omega)
		=a(\varrho)Y_{j,m}(\omega)\) with \(a\in V_j\) and
		\(1\le m\le h_{d,j}\). For each such \(g\),
		\[
		(Rg)(\varrho,\omega)=a(\varrho)Y_{j,m}(R^{-1}\omega).
		\]
		Since \(Y_{j,m}\circ R^{-1}\) is a linear combination of
		\(Y_{j,1},\ldots,Y_{j,h_{d,j}}\), we have \(Rg\in\mathcal S_j\).
		Hence \(R\mathcal S_j\subseteq\mathcal S_j\). Applying this
		inclusion to \(R^{-1}\) gives \(R\mathcal S_j=\mathcal S_j\).
		Every finite direct sum
		\(\mathcal S(I)\) is therefore steerable.
	\end{proof}
	
	We now use the angular expansion to construct a single steerable
	subspace that approximates every rotation of every template.
	
	\begin{proof}[Proof of Proposition~\ref{prop:ambient_uniform_coverage}]
		Fix \(\varepsilon>0\). For each template \(u_\ell\), let
		\(c_{\ell,j,m}\) be its angular coefficients in
		\eqref{S-eq:angular_coefficients}, and let \(u_{\ell,J}\) be its
		partial sum through degree \(J\). By
		Lemma~\ref{S-lem:real_angular_expansion},
		\(u_{\ell,J}\to u_\ell\) in \(L^2(D_{B/2})\) as \(J\to\infty\).
		Since there are finitely many templates, we can choose a single
		\(J\) such that
		\[
		\|u_\ell-u_{\ell,J}\|_{L^2(D_{B/2})}\le\varepsilon
		\qquad\text{for every }1\le\ell\le n.
		\]
		For \(0\le j\le J\), define
		\[
		V_j:=\spn\{c_{\ell,j,m}:1\le\ell\le n,\ 1\le m\le h_{d,j}\},
		\]
		and let \(\mathcal S=\mathcal S(J;V_0,\ldots,V_J)\) as in
		\eqref{eq:ambient_space}. 
		By Lemma~\ref{S-lem:real_angular_expansion}, \(\mathcal S\) is steerable and
		Lemma~\ref{S-lem:proj_rot_appendix} together with invariance of the \(L^2\) norm
		under rotations therefore give, for every \(R\in\SO(d)\),
		\[
		\|Ru_\ell-P_{\mathcal S}(Ru_\ell)\|_{L^2(D_{B/2})}
		=\|u_\ell-P_{\mathcal S}u_\ell\|_{L^2(D_{B/2})}.
		\]
		The projection \(P_{\mathcal S}u_\ell\) is the element of
		\(\mathcal S\) closest to \(u_\ell\) in the \(L^2\) norm. Since
		\(u_{\ell,J}\in\mathcal S\),
		\[
		\|u_\ell-P_{\mathcal S}u_\ell\|_{L^2(D_{B/2})}
		\le \|u_\ell-u_{\ell,J}\|_{L^2(D_{B/2})}
		\le \varepsilon.
		\]
		The bound holds for every \(\ell=1,\ldots,n\) and every \(R\in\SO(d)\),
		which proves \eqref{eq:uniform_template_approximation}.
	\end{proof}
	
	\section{Proofs for candidate detection subspaces and pSNR selection}
	\label{S-app:psnr_selection_proofs}
	
	Each candidate subspace \(\mathcal V\in\mathfrak S\) in
	\eqref{eq:candidate_family} is a direct sum of angular components
	\(\mathcal S_j\). Let \(S_{\mathcal V}^x\) and
		\(S_{\mathcal V}^z\) denote the score maps in~\eqref{eq:Sx}
		and~\eqref{eq:noise_score_map}, computed using an orthonormal
		basis of \(\mathcal V\).
	To prove Proposition~\ref{prop:template_localization_condition},
	we need to bound the clean score at a point \(t\) outside all
	detection cubes. The ball \(t+D_{B/2}\) may still intersect the
	supports of several instances. We first bound the contribution of
	one instance in terms of the distance from its center to \(t\).
	Fix \(\mathcal V\in\mathfrak S\) such that
	\(P_{\mathcal V}u_\ell\ne0\) for every \(\ell\). Define
	\begin{equation}\label{S-eq:radial_projection_envelope}
		e_{\mathcal V}(\lambda)
		:=\max_{1\le\ell\le n}\sup_{\|h\|=\lambda}
		\frac{\|P_{\mathcal V}[u_\ell(h+\cdot)]\|_{L^2(D_{B/2})}}
		{\|P_{\mathcal V}u_\ell\|_{L^2(D_{B/2})}},
		\qquad 0\le\lambda\le B.
	\end{equation}
	The same bound applies to a rotated template \(Ru_\ell\). Indeed,
	\[
	(Ru_\ell)(h+\cdot)=R\bigl[u_\ell(R^{-1}h+\cdot)\bigr].
	\]
	By Lemma~\ref{S-lem:proj_rot_appendix}, the ratio of projection norms
	for \(Ru_\ell\) with shift \(h\) equals the ratio for \(u_\ell\)
	with shift \(R^{-1}h\). Since \(\|R^{-1}h\|=\|h\|\), the ratio
	is bounded by \(e_{\mathcal V}(\|h\|)\).
	Suppose that the ball \(t+D_{B/2}\) intersects the supports of
	instances centered at \(\tau_{i_1},\ldots,\tau_{i_k}\). For each
	\(m\), the shift in \eqref{S-eq:radial_projection_envelope} is
	\(v_m:=t-\tau_{i_m}\), the vector from the instance center to \(t\). For \(1\le k\le 2^d\), let \(\mathcal C_k(B,\delta)\)
	be the set of tuples \((v_1,\ldots,v_k)\in(\mathbb R^d)^k\)
	satisfying
	\begin{equation}\label{S-eq:admissible_displacement_tuples}
		\begin{aligned}
			\|v_m\|\le B,\qquad \|v_m\|_\infty\ge\delta/2
			&\qquad (1\le m\le k),\\
			\|v_m-v_q\|_\infty\ge B+3\delta/2
			&\qquad (1\le m<q\le k).
		\end{aligned}
	\end{equation}
	The first bound follows because the ball and the ball supporting
	each instance both have radius \(B/2\) and must intersect. The second
	follows because \(t\) lies outside each detection cube. The pairwise
	bound follows from \eqref{eq:separation}, since
	\(v_m-v_q=\tau_{i_q}-\tau_{i_m}\). In particular, there are at most
	\(2^d\) such vectors: they lie in \([-B,B]^d\), which can be covered
	by \(2^d\) cubes of side length \(B\), each containing at most one
	vector.
	
	To bound the sum of the projection norm ratios for instances that
	contribute to the same ball, define
	\begin{equation}\label{S-eq:chi_V_def}
		\chi_{\mathcal V}(\delta)
		:=\max_{1\le k\le 2^d}
		\sup_{(v_1,\ldots,v_k)\in\mathcal C_k(B,\delta)}
		\sum_{m=1}^k e_{\mathcal V}(\|v_m\|),
	\end{equation}
	with the supremum over an empty set taken to be zero.
	We now use \eqref{S-eq:chi_V_def} to compare the clean score at each
	instance center with its values on the comparison set in
	\eqref{eq:comparison set}.
	
	\begin{proof}[Proof of Proposition~\ref{prop:template_localization_condition}]
		Fix \(0<\delta\le 2B\) and a candidate subspace \(\mathcal V\)
		satisfying the assumptions of the proposition. Suppose that
		\(\chi_{\mathcal V}(\delta)<1\). We verify
		Assumption~\ref{def:localization_property} with \(q_i=\tau_i\)
		for every \(i\).
		For \(s\in D_{B/2}\), the instance model
		\eqref{eq:template_instance_model} gives
		\[
		x_{j,t}(s):=x_j(t+s)=b_j(R_j u_{\ell_j})(t-\tau_j+s).
		\]
		By orthonormality and \eqref{eq:x},
		\(\|a_j\|=b_j\|P_{\mathcal V}(R_j u_{\ell_j})\|_{L^2(D_{B/2})}\).
		Applying the bound in \eqref{S-eq:radial_projection_envelope}
		to the rotated template \(R_j u_{\ell_j}\) with shift \(t-\tau_j\)
		therefore gives
		\begin{equation}\label{S-eq:shifted_instance_projection_bound}
			\|P_{\mathcal V}x_{j,t}\|_{L^2(D_{B/2})}
			\le e_{\mathcal V}(\|t-\tau_j\|)\|a_j\|,
			\qquad \|t-\tau_j\|\le B.
		\end{equation}
		Relabel the instances so that
		\(\|a_1\|\ge\cdots\ge\|a_N\|\), as in
		Assumption~\ref{def:localization_property}. Fix \(i\) and set
		\(q_i:=\tau_i\). Since \(\tau_i\in\overline{\mathcal B}(0,L)\)
		and \(\delta>0\), the point \(q_i\) satisfies \eqref{eq:qi}.
		Since \(q_i=\tau_i\), applying \eqref{S-eq:candidate_center_score}
		to the basis of \(\mathcal V\) gives
		\(S_{\mathcal V}^x(q_i)=\|a_i\|^2\).
		Now let \(t\) belong to the comparison set
		\eqref{eq:comparison set} for instance \(i\). For \(j<i\),
		we have \(\|t-\tau_j\|_\infty\ge B\), so the balls
		\(t+D_{B/2}\) and \(\tau_j+D_{B/2}\) intersect at most on a
		set of measure zero. Hence \(x_{j,t}=0\) in \(L^2(D_{B/2})\).
		For \(j\ge i\), we have \(\|t-\tau_j\|_\infty\ge\delta/2\).
		Define
		\[
		J_{i,t}:=\{j\in\{i,\ldots,N\}:x_{j,t}\ne0
		\text{ in }L^2(D_{B/2})\},
		\qquad v_{j,t}:=t-\tau_j.
		\]
		For \(j\in J_{i,t}\), the two balls intersect in positive
		measure, so \(\|v_{j,t}\|<B\). For distinct \(j,m\in J_{i,t}\),
		\eqref{eq:separation} gives
		\[
		\|v_{j,t}-v_{m,t}\|_\infty
		=\|\tau_j-\tau_m\|_\infty>B+\frac{3\delta}{2}.
		\]
		The vectors therefore satisfy the constraints in
		\eqref{S-eq:admissible_displacement_tuples}, and
		\(|J_{i,t}|\le 2^d\). If \(J_{i,t}\) is nonempty, it follows that
		\((v_{j,t})_{j\in J_{i,t}}\in\mathcal C_{|J_{i,t}|}(B,\delta)\).
		By \eqref{eq:Sx} and orthonormality, the clean score at \(t\)
		is the squared projection norm of the image in the ball
		\(t+D_{B/2}\). If \(J_{i,t}\) is empty, the score is zero.
		Otherwise, the triangle inequality and
		\eqref{S-eq:shifted_instance_projection_bound} give
		\[
		\begin{aligned}
			\sqrt{S_{\mathcal V}^x(t)}
			&=\left\|P_{\mathcal V}\sum_{j\in J_{i,t}}x_{j,t}
			\right\|_{L^2(D_{B/2})}\\
			&\le\sum_{j\in J_{i,t}}e_{\mathcal V}(\|v_{j,t}\|)\|a_j\|\\
			&\le\|a_i\|\sum_{j\in J_{i,t}}e_{\mathcal V}(\|v_{j,t}\|)
			\le\|a_i\|\chi_{\mathcal V}(\delta).
		\end{aligned}
		\]
		The last line uses the ordering \(\|a_j\|\le\|a_i\|\) for
		\(j\ge i\) and the definition of \(\chi_{\mathcal V}(\delta)\)
		in \eqref{S-eq:chi_V_def}. Thus, for every \(t\) in the comparison
		set of instance \(i\),
		\[
		S_{\mathcal V}^x(q_i)-S_{\mathcal V}^x(t)
		\ge \bigl(1-\chi_{\mathcal V}(\delta)^2\bigr)\|a_i\|^2.
		\]
		Hence the clean score map \(S_{\mathcal V}^x\) satisfies
			Assumption~\ref{def:localization_property} with
		\(\rho=1-\chi_{\mathcal V}(\delta)^2>0\). The constant is
		independent of \(L\), since the templates, \(B\), \(\mathcal V\),
		and \(\delta\) are fixed.
	\end{proof}
	
	\begin{remark}[The ambient space for \(\delta\) near \(2B\)]
		\label{S-rem:nonempty_retained_family}
		Suppose that \(\mathcal S\in\mathfrak S\). Since the templates
		have unit norm, \eqref{eq:uniform_template_approximation} with
		\(0<\varepsilon<1\) gives \(P_{\mathcal S}u_\ell\ne0\) for every
		\(\ell\). To show that \(\mathcal S\in\mathfrak S_{\rm loc}\) for
		\(\delta\) sufficiently close to \(2B\), it remains to verify
		\(\chi_{\mathcal S}(\delta)<1\).
		For each template \(u_\ell\), the translate \(u_\ell(h+\cdot)\)
		depends continuously on \(h\) in \(L^2(\R^d)\). Restricting two
		translates to \(D_{B/2}\) and projecting them onto \(\mathcal S\)
		does not increase their \(L^2\) distance. Since the denominator in
		\eqref{S-eq:radial_projection_envelope} is fixed and positive, the
		ratio for \(u_\ell\) is continuous in \(h\). With
		\(h=\lambda\omega\), the ratio is uniformly continuous on the
		compact set \([0,B]\times S^{d-1}\). Its maximum over \(\omega\)
		is therefore continuous in \(\lambda\). Taking the maximum over
		the finitely many templates gives continuity of \(e_{\mathcal S}\)
		on \([0,B]\).
		When \(\|h\|=B\), the ball \(D_{B/2}\) and the ball
		\(-h+D_{B/2}\) supporting \(u_\ell(h+\cdot)\) intersect only
		on a set of measure zero. Hence
		\begin{equation}\label{S-eq:radial_envelope_at_B}
			e_{\mathcal S}(B)=0.
		\end{equation}
		By \eqref{S-eq:admissible_displacement_tuples}, every vector in a
		tuple contributing to \eqref{S-eq:chi_V_def} has Euclidean norm
		between \(\delta/2\) and \(B\). Thus
		\[
		0\le\chi_{\mathcal S}(\delta)
		\le 2^d\sup_{\delta/2\le\lambda\le B}e_{\mathcal S}(\lambda).
		\]
		By continuity and \eqref{S-eq:radial_envelope_at_B}, the upper bound
		tends to zero as \(\delta\uparrow2B\) and equals zero at
		\(\delta=2B\). Hence there exists \(\delta_0\in(0,2B)\) such that
		\(\mathcal S\in\mathfrak S_{\rm loc}\) for every
		\(\delta\in[\delta_0,2B]\).
	\end{remark}
	
	\begin{lemma}[The noise assumption for candidate subspaces]
		\label{S-lem:noise_reg_inheritance}
		If Assumption~\ref{ass:noise_reg_main} holds for an
			orthonormal basis of \(\mathcal S\), then it holds for
		any fixed orthonormal basis of each
		\(\mathcal V\in\mathfrak S\).
	\end{lemma}
	
	\begin{proof}
		For every nonzero \(f\in\mathcal V\subseteq\mathcal S\),
			Assumption~\ref{ass:noise_reg_main} for the basis of
			\(\mathcal S\) gives
		\[
		\Var\!\left(\int_{D_{B/2}}z(s)f(s)\,\dd s\right)>0.
		\]
		Every nonzero linear combination of an orthonormal basis
			of \(\mathcal V\) is such a function \(f\). Thus the covariance
			matrix of the corresponding projected noise coordinates at
			\(0\) is positive definite.
	\end{proof}
	
	\begin{proof}[Proof of Corollary~\ref{cor:asymptotic_selected}]
		Every \(\mathcal V\in\mathfrak S_{\rm loc}\) is steerable by
		Lemma~\ref{S-lem:real_angular_expansion}.
		Proposition~\ref{prop:template_localization_condition}
		and Lemma~\ref{S-lem:noise_reg_inheritance} verify the localization
		and noise assumptions of
		Propositions~\ref{prop:asymptotic_error_power}
		and~\ref{prop:bh_density_strength} for \(\mathcal V\).
		By Lemma~\ref{S-lem:noise_maxima_fixed_growing},
		\(0<\eta_{\mathrm{peak}}(\mathcal V)<\infty\).
		By Lemma~\ref{S-lem:proj_rot_appendix} and
			\eqref{eq:aL_V_def},
		\[
			a_L(\mathcal V)^2
			=\min_{1\le i\le N}b_i^2
			\|P_{\mathcal V}u_{\ell_i}\|_{L^2(D_{B/2})}^2
			\ge b_L^2\min_{1\le\ell\le n}
			\|P_{\mathcal V}u_\ell\|_{L^2(D_{B/2})}^2.
		\]
		Hence, by~\eqref{eq:psnr_V},
		\begin{equation}\label{S-eq:psnr_retained_bound}
				\frac{a_L(\mathcal V)^2}{\eta_{\mathrm{peak}}(\mathcal V)}
				\ge b_L^2\mathrm{pSNR}(\mathcal V).
		\end{equation}
		Take positive constants
		\(A_{\rm BH}(\mathcal V),C(\mathcal V),c(\mathcal V)\) from
		\eqref{S-eq:bh_power_explicit_constants}, using the same absolute
		constant \(C_{\rm BH}\) for every subspace. Define
		\begin{equation}\label{S-eq:selected_uniform_constants}
			\begin{aligned}
				A_\star
				&:=\max_{\mathcal V\in\mathfrak S_{\rm loc}}
				\frac{A_{\rm BH}(\mathcal V)}
				{\eta_{\mathrm{peak}}(\mathcal V)},\\
				C_\star
				&:=\max_{\mathcal V\in\mathfrak S_{\rm loc}}C(\mathcal V),\\
				c_\star
				&:=\min_{\mathcal V\in\mathfrak S_{\rm loc}}c(\mathcal V).
			\end{aligned}
		\end{equation}
		Since \(\mathfrak S_{\rm loc}\) is finite, fixed, and nonempty, these
		constants are finite and positive, and are independent of \(L\),
		\(N\), and the unknown instance parameters.
		Fix \(\mathcal V\in\mathfrak S_{\rm loc}\) and
			\(0<\beta<1\), and suppose that
			\eqref{eq:selected_bh_target_power} holds.
			By~\eqref{S-eq:psnr_retained_bound},
		\[
		\frac{a_L(\mathcal V)^2}{\eta_{\mathrm{peak}}(\mathcal V)}
		\ge b_L^2\mathrm{pSNR}(\mathcal V)
		\ge A_\star
		\ge\frac{A_{\rm BH}(\mathcal V)}
		{\eta_{\mathrm{peak}}(\mathcal V)}.
		\]
		Thus \(a_L(\mathcal V)^2\ge A_{\rm BH}(\mathcal V)\).
		The same bound also gives the logarithmic condition in
		\eqref{eq:bh_finite_power_conditions}.
		Proposition~\ref{prop:bh_density_strength} and
		\eqref{S-eq:selected_uniform_constants} yield
		\[
		\begin{aligned}
			1-\Power(u_{\rm BH})
			&\le C(\mathcal V)\exp\!\left\{
			-c(\mathcal V)
			\frac{a_L(\mathcal V)^2}{\eta_{\mathrm{peak}}(\mathcal V)}
			\right\}\\
			&\le C_\star\exp\!\left\{
			-c_\star b_L^2\mathrm{pSNR}(\mathcal V)
			\right\}\le\beta,
		\end{aligned}
		\]
		where the last inequality follows from
			\eqref{eq:selected_bh_target_power}. Here \(u_{\rm BH}\) is
			the BH threshold for the procedure using \(\mathcal V\).
			This proves part~\textup{(i)}.
		For simultaneous recovery, \eqref{S-eq:psnr_retained_bound}
			and~\eqref{eq:selected_simultaneous_recovery} give
		\[
			\begin{aligned}
				\liminf_{L\to\infty}
				\frac{a_L(\mathcal V)^2}
				{\eta_{\mathrm{peak}}(\mathcal V)\log M_L}
				&\ge \mathrm{pSNR}(\mathcal V)
				\liminf_{L\to\infty}\frac{b_L^2}{\log M_L}\\
				&>\Gamma_\star\ge\Gamma_{\mathcal V},
		\end{aligned}
		\]
		since \(\mathrm{pSNR}(\mathcal V)>0\) is independent
			of \(L\) and
			\(\Gamma_\star=\max_{\mathcal V\in\mathfrak S_{\rm loc}}
			\Gamma_{\mathcal V}\).
		Corollary~\ref{cor:simultaneous_recovery_asymptotic} gives
		simultaneous recovery at both thresholds and the stated
		Bonferroni FWER bound.
		The right-hand sides of
			\eqref{eq:selected_bh_target_power}
			and~\eqref{eq:selected_simultaneous_recovery} depend on
			\(\mathcal V\) only through \(1/\mathrm{pSNR}(\mathcal V)\).
			The choice of \(\mathcal S^\star\) in~\eqref{eq:Sstar_def}
			therefore minimizes both right-hand sides.
		With BH power tending to one, Proposition~\ref{prop:bh_fdr_control}
		gives
		\[
		\limsup_{L\to\infty}\FDR(u_{\rm BH})\le\alpha
		\]
		under either of the two additional conditions stated in
		that proposition.
	\end{proof}

\bibliographystyle{imsart-number}
\bibliography{refs}

\begin{thebibliography}{30}

\bibitem{abrahamCastilloRoquain2024}
\begin{barticle}[author]
\bauthor{\bsnm{Abraham},~\bfnm{Kweku}\binits{K.}},
  \bauthor{\bsnm{Castillo},~\bfnm{Isma{\"e}l}\binits{I.}} \AND
  \bauthor{\bsnm{Roquain},~\bfnm{{\'E}tienne}\binits{{\'E}.}}
(\byear{2024}).
\btitle{Sharp Multiple Testing Boundary for Sparse Sequences}.
\bjournal{The Annals of Statistics}
\bvolume{52}
\bpages{1564--1591}.
\bdoi{10.1214/24-AOS2404}
\end{barticle}
\endbibitem

\bibitem{adlerSubagTaylor2012}
\begin{barticle}[author]
\bauthor{\bsnm{Adler},~\bfnm{Robert~J.}\binits{R.~J.}},
  \bauthor{\bsnm{Subag},~\bfnm{Eliran}\binits{E.}} \AND
  \bauthor{\bsnm{Taylor},~\bfnm{Jonathan~E.}\binits{J.~E.}}
(\byear{2012}).
\btitle{Rotation and Scale Space Random Fields and the Gaussian Kinematic
  Formula}.
\bjournal{The Annals of Statistics}
\bvolume{40}
\bpages{2910--2942}.
\bdoi{10.1214/12-AOS1055}
\end{barticle}
\endbibitem

\bibitem{adlerTaylor2007}
\begin{bbook}[author]
\bauthor{\bsnm{Adler},~\bfnm{Robert~J.}\binits{R.~J.}} \AND
  \bauthor{\bsnm{Taylor},~\bfnm{Jonathan~E.}\binits{J.~E.}}
(\byear{2007}).
\btitle{Random Fields and Geometry}.
\bpublisher{Springer}.
\end{bbook}
\endbibitem

\bibitem{almiraEtAl2024}
\begin{barticle}[author]
\bauthor{\bsnm{Almira},~\bfnm{Jos{\'e}~Mar{\'i}a}\binits{J.~M.}},
  \bauthor{\bsnm{Phelippeau},~\bfnm{Harold}\binits{H.}} \AND
  \bauthor{\bsnm{Martinez-Sanchez},~\bfnm{Antonio}\binits{A.}}
(\byear{2024}).
\btitle{Fast Normalized Cross-Correlation for Template Matching with
  Rotations}.
\bjournal{Journal of Applied Mathematics and Computing}
\bvolume{70}
\bpages{4937--4969}.
\bdoi{10.1007/s12190-024-02157-6}
\end{barticle}
\endbibitem

\bibitem{benjaminiHochberg1995}
\begin{barticle}[author]
\bauthor{\bsnm{Benjamini},~\bfnm{Yoav}\binits{Y.}} \AND
  \bauthor{\bsnm{Hochberg},~\bfnm{Yosef}\binits{Y.}}
(\byear{1995}).
\btitle{Controlling the false discovery rate: a practical and powerful approach
  to multiple testing}.
\bjournal{Journal of the Royal Statistical Society: Series B (Methodological)}
\bvolume{57}
\bpages{289--300}.
\bdoi{10.1111/j.2517-6161.1995.tb02031.x}
\end{barticle}
\endbibitem

\bibitem{bergelsonRosenblatt1988}
\begin{barticle}[author]
\bauthor{\bsnm{Bergelson},~\bfnm{Vitaly}\binits{V.}} \AND
  \bauthor{\bsnm{Rosenblatt},~\bfnm{Joseph}\binits{J.}}
(\byear{1988}).
\btitle{Mixing Actions of Groups}.
\bjournal{Illinois Journal of Mathematics}
\bvolume{32}
\bpages{65--80}.
\bdoi{10.1215/ijm/1255989229}
\end{barticle}
\endbibitem

\bibitem{carion2020detr}
\begin{binproceedings}[author]
\bauthor{\bsnm{Carion},~\bfnm{Nicolas}\binits{N.}},
  \bauthor{\bsnm{Massa},~\bfnm{Francisco}\binits{F.}},
  \bauthor{\bsnm{Synnaeve},~\bfnm{Gabriel}\binits{G.}},
  \bauthor{\bsnm{Usunier},~\bfnm{Nicolas}\binits{N.}},
  \bauthor{\bsnm{Kirillov},~\bfnm{Alexander}\binits{A.}} \AND
  \bauthor{\bsnm{Zagoruyko},~\bfnm{Sergey}\binits{S.}}
(\byear{2020}).
\btitle{End-to-End Object Detection with Transformers}.
In \bbooktitle{Computer Vision -- {ECCV} 2020}.
\bseries{Lecture Notes in Computer Science}
\bvolume{12346}
\bpages{213--229}.
\bpublisher{Springer}.
\bdoi{10.1007/978-3-030-58452-8_13}
\end{binproceedings}
\endbibitem

\bibitem{chengSchwartzman2017}
\begin{barticle}[author]
\bauthor{\bsnm{Cheng},~\bfnm{Dan}\binits{D.}} \AND
  \bauthor{\bsnm{Schwartzman},~\bfnm{Armin}\binits{A.}}
(\byear{2017}).
\btitle{Multiple Testing of Local Maxima for Detection of Peaks in Random
  Fields}.
\bjournal{The Annals of Statistics}
\bvolume{45}
\bpages{529--556}.
\bdoi{10.1214/16-AOS1458}
\end{barticle}
\endbibitem

\bibitem{daiXu2013spheres}
\begin{bbook}[author]
\bauthor{\bsnm{Dai},~\bfnm{Feng}\binits{F.}} \AND
  \bauthor{\bsnm{Xu},~\bfnm{Yuan}\binits{Y.}}
(\byear{2013}).
\btitle{Approximation Theory and Harmonic Analysis on Spheres and Balls}.
\bseries{Springer Monographs in Mathematics}.
\bpublisher{Springer}.
\bdoi{10.1007/978-1-4614-6660-4}
\end{bbook}
\endbibitem

\bibitem{fageot2021steerability}
\begin{barticle}[author]
\bauthor{\bsnm{Fageot},~\bfnm{Julien}\binits{J.}},
  \bauthor{\bsnm{Uhlmann},~\bfnm{Virginie}\binits{V.}},
  \bauthor{\bsnm{P{\"u}sp{\"o}ki},~\bfnm{Zsuzsanna}\binits{Z.}},
  \bauthor{\bsnm{Beck},~\bfnm{Benjamin}\binits{B.}},
  \bauthor{\bsnm{Unser},~\bfnm{Michael}\binits{M.}} \AND
  \bauthor{\bsnm{Depeursinge},~\bfnm{Adrien}\binits{A.}}
(\byear{2021}).
\btitle{Principled Design and Implementation of Steerable Detectors}.
\bjournal{IEEE Transactions on Image Processing}
\bvolume{30}
\bpages{4465--4478}.
\bdoi{10.1109/TIP.2021.3072499}
\end{barticle}
\endbibitem

\bibitem{freemanAdelson1991}
\begin{barticle}[author]
\bauthor{\bsnm{Freeman},~\bfnm{William~T.}\binits{W.~T.}} \AND
  \bauthor{\bsnm{Adelson},~\bfnm{Edward~H.}\binits{E.~H.}}
(\byear{1991}).
\btitle{The Design and Use of Steerable Filters}.
\bjournal{IEEE Transactions on Pattern Analysis and Machine Intelligence}
\bvolume{13}
\bpages{891--906}.
\bdoi{10.1109/34.93808}
\end{barticle}
\endbibitem

\bibitem{hoffmannJorgensenSheppDudley1979}
\begin{barticle}[author]
\bauthor{\bsnm{Hoffmann-J{\o}rgensen},~\bfnm{J.}\binits{J.}},
  \bauthor{\bsnm{Shepp},~\bfnm{L.~A.}\binits{L.~A.}} \AND
  \bauthor{\bsnm{Dudley},~\bfnm{R.~M.}\binits{R.~M.}}
(\byear{1979}).
\btitle{On the Lower Tail of Gaussian Seminorms}.
\bjournal{The Annals of Probability}
\bvolume{7}
\bpages{319--342}.
\bdoi{10.1214/aop/1176995091}
\end{barticle}
\endbibitem

\bibitem{holm1979}
\begin{barticle}[author]
\bauthor{\bsnm{Holm},~\bfnm{Sture}\binits{S.}}
(\byear{1979}).
\btitle{A simple sequentially rejective multiple test procedure}.
\bjournal{Scandinavian Journal of Statistics}
\bvolume{6}
\bpages{65--70}.
\end{barticle}
\endbibitem

\bibitem{kallenberg2021}
\begin{bbook}[author]
\bauthor{\bsnm{Kallenberg},~\bfnm{Olav}\binits{O.}}
(\byear{2021}).
\btitle{Foundations of Modern Probability},
\bedition{3rd} ed.
\bseries{Probability Theory and Stochastic Modelling}
\bvolume{99}.
\bpublisher{Springer}, \baddress{Cham}.
\bdoi{10.1007/978-3-030-61871-1}
\end{bbook}
\endbibitem

\bibitem{kay1998detection}
\begin{bbook}[author]
\bauthor{\bsnm{Kay},~\bfnm{Steven~M.}\binits{S.~M.}}
(\byear{1998}).
\btitle{Fundamentals of Statistical Signal Processing: Detection Theory}.
\bseries{Fundamentals of Statistical Signal Processing}
\bvolume{2}.
\bpublisher{Prentice Hall}.
\end{bbook}
\endbibitem

\bibitem{merhav2011threshold}
\begin{barticle}[author]
\bauthor{\bsnm{Merhav},~\bfnm{Neri}\binits{N.}}
(\byear{2011}).
\btitle{Threshold effects in parameter estimation as phase transitions in
  statistical mechanics}.
\bjournal{IEEE Transactions on Information Theory}
\bvolume{57}
\bpages{7000--7010}.
\bdoi{10.1109/TIT.2011.2165824}
\end{barticle}
\endbibitem

\bibitem{wakninEtAl2026object}
\begin{barticle}[author]
\bauthor{\bsnm{Mor~Waknin},~\bfnm{Keren}\binits{K.}},
  \bauthor{\bsnm{Eldar},~\bfnm{Amitay}\binits{A.}},
  \bauthor{\bsnm{Davenport},~\bfnm{Samuel}\binits{S.}},
  \bauthor{\bsnm{Bendory},~\bfnm{Tamir}\binits{T.}},
  \bauthor{\bsnm{Schwartzman},~\bfnm{Armin}\binits{A.}} \AND
  \bauthor{\bsnm{Shkolnisky},~\bfnm{Yoel}\binits{Y.}}
(\byear{2026}).
\btitle{Object Detection under the Linear Subspace Model with Application to
  {Cryo-EM} Images}.
\bjournal{The Annals of Statistics}
\bvolume{54}
\bpages{910--928}.
\bdoi{10.1214/25-AOS2595}
\end{barticle}
\endbibitem

\bibitem{passeggeriVeraart2019}
\begin{barticle}[author]
\bauthor{\bsnm{Passeggeri},~\bfnm{Riccardo}\binits{R.}} \AND
  \bauthor{\bsnm{Veraart},~\bfnm{Almut E.~D.}\binits{A.~E.~D.}}
(\byear{2019}).
\btitle{Mixing Properties of Multivariate Infinitely Divisible Random Fields}.
\bjournal{Journal of Theoretical Probability}
\bvolume{32}
\bpages{1845--1879}.
\bdoi{10.1007/s10959-018-0864-7}
\end{barticle}
\endbibitem

\bibitem{piterbarg1996}
\begin{bbook}[author]
\bauthor{\bsnm{Piterbarg},~\bfnm{Vladimir~I.}\binits{V.~I.}}
(\byear{1996}).
\btitle{Asymptotic Methods in the Theory of Gaussian Processes and Fields}.
\bseries{Translations of Mathematical Monographs}
\bvolume{148}.
\bpublisher{American Mathematical Society}.
\end{bbook}
\endbibitem

\bibitem{puspokiEtAl2019}
\begin{barticle}[author]
\bauthor{\bsnm{P{\"u}sp{\"o}ki},~\bfnm{Zsuzsanna}\binits{Z.}},
  \bauthor{\bsnm{Fageot},~\bfnm{Julien}\binits{J.}},
  \bauthor{\bsnm{Amini},~\bfnm{Arash}\binits{A.}},
  \bauthor{\bsnm{Ward},~\bfnm{John~Paul}\binits{J.~P.}} \AND
  \bauthor{\bsnm{Unser},~\bfnm{Michael}\binits{M.}}
(\byear{2019}).
\btitle{Angular Accuracy of Steerable Feature Detectors}.
\bjournal{SIAM Journal on Imaging Sciences}
\bvolume{12}
\bpages{344--371}.
\bdoi{10.1137/18M1190690}
\end{barticle}
\endbibitem

\bibitem{redmon2016yolo}
\begin{binproceedings}[author]
\bauthor{\bsnm{Redmon},~\bfnm{Joseph}\binits{J.}},
  \bauthor{\bsnm{Divvala},~\bfnm{Santosh}\binits{S.}},
  \bauthor{\bsnm{Girshick},~\bfnm{Ross}\binits{R.}} \AND
  \bauthor{\bsnm{Farhadi},~\bfnm{Ali}\binits{A.}}
(\byear{2016}).
\btitle{You Only Look Once: Unified, Real-Time Object Detection}.
In \bbooktitle{Proceedings of the IEEE Conference on Computer Vision and
  Pattern Recognition}
\bpages{779--788}.
\bdoi{10.1109/CVPR.2016.91}
\end{binproceedings}
\endbibitem

\bibitem{ren2015faster}
\begin{binproceedings}[author]
\bauthor{\bsnm{Ren},~\bfnm{Shaoqing}\binits{S.}},
  \bauthor{\bsnm{He},~\bfnm{Kaiming}\binits{K.}},
  \bauthor{\bsnm{Girshick},~\bfnm{Ross}\binits{R.}} \AND
  \bauthor{\bsnm{Sun},~\bfnm{Jian}\binits{J.}}
(\byear{2015}).
\btitle{Faster {R-CNN}: Towards Real-Time Object Detection with Region Proposal
  Networks}.
In \bbooktitle{Advances in Neural Information Processing Systems}
\bvolume{28}
\bpages{91--99}.
\end{binproceedings}
\endbibitem

\bibitem{scharfFriedlander1994}
\begin{barticle}[author]
\bauthor{\bsnm{Scharf},~\bfnm{Louis~L.}\binits{L.~L.}} \AND
  \bauthor{\bsnm{Friedlander},~\bfnm{Benjamin}\binits{B.}}
(\byear{1994}).
\btitle{Matched subspace detectors}.
\bjournal{IEEE Transactions on Signal Processing}
\bvolume{42}
\bpages{2146--2157}.
\bdoi{10.1109/78.301849}
\end{barticle}
\endbibitem

\bibitem{schwartzmanGavrilovAdler2011}
\begin{barticle}[author]
\bauthor{\bsnm{Schwartzman},~\bfnm{Armin}\binits{A.}},
  \bauthor{\bsnm{Gavrilov},~\bfnm{Yulia}\binits{Y.}} \AND
  \bauthor{\bsnm{Adler},~\bfnm{Robert~J.}\binits{R.~J.}}
(\byear{2011}).
\btitle{Multiple Testing of Local Maxima for Detection of Peaks in {1D}}.
\bjournal{The Annals of Statistics}
\bvolume{39}
\bpages{3290--3319}.
\bdoi{10.1214/11-AOS943}
\end{barticle}
\endbibitem

\bibitem{shafieSigalSiegmundWorsley2003}
\begin{barticle}[author]
\bauthor{\bsnm{Shafie},~\bfnm{Khalil}\binits{K.}},
  \bauthor{\bsnm{Sigal},~\bfnm{Bronislava}\binits{B.}},
  \bauthor{\bsnm{Siegmund},~\bfnm{David}\binits{D.}} \AND
  \bauthor{\bsnm{Worsley},~\bfnm{Keith~J.}\binits{K.~J.}}
(\byear{2003}).
\btitle{Rotation Space Random Fields with an Application to fMRI Data}.
\bjournal{The Annals of Statistics}
\bvolume{31}
\bpages{1732--1771}.
\bdoi{10.1214/aos/1074290326}
\end{barticle}
\endbibitem

\bibitem{siegmundWorsley1995}
\begin{barticle}[author]
\bauthor{\bsnm{Siegmund},~\bfnm{David~O.}\binits{D.~O.}} \AND
  \bauthor{\bsnm{Worsley},~\bfnm{Keith~J.}\binits{K.~J.}}
(\byear{1995}).
\btitle{Testing for a Signal with Unknown Location and Scale in a Stationary
  Gaussian Random Field}.
\bjournal{The Annals of Statistics}
\bvolume{23}
\bpages{608--639}.
\bdoi{10.1214/aos/1176324539}
\end{barticle}
\endbibitem

\bibitem{vanTreesBellTian2013}
\begin{bbook}[author]
\bauthor{\bsnm{Van~Trees},~\bfnm{Harry~L.}\binits{H.~L.}},
  \bauthor{\bsnm{Bell},~\bfnm{Kristine~L.}\binits{K.~L.}} \AND
  \bauthor{\bsnm{Tian},~\bfnm{Zhi}\binits{Z.}}
(\byear{2013}).
\btitle{Detection, Estimation, and Modulation Theory, Part I: Detection,
  Estimation, and Filtering Theory}.
\bpublisher{Wiley}.
\end{bbook}
\endbibitem

\bibitem{zakaiZiv1969threshold}
\begin{barticle}[author]
\bauthor{\bsnm{Zakai},~\bfnm{Moshe}\binits{M.}} \AND
  \bauthor{\bsnm{Ziv},~\bfnm{Jacob}\binits{J.}}
(\byear{1969}).
\btitle{On the threshold effect in radar range estimation}.
\bjournal{IEEE Transactions on Information Theory}
\bvolume{15}
\bpages{167--170}.
\bdoi{10.1109/TIT.1969.1054247}
\end{barticle}
\endbibitem

\bibitem{zhaoBlu2020}
\begin{barticle}[author]
\bauthor{\bsnm{Zhao},~\bfnm{Tianle}\binits{T.}} \AND
  \bauthor{\bsnm{Blu},~\bfnm{Thierry}\binits{T.}}
(\byear{2020}).
\btitle{The {Fourier-Argand} Representation: An Optimal Basis of Steerable
  Patterns}.
\bjournal{IEEE Transactions on Image Processing}
\bvolume{29}
\bpages{6357--6371}.
\bdoi{10.1109/TIP.2020.2990483}
\end{barticle}
\endbibitem

\bibitem{zivZakai1969}
\begin{barticle}[author]
\bauthor{\bsnm{Ziv},~\bfnm{Jacob}\binits{J.}} \AND
  \bauthor{\bsnm{Zakai},~\bfnm{Moshe}\binits{M.}}
(\byear{1969}).
\btitle{Some lower bounds on signal parameter estimation}.
\bjournal{IEEE Transactions on Information Theory}
\bvolume{15}
\bpages{386--391}.
\bdoi{10.1109/TIT.1969.1054301}
\end{barticle}
\endbibitem

\end{thebibliography}
\end{document}